\magnification=\magstep1
\input amstex
\documentstyle{amsppt}

\define\defeq{\overset{\text{def}}\to=}
\define\ab{\operatorname{ab}}
\define\pr{\operatorname{pr}}
\define\Gal{\operatorname{Gal}}
\define\diag{\operatorname{diag}}
\define\Hom{\operatorname{Hom}}
\define\sep{\operatorname{sep}}
\def \isom {\overset \sim \to \rightarrow}
\define\cor{\operatorname{cor}}
\define\Pic{\operatorname{Pic}}
\define\Br{\operatorname{Br}}
\define\norm{\operatorname{norm}}
\define\tame{\operatorname{tame}}
\define\period{\operatorname{period}}
\define\index{\operatorname{index}}
\define\Spec{\operatorname{Spec}}
\define\id{\operatorname{id}}
\define\cn{\operatorname{cn}}
\define\inn{\operatorname{inn}}
\define\Ker{\operatorname{Ker}}
\define\ob{\operatorname{ob}}
\def \c{\operatorname {c}}
\def \Aut{\operatorname {Aut}}
\def \Out{\operatorname {Out}}

\def \res{\operatorname {res}}
\def \Jac{\operatorname {Jac}}
\def \Ind{\operatorname {Ind}}
\def \out{\operatorname {out}}
\def \et{\operatorname {et}}
\def \Sec{\operatorname {Sec}}
\def\Sp{\operatorname {Sp}}

\def\Aut{\operatorname{Aut}}
\def\Out{\operatorname{Out}}
\def\Inn{\operatorname{Inn}}
\def\Im{\operatorname{Im}}
\def\char{\operatorname{char}}
\def\GASC{\operatorname{GASC}}
\def\Conjecture{\operatorname{Conjecture}}
\def\adic{\operatorname{adic}}

\define\Primes{\frak{Primes}}
\NoRunningHeads
\NoBlackBoxes
\topmatter

\title
Good Sections of Arithmetic Fundamental Groups
\endtitle

\author
Mohamed Sa\"\i di
\endauthor

\abstract
In this paper we exhibit the notion of (uniformly) good sections of arithmetic fundamental groups.
We introduce and investigate the problem of cuspidalisation of sections of arithmetic fundamental groups,
its ultimate aim is to reduce the solution of the Grothendieck anabelian section conjecture to the solution of its birational version. 
We show that (uniformly) good sections of arithmetic fundamental groups of smooth, proper, and geometrically
connected hyperbolic curves over slim (and regular) fields can be lifted to sections of cuspidally abelian absolute
Galois groups. As an application we prove a (pro-$p$) 
$p$-adic version of the Grothendieck anabelian section conjecture
for hyperbolic curves, under the assumption that the existence of sections of arithmetic fundamental groups,
and cuspidally abelian Galois groups, implies the existence
of tame points. We also prove that the existence of uniformly good sections of arithmetic fundamental groups for hyperbolic curves 
over number fields implies the existence of divisors of degree $1$, under a finiteness condition of the Tate-Shafarevich group of the 
jacobian of the curve.
\endabstract
\toc

\subhead
\S 0. Introduction
\endsubhead

\subhead
\S 1. Good Sections of Arithmetic Fundamental Groups
\endsubhead

\subhead
\S 2. Cuspidalisation of Good Sections of Arithmetic Fundamental Groups over Slim Fields
\endsubhead

\subhead
\S 3. Applications to the Grothendieck Anabelian Section Conjecture
\endsubhead

\subhead
\S 4. Sections of Geometrically pro-$\Sigma$ Arithmetic Fundamental Groups of Curves over
$p$-adic Local Field: $p\notin \Sigma$
\endsubhead

\subhead
\S 5. Appendix
\endsubhead

\endtoc

\endtopmatter

\document

\subhead
\S 0. Introduction
\endsubhead
Let $k$ be a field, and $X$ a proper, smooth, geometrically connected, and hyperbolic algebraic curve over
$k$. 

Let $\pi_1(X)$ be the arithmetic \'etale fundamental group of $X$, which sits in the following exact sequence
$$1\to \pi_1(\overline X)\to \pi_1(X) @>{\pr}>> G_k\to 1,$$
where $G_k$ is the absolute Galois group of $k$, and $\pi_1(\overline X)$ is the geometric \'etale fundamental group of $X$.

In this paper we investigate continuous group-theoretic sections $s:G_k\to \pi_1(X)$ of the natural projection
$\pr:\pi_1(X)\twoheadrightarrow G_k$, meaning that $\pr\circ s=\id_{G_k}$, or equivalently splittings of the above exact sequence, which
we will refer to as sections of the \'etale fundamental group $\pi_1(X)$. 

A section of $\pi_1(X)$ is uniquely determined by a closed subgroup of $\pi_1(X)$, namely the image of the section, which maps isomorphically to $G_k$ via the projection
$\pr:\pi(X)\twoheadrightarrow G_k$. 

Sections of $\pi_1(X)$ arise naturally from rational points of $X$. 

More precisely, a closed point $x\in X$ determines a decomposition subgroup
$D_x\subset \pi_1(X)$, which is defined only modulo conjugation by the elements of $\pi_1(\overline X)$, and which maps isomorphically to the open subgroup
$G_{k(x)}$ of $G_k$ via the projection $\pr:\pi_1(X) \twoheadrightarrow G_k$, where $k(x)$ is the residue field at $x$.

Thus, the decomposition group $D_x$ associated to a rational point $x\in X(k)$ determines a group-theoretic
section $s_x:G_k\to \pi _1(X)$ 
of $\pi _1(X)$, which is only defined up to conjugation by the elements of $\pi _1(\overline X)$.
We will refer to such a section of $\pi_1(X)$ as point-theoretic.

We have a natural set-theoretic map
$$\varphi_X: X(k)\to  \overline {\Sec} _{\pi _1(X)},$$
$$x\mapsto \varphi_X (x)=[s_x],$$
where $\overline {\Sec} _{\pi_1(X)}$ is the set of conjugacy classes of all continuous group-theoretic sections 
of $\pi _1(X)$, modulo inner conjugation by the elements of $\pi _1(\overline X)$,
and $[s_x]$ denotes the image, i.e. conjugacy class, of a section $s_x$ associated to the rational point $x$ in $\overline {\Sec} _{\pi _1(X)}$.

In his seminal letter to Faltings, Grothendieck formulated the following fundamental conjecture (cf. [Grothendieck]).

\definition {Grothendieck's Anabelian Section Conjecture (GASC)}  Assume that $k$ is finitely generated over the prime field
$\Bbb Q$. Then the map $\varphi_X: X(k)\to  \overline {\Sec} _{\pi _1(X)}$ is bijective.
\enddefinition

In investigating the GASC, one
is naturally led to formulate, and investigate, an analogous conjecture over $p$-adic local fields (finite extension of $\Bbb Q_p$),
which we will refer to as the $p$-adic GASC (cf. 3.1. for more details).

The injectivity of the map $\varphi_X$, if $k$ is finitely generated over the prime field $\Bbb Q$, or a $p$-adic local field, is well-known. 
So the statements of the GASC, and $p$-adic GASC, are equivalent to the
surjectivity of the map $\varphi_X$, i.e. that every group-theoretic section of $\pi_1(X)$ is a section associated to a rational point, under the above assumptions 
on the field $k$.

One can also consider sections of absolute Galois groups of function fields of curves.
 
Let $K_X$ be the function field of $X$, and $K_X^{\sep}$ a separable closure of $K_X$. Let 
$\overline G_X\defeq \Gal (K_X^{\sep}/K_X.\bar k)$, and
$G_X\defeq \Gal (K_X^{\sep}/K_X)$.
Thus, $G_X$ sits naturally
in the following exact sequence
$$1\to \overline G_X\to G_X\to G_k\to 1.$$

Similarly, as above, let $x\in X(k)$ be a rational point. Then $x$ determines a decomposition subgroup
$D_x\subset G_X$, which is only defined up to conjugation by the elements of $\overline G_X$, and which maps surjectively
to $G_k$ via the natural projection $G_X\twoheadrightarrow G_k$. More precisely, $D_x$ sits naturally in the following exact sequence
$$1\to \hat \Bbb Z(1)\to D_x\to G_k\to 1.$$

The above group extension is known to be split. 
Each section $G_k\to D_x$ of the natural projection $D_x\twoheadrightarrow G_k$
determines naturally a section $G_k\to G_X$ of the natural projection $G_X\twoheadrightarrow G_k$, whose image is contained in $D_x$.

In light of the GASC, it is natural to formulate the following birational version.

\definition {The Birational Grothendieck Anabelian Section Conjecture (BGASC)} 
\newline
Assume that $k$ is finitely generated over the
prime field $\Bbb Q$. Let $s:G_k\to G_X$ be a group-theoretic section of the natural projection
$G_X \twoheadrightarrow G_k$. Then the image $s(G_k)$ is contained in a decomposition subgroup $D_x$ associated to a unique
rational point $x\in X(k)$. In particular, the existence of the section $s$ implies that $X(k)\neq \varnothing$.
\enddefinition

One can also formulate an analog of the BGASC over $p$-adic local fields, which we will refer to as the $p$-adic BGASC.

Special examples of the GASC were investigated in [Harari-Szamuely], and [Stix]. See also  [Esnault-Wittenberg1], and [Harari-Stix], 
for some partial results concerning the BGASC.

The GASC, and its $p$-adic version, are still widely open.

An important progress, in the last ten years, around the  birational version of the Grothendieck anabelian section conjecture is the proof by 
Koenigsmann, refined by Pop, that the $p$-adic BGASC holds true (cf. [Koenigsmann], and [Pop]).

The main motivation of this paper is to establish a theory for group-theoretic sections of $\pi_1(X)$, whose ultimate aim is to reduce the solution of the GASC 
(resp. $p$-adic GASC) to the solution of its birational version, the birational GASC (resp. $p$-adic BGASC). 
In this paper we introduce, and investigate, the theory of cuspidalisation of sections of $\pi_1(X)$ for this purpose.

In $\S1$ we exhibit a necessary condition for a group-theoretic section $s:G_k\to \pi_1(X)$
of $\pi_1(X)$ to be point-theoretic, which we call the (uniform) goodness condition (cf. Definition 1.4.1).  

Goodness of the section $s:G_k\to \pi_1(X)$ 
means that the natural pull back homomorphism of cohomology classes
$$s^{\star}:H_{\et}^2(X,\hat \Bbb Z(1))\to H_{\et}^2(G_k,\hat \Bbb Z(1)),$$ 
which is naturally induced by the section $s$, annihilates the Chern classes of line bundles. 
A similar condition should hold for every neighbourhood of the section, and after finite extensions of the base field $k$ (cf. loc. cit.). 

The goodness condition is equivalent, in the case where $k$ is a $p$-adic local field, to the existence of tame points, 
i.e. $X(k^{\tame})\neq \varnothing$, where $k^{\tame}$ is the maximal tamely ramified extension of $k$ (cf. Proposition 1.6.6).

If $k$ is a number field, the section $s$ gives rise naturally to sections
$s_v:G_{k_v}\to \pi_1(X_{k_v})$ of $\pi_1(X_{k_v})$, for every place $v$ of $k$, where $k_v$ is the completion of $k$ at $v$, and
$X_v\defeq X\times _k k_v$. 

A local-global principle for
goodness holds. Namely the section $s$ is (uniformly) good if and only if the sections $s_v$ are (uniformly) good for all places $v$ of $k$
(cf. Proposition 1.8.1).

In $\S2$ we introduce, and investigate, the problem of cuspidalisation for sections of arithmetic fundamental groups, 
which is formulated as follows.

\definition {The Cuspidalisation Problem for Sections of $\pi_1(X)$}
Given a group-theoretic section $s:G_k\to \pi_1(X)$
of the natural projection $\pr:\pi_1(X)\twoheadrightarrow G_k$, is it possible to lift it to a section  $\tilde s:G_k\to G_{K_X}$
of the natural projection $G_{K_X}\twoheadrightarrow G_k$? i.e. is it possible to construct a section $\tilde s$ such that
the following diagram is commutative
$$
\CD
G_k @>{\tilde s}>> G_{X} \\
@V{\id}VV   @VVV \\
G_k @>{s}>>  \pi _1(X)
\endCD
$$
where the right vertical map is the natural (surjective) homomorphism $G_{X}\twoheadrightarrow \pi_1(X)$?
\enddefinition

Note that the cuspidalisation problem has a positive solution if the section $s$ is point-theoretic.

A positive solution to both the BGASC (resp. $p$-adic BGASC), and the cuspidalisation problem,
in the case where $k$ is finitely generated over the prime field (resp. a $p$-adic local field), would imply a positive solution to the GASC (resp. $p$-adic GASC)
(cf. Remark 3.2.1).

Our main result in this paper, concerning the cuspidalisation problem, is the following (cf. Theorem 2.4, Theorem 2.6, and Theorem 2.7).

\proclaim {Theorem A} Assume that the field $k$ is a number field, or a $p$-adic local field.
More generally, assume that $k$ is slim, and regular (cf. Definition 2.3.1). 
Then a section $s:G_k\to \pi_1(X)$ of $\pi_1(X)$ is
uniformly good if and only if it can be lifted to a section  
$$s^{\c-\ab}:G_k\to G_{X}^{\c-\ab}$$ 
of the natural projection
$G_{X}^{\c-\ab}\twoheadrightarrow G_k$, where $G_{X}^{\c-\ab}$ is the maximal cuspidally abelian
quotient of $G_{X}$,  i.e. if and only if one can construct a section $s^{\c-\ab}$ which
inserts into the following commutative diagram:
$$
\CD
G_k @>s^{\c-\ab}>>   G_{X}^{\c-\ab}\\
@V{\id}VV     @VVV  \\
G_k   @>{s}>> \Pi_X
\endCD
$$

Here, $G_{X}^{\c-\ab}$ is the maximal quotient $H$ of $G_{X}$, satisfying $G_X\twoheadrightarrow  H \twoheadrightarrow \pi_1(X)$, and 
such that the kernel of the natural homomorphism $H\twoheadrightarrow \pi_1(X)$, which is generated by inertia subgroups at the geometric points of $X$,
is abelian (cf. 2.1.2 for more details).
\endproclaim

Moreover, we prove a pro-$\Sigma$ version of Theorem A, where $\Sigma$ is a non-empty set of prime integers (cf. loc. cit.)

In $\S3$ we apply our results on the cuspidalisation problem to the Grothendieck anabelian section conjecture.

As an application of Theorem A, we prove a (pro-$p$) $p$-adic version of the Grothendieck anabelian section conjecture,
under the assumption that the existence of sections of arithmetic fundamental groups, and sections of 
cuspidally abelian absolute Galois groups, 
implies the existence of tame points, by reducing to the birational version of the Grothendieck anabelian
section conjecture over $p$-adic local fields that was proven by Pop in [Pop].

More precisely, we prove the following (cf. Theorem 3.3.3, and Corollary 3.3.5).

\proclaim {Theorem B} Assume that $k$ is a $p$-adic local field.  
Consider the following two properties.

(i)\ For every proper, smooth, geometrically connected, and 
hyperbolic curve over $k$.
Every group-theoretic section $s:G_k\to \pi_1(X)$ of $\pi_1(X)$ is good.

(ii)\ For every proper, smooth, geometrically connected, and 
hyperbolic curve over $k$.
Every group-theoretic section $\tilde s :G_k\to G_X^{\c-\ab}$ of the natural projection
$G_X^{\c-\ab}\twoheadrightarrow G_k$ is tame-point theoretic (cf. Definition 1.7.1), meaning the following.
 Let $\Tilde L/K_X$ be the sub-extension of $K_X^{\sep}/K_X$ with Galois group $G_X^{\c-\ab}$, and $L/K$ 
 the sub-extension $\Tilde L/K_X$ corresponding to the closed subgroup $\tilde s(G_k)$.
Then the $p$-primary part of the kernel of the natural map $\Br (k)\to \Br(L)$ between Brauer groups is trivial.

Then we have the following implication:

$$ \{(i)+(ii)\}\Longrightarrow \{p-{\adic}\ \  \GASC\}.$$
\endproclaim

Over number fields we prove that the existence of uniformly good sections of arithmetic fundamental groups
implies the existence of divisors of degree $1$, under the condition that the Tate-Shafarevich group of the jacobian of the curve is finite,
by reducing to an analogous birational result proven in [Esnault-Wittenberg1]. 

More precisely, we prove the following (cf. Theorem 3.4.1).

\proclaim {Theorem C} Assume that $k$ is a number field, i.e. a finite extension of the prime field $\Bbb Q$. 
Let $s:G_k\to \pi_1(X)$ be a group theoretic section of the natural projection  $\pi_1(X) \twoheadrightarrow G_k$. 
Assume that $s$ is a uniformly good group-theoretic section, and that the jacobian variety
of $X$ has a finite Tate-Shafarevich group. Then there exists a divisor of degree $1$ on $X$.
\endproclaim

In $\S4$ we investigate sections of geometrically pro-$\Sigma$ arithmetic fundamental groups of hyperbolic algebraic
curves over $p$-adic local field, in the case where $p\notin \Sigma$. We give
examples of such sections which are not point-theoretic (cf. Proposition 4.2.1), despite the fact that such sections are good sections (cf. Proposition 4.3.1).

\definition{Acknowledgment} This paper grew up from several lengthy discussions I had with Akio Tamagawa
during my visit to the Research Institute for Mathematical Sciences (RIMS) of Kyoto university during the summer 2008.
I would like to express my gratitude to the research staff of RIMS for inviting me, and for the wonderful working atmosphere.
I would like very much to thank Akio Tamagawa for  the very fruitful discussions we had, and for sharing with me his 
knowledge on the section conjecture, and beyond. He especially helped with several technical points in this paper.
\enddefinition

\subhead
\S 1. Good Sections of Arithmetic Fundamental Groups
\endsubhead

In $\S1$ we introduce the notion of (uniformly) good sections of arithmetic fundamental groups, and we investigate some of their
properties over $p$-adic local base fields, and number fields.

\subhead {1.1}
\endsubhead
In this subsection we recall some general facts on arithmetic fundamental groups, and their group-theoretic sections.
We also fix some notations, that will be used throughout this paper.

Let $k$ be a field of characteristic $\char (k)=l\ge 0$. Let $X$ be a proper, smooth, geometrically connected,
hyperbolic algebraic curve over $k$, and $K\defeq K_X$ the function field of $X$.

Let $\eta $ be a geometric point of $X$ above
the generic point of $X$. Then $\eta$ determines naturally an algebraic closure $\bar k$
of $k$, a separable closure $K_X^{\sep}$ of $K_X$,
and a geometric point $\bar \eta$ of $\overline X\defeq X\times _k \bar k$. 

There exists
a canonical exact sequence of profinite groups

$$1\to \pi_1(\overline X,\bar \eta)\to \pi_1(X, \eta) @>{\pr_X}>> G_k\to 1.\tag {$1.1$}$$
Here, $\pi_1(X, \eta)$ denotes the arithmetic \'etale fundamental group of $X$ with base
point $\eta$, $\pi_1(\overline X,\bar \eta)$ the \'etale fundamental group of $\overline X$ with base
point $\bar \eta$, and $G_k\defeq \Gal (\bar k/k)$ the absolute Galois group of $k$.

Throughout this paper $\Primes$ denotes the set of all prime integers. 

We will consider the following variant of the exact sequence (1.1). 
Let 
$$\Sigma \subseteq \Primes$$
be a non-empty set of prime integers.
In the case where $\char (k)=l >0$, we will assume that $l\notin \Sigma$. 

Write
$$\Delta_X\defeq \pi_1(\overline X,\bar \eta)^{\Sigma}$$
for the maximal pro-$\Sigma$ quotient of $\pi_1(\overline X,\bar \eta)$,
and
$$\Pi_X\defeq  \pi_1(X, \eta)/ \Ker  (\pi_1(\overline X,\bar \eta)\twoheadrightarrow
\pi_1(\overline X,\bar \eta)^{\Sigma})$$
for the quotient of  $\pi_1(X, \eta)$ by the kernel of the natural surjective homomorphism
$\pi_1(\overline X,\bar \eta)\twoheadrightarrow \pi_1(\overline X,\bar \eta)^{\Sigma}$, which is a normal subgroup
of $\pi _1(X,\eta)$. Thus, we have a natural exact sequence of profinite groups
$$1\to \Delta_X\to \Pi_X @>{\pr_{X,\Sigma}}>> G_k\to 1.\tag {$1.2$}$$

We shall refer to 
$$\pi_1(X, \eta)^{(\Sigma)}\defeq \Pi_X$$
as the (maximal)  geometrically pro-$\Sigma$ quotient of $\pi_1(X, \eta)$, or the geometrically pro-$\Sigma$ 
arithmetic fundamental group of $X$. 

Let
$$s:G_k\to \Pi_X$$
be a continuous group-theoretic section of the natural projection $\pr\defeq \pr_{X,\Sigma}:\Pi_X
\twoheadrightarrow G_k$ (cf. exact sequence (1.2)),
meaning that $\pr \circ s :G_k\to G_k$ is the identity homomorphism. 
We will refer to $s:G_k\to \Pi_X$ as a section of the geometrically pro-$\Sigma$ 
arithmetic fundamental group $\Pi_X$.

Every inner automorphism
$\inn ^{g}:\Pi_X\to \Pi_X$ of $\Pi_X$, by an element $g \in \Delta _X$, gives rise to a
conjugate section 
$$\inn ^{g}\circ s:G_k\to \Pi_X$$
of $\Pi_X$.
We shall refer to the set
$$\{\inn ^{g}\circ s:G_k\to \Pi_X\}_{g \in \Delta _X}$$
as the set of conjugacy classes of the section $\sigma$. 

Write
$X\times X\defeq X\times _kX$, and $\iota : X\to X\times X$
for the natural diagonal embedding. The geometric point $\eta$ of $X$ determines naturally (via $\iota$)
a geometric point,
which we will also denote $\eta$, of $X\times X$. 

There exists a natural exact sequence of
profinite groups
$$1\to \pi_1(\overline {X\times X},\bar \eta)\to \pi_1(X\times X, \eta) @>{\pr_{X\times X}}>> G_k\to 1.
\tag {$1.3$}$$

Here, $\pi_1(X\times X, \eta)$ denotes the arithmetic \'etale fundamental group of $X\times X$
with base point $\eta$,  which is naturally identified with the fibre product
$\pi_1(X, \eta)\times _{G_k}\pi_1(X, \eta)$, and $\pi_1(\overline {X\times X},\bar \eta)$ is the
\'etale fundamental group of $\overline {X\times X}\defeq (X\times X)\times _k  \bar k$
with base point $\bar \eta$ (the base point $\bar \eta$ is naturally induced by $\eta$),  which is naturally
identified with the product $\pi_1(\overline X, \bar \eta)\times \pi_1(\overline X, \bar \eta)$.

Similarly, as above, we consider the maximal pro-$\Sigma$ quotient
$$\Delta_{X\times X}\defeq \pi_1(\overline {X\times X},\bar \eta)^{\Sigma}$$ 
of $\pi_1(\overline {X\times X},\bar \eta)$,
which is naturally identified with $\Delta_X\times \Delta _X$,
and the (maximal) geometrically pro-$\Sigma$ quotient
$$ \Pi_{X\times X}\defeq \pi_1(X\times X, \eta)^{(\Sigma)}\defeq \pi_1(X\times X, \eta)/\Ker (\pi_1(\overline {X\times X},\bar \eta)
\twoheadrightarrow \pi_1(\overline {X\times X},\bar \eta)^{\Sigma})$$ 
of $\pi_1(X\times X, \eta)$, which is
naturally identified with $\Pi_X\times _{G_k}\Pi_X$.

Thus, we have a natural exact sequence
$$1\to \Delta_{X\times X} \to  \Pi_{X\times X} \to G_k\to 1.$$

\subhead {1.2}
\endsubhead
Next, we recall the definition of the \'etale Chern class associated to a section of the arithmetic fundamental group $\Pi_X$.

In [Esnault-Wittenberg] \'etale cohomology classes were associated to sections of 
arithmetic fundamental groups in a quite general setting. In our setting one can define
these \'etale cohomology classes differently as follows. The following definition doesn't
work in the non-proper case, which is treated in loc. cit.. 

In what follows all scheme cohomology
groups are \'etale cohomology groups. 

First, let ${\hat \Bbb Z}^{\Sigma}$ be the maximal pro-$\Sigma$ quotient of
$\hat \Bbb Z$, and
$$M_X\defeq \Hom (\Bbb Q/\Bbb Z,(K_X^{\sep})^{\times})\otimes_{\hat \Bbb Z}{\hat \Bbb Z^{\Sigma}}.$$
Note that $M_X$ is a free ${\hat \Bbb Z}^{\Sigma}$-module of rank one, and
has a natural structure of $G_k$-module which is  isomorphic to the $G_k$-module
${\hat \Bbb Z}(1)^{\Sigma}$, where the ``(1)'' denotes a Tate twist, i.e. $G_k$ acts on
${\hat \Bbb Z}(1)^{\Sigma}$ via the $\Sigma$-part of the cyclotomic character. The $G_k$-module $M_X$ is the
module of roots of unity attached to $X$, relative to the set of primes $\Sigma$. 

Let
$$\eta _X^{\diag} \in H^2(X\times X,M_X)$$
be the \'etale Chern class which is associated to the diagonal embedding $\iota$,
or alternatively the first Chern class of the line bundle $\Cal O_{X\times X}(\iota (X))$
(cf. [Mochizuki], Proposition 1.6, for the group-theorecity of the Chern class $\eta _X^{\diag}$
in the case where the base field $k$ is a $p$-adic local field, or a finite field).

There exists a natural identification (cf. [Mochizuki], Proposition 1.1)
$$H^2(X\times X,M_X)\isom H^2(\Pi_{X\times X} ,M_X).$$
The Chern class $\eta _X^{\diag}$ corresponds, via the above identification, to an extension class
$$\eta _X^{\diag}\in  H^2(\Pi_{X\times X},M_X).$$
We shall refer to the extension class $\eta _X^{\diag}$
as the extension class of the diagonal. 

Let $s:G_k\to \Pi_X$ be a group-theoretic section of the natural projection $\Pi_X\twoheadrightarrow G_k$ (cf. 1.1). 

Let
$$1\to M_X\to \Cal D\to \Pi_{X\times X}  \to 1\tag {$1.4$}$$
be a group extension whose class in  $H^2(\Pi_{X\times X},M_X)$ coincides with
the diagonal class $\eta _X^{\diag}$. 

By pulling back the group extension $(1.4)$ by the continuous
injective homomorphism
$$(s,\id): G_k\times_{G_k} \Pi_X    \to \Pi_{X\times X}$$
we obtain a natural commutative diagram:

$$
\CD
1@>>>   M_X    @>>>     \Cal D_s   @>>>  G_k \times _{G_k} \Pi_X        @>>> 1\\
  @.        @V{\id}VV           @VVV              @V{(s,\id)} VV \\
1     @>>> M_X  @>>>   \Cal D    @>>>  \Pi_{X\times X}  @>>> 1
\endCD
$$
where the right square is cartesian.
Further, via the natural identification $G_k\times _{G_k} \Pi_X \isom \Pi_X$, 
the upper group extension $\Cal D_s$ in the above diagram corresponds to a group extension
(which we denote also $\Cal D_s$)
$$1 @>>>  M_X    @>>>     \Cal D_s   @>>>  \Pi_X\to 1.$$ 

We shall refer to the class $[\Cal D_s]$
of the extension $\Cal D_s$ in $H^2(\Pi_X,M_X)$ as the extension class
associated to the section $s$.

\definition{Definition 1.2.1 (The ($\Sigma$)-\'Etale Chern Class associated to a Section)}
We define the ($\Sigma$)-\'etale Chern class $c(s)\in  H^2(X,M_X)$
associated to the section $s$ as the element corresponding
to the above extension class $[\Cal D_s]$, which is associated to the section $s$,
via the natural identification
$H^2(\Pi_X,M_X)\isom H^2(X,M_X)$ (cf. loc. cit.).
\enddefinition

\definition{Remark 1.2.2} It is easy to verify that our definition of the Chern class $c(s)$
associated to the section $s$, which can be carried out in the more general setting of a proper and
smooth variety $X$ over $k$, coincides in this case, i.e. under the assumption of being proper,
with the definition in  [Esnault-Wittenberg].
\enddefinition

\subhead {1.3}
\endsubhead
Next, we recall the definition of a system of neighbourhoods  of a group-theoretic
section of the arithmetic fundamental group $\Pi _X$.

The profinite group $\Delta _X$ being topologically finitely generated, there
exists a sequence of characteristic open subgroups
$$...\subseteq \Delta _X[i+1]\subseteq \Delta _X[i]\subseteq...\subseteq \Delta _X[1]\defeq \Delta _X$$
of $\Delta_X$, where $i\ge 1$ ranges over all positive integers, such that 
$$\bigcap _{i\ge 1}\Delta _X[i]=\{1\}.$$

In particular, given a group-theoretic section $s:G_k\to \Pi_X$ of $\Pi_X$, we obtain open subgroups
$$\Pi _X[i,s]\defeq s(G_k).\Delta _X[i] \subseteq \Pi_X$$
($s(G_k)$ denotes the image of $G_k$ in $\Pi_X$ via the section $s$) of $\Pi_X$, 
whose intersection coincide with $s(G_k)$,
and which correspond to a tower of finite \'etale (not necessarily Galois) covers
$$...\to X_{i+1}[s]\to X_{i}[s]\to ...\to X_1[s]\defeq X$$
defined over $k$. We will refer to the set 
$$\{X_i[s]\}_{i\ge 1}$$
as a system of neighbourhoods of the section $s$.

Note that for every positive integer $i$, the open subgroup $\Pi _X[i,s]$ of $\Pi_X$
is naturally identified with the geometrically pro-$\Sigma$ arithmetic \'etale fundamental group
$\pi_1(X_i[s], \eta_i)^{(\Sigma)}$ of $X_i[s]$, the geometric point $\eta_i$ of $X_i[s]$ being 
naturally induced by the geometric point $\eta$ of $X$, and
sits naturally in the following exact sequence
$$1\to \Delta _X[i,s]\to \Pi_X[i,s] \to G_k\to 1,$$
which inserts in the following commutative diagram:
$$
\CD
1 @>>>  \Delta _X[i,s]     @>>> \Pi_X[i,s]     @>>> G_k @>>> 1\\
  @.        @VVV           @VVV              @V{\id}VV \\
1 @>>>  \Delta _X       @>>> \Pi_X       @>>> G_k   @>>> 1
\endCD
$$
where the two left vertical homomorphisms are the natural inclusions.

In particular, by the very definition of $\Pi_X[i,s]$,
the section $s$ restricts naturally to a group-theoretic section
$$s_i:G_k\to \Pi_X[i,s]$$
of the natural projection $\Pi_X[i,s]\twoheadrightarrow G_k$,
which fits into the following commutative diagram:
$$
\CD
G_k @>{s_i}>> \Pi_X[i,s] \\
@V{\id}VV     @VVV  \\
G_k @>{s}>> \Pi_X
\endCD
$$
where the right vertical homomorphism is the natural inclusion.

\proclaim {Lemma 1.3.1} For each positive integer $i$, the image of the Chern class
$c(s_{i+1})\in H^2 (X_{i+1}[s],M_X)$
associated to the section $s_{i+1}:G_k\to \Pi_X[i+1,s]$ in $H^2 (X[i,s],M_X)$, via the
corestriction homomorphism
$\cor:H^2 (X_{i+1}[s],M_X)\to H^2 (X_i[s],M_X)$, 
coincides with
the Chern class $c(s_i)$ associated to the section $s_i:G_k\to \Pi_X[i,s]$.
\endproclaim

\demo {Proof} Follows from the various above definitions, and the fact that the image of the Chern class
$c(\eta ^{\diag}_{X_{i+1}})$ of the
diagonal morphism $\iota _{i+1}:X_{i+1}\to X_{i+1}\times X_{i+1}$ in $H^2(X_{i}\times X_{i},M_X)$,
via the corestriction homomorphism $\cor:H^2(X_{i+1}\times X_{i+1},M_X)\to H^2(X_{i}\times X_{i},M_X)$,
coincides with the Chern class $c(\eta ^{\diag}_{X_{i}})$ of the diagonal  morphism
$\iota _{i}:X_{i}\to X_{i}\times X_{i}$.
\qed
\enddemo

\definition{Definition 1.3.2  (The pro-($\Sigma$)-\'Etale Chern Class associated to a Section)} Let
$\underset{i\ge 1}\to{\varprojlim}\ H^2 (X_i[s],M_X)$
be the projective limit of the  $H^2 (X_i[s],M_X)$'s, where the transition homomorphisms are the
corestriction homomorphisms. We define the pro-$(\Sigma)$-\'etale Chern class associated to the section $s$,
relative to the system of neighbourhoods $\{X_i[s]\}_{i\ge 1}$, as the element
$$\hat c(s)\defeq (c(s_i))_{i\ge 1}\in \underset{i\ge 1}\to{\varprojlim}\ H^2 (X_i[s],M_X)$$
(cf. Lemma 1.3.1).
\enddefinition

\subhead {1.4}
\endsubhead
Next, we will introduce the notion of (uniformly) good sections of arithmetic fundamental groups.
 
For each positive $\Sigma$-integer $n$, meaning that $n$ is an integer which is divisible only by primes in
$\Sigma$, the Kummer exact sequence in \'etale topology
$$1\to \mu_n\to \Bbb G_m @>n>> \Bbb G_m\to 1,$$
induces naturally, for each positive integer $i$, an exact sequence of abelian groups
$$0 \to \Pic (X_i[s])/n \Pic (X_i[s]) \to  H^2(X_i[s],\mu_{n}) \to _{n} \Br (X_i[s])
\to 0.\tag {$1.5$} $$

For positive $\Sigma$-integers $m$ and $n$, with $n$ divides $m$, we have a commutative diagram:
$$
\CD
0@>>> \Pic (X_i[s])/m \Pic (X_i[s])      @>>>     H^2(X_i[s],\mu_{m})  @>>> _{m} \Br (X_i[s]) @>>> 0\\
  @.        @VVV           @VVV              @VVV \\
0@>>> \Pic (X_i[s])/n \Pic (X_i[s])      @>>>     H^2(X_i[s],\mu_{n})  @>>> _{n} \Br (X_i[s])    @>>> 0
\endCD
$$
where the lower and upper horizontal sequences are the above exact sequence $(1.5)$,
and the vertical homomorphisms are the natural homomorphisms. 

Here
$\Pic\defeq H^1(\ ,\Bbb G_m)$ denotes the Picard group, $\Br\defeq H^2(\ ,\Bbb G_m)$
the Brauer-Grothendieck cohomological group, and for a positive integer $n$: $_{n} \Br\subseteq \Br$
denotes the subgroup of $\Br$ which is annihilated by $n$. 

By taking
projective limits, the above diagram induces naturally, for every positive integer $i$, the following exact sequence
$$0\to \underset{n\ \Sigma-\text {integer}}\to{\varprojlim}  \Pic (X_i[s])/n \Pic (X_i[s]) \to H^2(X_i[s],M_X)\to
\underset{n\ \Sigma-\text {integer}}\to{\varprojlim} _n \Br (X_i[s])\to 0.$$

We will denote by
$$\Pic (X_i[s])^{\wedge,\Sigma}\defeq \underset{n\ \Sigma-\text {integer}}
\to{\varprojlim} \Pic (X_i[s])/n \Pic (X_i[s])$$
the $\Sigma$-adic completion of the Picard group $\Pic (X_i[s])$, and
$$T_{\Sigma}\Br(X_i[s]) \defeq \underset{n\ \Sigma-\text{integer}}\to{\varprojlim} _{n} \Br (X_i[s])$$
the $\Sigma$-Tate module of the Brauer group $\Br (X_i[s])$. Thus, we have a natural exact sequence
$$0\to\Pic (X_i[s])^{\wedge,\Sigma}\to H^2(X_i[s],M_X)\to T_{\Sigma}\Br(X_i[s])\to 0. \tag {$1.6$}$$

In what follows we will identify $\Pic (X_i[s])^{\wedge,\Sigma}$ with its image in  $H^2(X_i[s],M_X)$,
and refer to it as the Picard part of $H^2(X_i[s],M_X)$.

Let 
$$s:G_k\to \Pi_X[1,s]\defeq \Pi_X$$
be a continuous group-theoretic section of $\Pi_X$ as above. For each positive integer $i$,
let 
$$s_i:G_k\to \Pi_X[i,s]$$
be the induced group-theoretic section of $\Pi_X[i,s]$.

By pulling back cohomology classes via the section $s_i$,
and bearing in mind the natural identification 
$$H^2(\Pi_X[i,s],M_X)\isom H^2(X_i[s],M_X)$$
(cf. [Mochizuki], Proposition 1.1), we obtain a natural restriction homomorphism
$$s_i^{\star}: H^2(X_i[s],M_X)\to H^2(G_k,M_X).$$

Finally, observe that if $k'/k$ is a finite extension, and $X_{k'}\defeq X\times _kk'$,
then we have a natural commutative diagram:
$$
\CD
1 @>>>  \Delta _X     @>>> \Pi_{X_{k'}}\defeq \pi_1(X_{k'} ,\eta)^{(\Sigma)}     @>>> G_{k'} @>>> 1\\
  @.        @V{\id}VV           @VVV              @VVV \\
1 @>>>    \Delta _X      @>>> \Pi_X       @>>> G_k   @>>> 1
\endCD
$$
where the vertical arrows are the natural inclusions, and the far right square is cartesian.
In particular, the section $s_k\defeq s:G_k\to \Pi_X$
induces naturally a group-theoretic section
$$s_{k'}:G_{k'}\to  \Pi_{X_{k'}}$$
of the arithmetic fundamental group $\Pi_{X_{k'}}$.

\definition{Definition 1.4.1 (Good and Uniformly Good Sections of
Arithmetic Fundamental Groups)}
We say that the section $s$ is a good group-theoretic section, relative to the system
of neighbourhoods $\{X_i[s]\}_{i\ge 1}$, if for every positive integer $i$
the above homomorphism
$s_i^{\star}: H^2(X_i[s],M_X)\to H^2(G_k,M_X)$ 
annihilates the Picard part $\Pic (X_i[s])^{\wedge,\Sigma}$ of
$H^2(X_i[s],M_X)$. In other words
the section $s$ is good if 
$$\Pic (X_i[s])^{\wedge,\Sigma}\subseteq \Ker s_i^{\star},$$
for every positive integer $i$.

We say that the section $s$ is uniformly good, relative to the system
of neighbourhoods $\{X_i[s]\}_{i\ge 1}$, if for every finite extension $k'/k$
the induced section
$s_{k'}:G_{k'}\to \Pi_{X_{k'}}$
is good, relative to the system of neighbourhoods of $s_{k'}$ which is naturally induced by the $\{X_i[s]\}_{i\ge 1}$.
\enddefinition

It is easy to see that the above definition is independent of the given system of neighbourhoods $\{X_i[s]\}_{i\ge 1}$
of the section $s$. We will refer to a section satisfying the conditions in Definition 1.4.1 as good, or uniformly good, without necessarily 
specifying a system of neighbourhoods of the section.

The above Definition 1.4.1 is motivated by the fact that a necessary condition for a group-theoretic
section $s:G_k\to \Pi_X$ of $\Pi_X$ to be point-theoretic, i.e. arises from a $k$-rational point
$x\in X(k)$ (cf. Definition 3.1.1), is that the section $s$ is uniformly good in the sense of
Definition 1.4.1 (cf. Proposition 1.5.2).

\subhead {1.5}
\endsubhead
Next, we will introduce the notion of well-behaved group-theoretic sections of $\Pi_X$. We use the same notations as above.

For every positive integer $i$, we have the following commutative diagram:
$$
\CD
0@>>> \Pic (X_{i+1}[s])^{\wedge,\Sigma}      @>>>     H^2(X_{i+1}[s],M_X)  @>>> T_{\Sigma}\Br (X_{i+1}[s])@>>> 0\\
  @.        @V{\norm}VV           @V{\cor}VV              @VVV \\
0@>>> \Pic (X_{i}[s])^{\wedge,\Sigma}      @>>>     H^2(X_i[s],M_X)  @>>> T_{\Sigma}\Br (X_i[s])    @>>> 0
\endCD
$$
where the horizontal sequences are the exact sequence $(1.6)$, the far left vertical map is naturally induced by
the norm homomorphism 
$\norm: \Pic (X_{i+1}[s])\to \Pic (X_{i}[s])$,
and the middle vertical map is the corestriction homomorphism. The far right vertical map is
naturally induced by the vertical middle map, and the commutativity of the above diagram.
Here, the $\{X_{i}[s]\}_{i\ge 1}$ form a system of
neighbourhoods of the section $s$.

In particular, by passing to the projective limit, via the above explained maps, we obtain a natural exact sequence
$$0 @>>> \underset{i\ge 1}\to{\varprojlim}\Pic (X_i[s])^{\wedge,\Sigma} @>>>   \underset{i\ge 1}
\to{\varprojlim}\ H^2(X_i[s],M_X)@>>> \underset{i\ge 1}\to{\varprojlim}\ T_{\Sigma}\Br (X_i[s]) .\tag {1.7}$$

We shall refer to the image of $\underset{i\ge 1}\to{\varprojlim}\Pic (X_i[s])^{\wedge,\Sigma}$ in
$\underset{i\ge 1}\to{\varprojlim}\ H^2(X_i[s],M_X)$, via the above homomorphism in (1.7), as the Picard part of
$\underset{i\ge 1}\to{\varprojlim}\ H^2(X_i[s],M_X)$.

The Chern classes $c(s_i)\in H^2(X_i[s],M_X)$ associated to the sections $s_i:G_k\to \Pi_X[i,s]$, for every positive integer
$i$, determine naturally an element $\hat c (s)\defeq (c(s_i))_{i\ge 1}\in \underset{i\ge 1}\to{\varprojlim}
H^2(X_i[s],M_X)$ (cf. Lemma 1.3.1): the pro-$(\Sigma)$-\'etale Chern class associated to the section $s$.

\definition{Definition 1.5.1 (Well-Behaved Sections of Arithmetic Fundamental Groups)}\ \ \ \ \ \
We say that the section $s$ is a well-behaved group-theoretic section, 
relative to the system
of neighbourhoods $\{X_i[s]\}_{i\ge 1}$,
if the image of the
pro-Chern class $\hat c (s)\in \underset{i\ge 1}\to{\varprojlim}\
H^2(X_i[s],M_X)$ in $\underset{i\ge 1}\to {\varprojlim}\ T_{\Sigma}\Br (X_i[s])$,
via the natural homomorphism
$\underset{i\ge 1}\to {\varprojlim}\ H^2(X_i[s],M_X)\to \underset{i\ge 1}\to {\varprojlim}\ T_{\Sigma}\Br (X_i[s])$
in $(1.7)$, equals $0$. In other words the section $s$ is well-behaved if the associated pro-Chern class
$\hat c (s)$ lies in the Picard part $\underset{i\ge 1}\to{\varprojlim}\ \Pic (X_i[s])^{\wedge,\Sigma}$ of
$\underset{i\ge 1}\to{\varprojlim}\ H^2(X_i[s],M_X)$.

We say that the section $s$ is uniformly well-behaved, relative to the system
of neighbourhoods $\{X_i[s]\}_{i\ge 1}$, if for every finite extension $k'/k$
the induced section $s_{k'}:G_{k'}\to \Pi_{X_{k'}}$ is well-behaved in the above sense, relative to the system
of neighbourhoods of $s_{k'}$ which is naturally induced by the system $\{X_i[s]\}_{i\ge 1}$.
\enddefinition

One easily verifies that the above definition is independent of the given system of neighbourhoods $\{X_i[s]\}_{i\ge 1}$
of the section $s$. We will refer to a section satisfying the conditions in Definition 1.5.1 as well-behaved, or uniformly well-behaved, without necessarily 
specifying a system of neighbourhoods of the section.

The above Definition 1.5.1 is motivated by the following. 

Assume that the section $s$ is point-theoretic, meaning that
$s\defeq s_x:G_k\to \Pi_X$
arises from a $k$-rational point $x\in X(k)$ (cf. Definition 3.1.1). Then there exists a compatible system of rational points
$\{x_i\in X_i[s](k)\}_{i\ge 1}$, i.e. $x_{i+1}$ maps to $x_i$ via the natural morphism $X_{i+1}[s]\to X_i[s]$.

For every positive integer $i$, let $\Cal O(x_i)\in \Pic (X_i[s])$
be the degree $1$ line bundle associated to $x_i$. The Chern class
$c(s_i) \in H^2(X_i[s],M_X)$ which is associated to the section $s_i$ (cf. Definition 1.2.1)
coincides with the \'etale Chern class
$c(x_i)\in H^2(X_i[s],M_X)$ associated to the line bundle $\Cal O(x_i)$
(cf. [Mochizuki3], Lemma 4.2).
Thus, the section $s=s_x$ is well-behaved in this case, since the pro-Chern class $\hat c(s)$
is the ``pro-Picard element'' induced by the $(\Cal O(x_i))_{i\ge 1}$.

The link between good and well-behaved group-theoretic sections is given in the following Proposition.

\proclaim{Proposition 1.5.2} Assume that $s:G_k\to \Pi_X$ is a well-behaved group-theoretic section.
Then $s$ is a good group-theoretic section. In particular, If the section $s$ is point-theoretic (cf. Definition 3.1.1), 
i.e. arises from a $k$-rational point $x\in X(k)$, then the section $s$ is a uniformly good
group-theoretic section.
\endproclaim

\demo{Proof} For every positive integer $i$, there exists  a natural pairing
$$\Pic (X_i[s])\times \Br (X_i[s])\to \Br (k)$$
(cf. [Lichtenbaum], 3, and the Appendix in this paper),
which induces, for every positive $\Sigma$-integer $n$, a natural pairing
$$\Pic (X_i[s])/n\Pic (X_i[s])\times \ _n \Br (X_i[s])\to \ _n\Br (k).$$

We have natural restriction homomorphisms, which are induced by pulling 
back cohomology classes via the sections 
$s_i$:
$$s_{i,n}^{\star}: H^2(X_i[s],M_X/nM_X)\to H^2(G_k,M_X/nM_X).$$

The proof of the Proposition follows from the observation (cf. Lemma A.4 in the Appendix) that the image of
the Picard part of  $H^2(X_i[s],M_X/nM_X)$ in $H^2(G_k,M_X/nM_X)$ via the map
$s_{i,n}^{\star}$, coincides with the image of
$\Pic (X_i[s])/n\Pic (X_i[s])$ in $_n\Br (k)$ via the above pairing, as a result of pairing the elements of
$\Pic (X_i[s])/n\Pic (X_i[s])$ with the image of
$c_n(s_i)\in H^2(X_i[s],M_X/nM_X)$ in $_n \Br (X_i[s])$ via the natural homomorphism
$H^2(X_i[s],M_X/nM_X) \to _n \Br (X_i[s])$. 

Here $c_n(s_i)\in H^2(X_i[s],M_X/nM_X)$ denotes the image of the
Chern class $c(s_i)$ via the natural homomorphism 
\newline
$H^2(X_i[s],M_X)\to H^2(X_i[s],M_X/nM_X)$. This later image is $0$
if the section $s$ is well-behaved by definition.

Finally the last assertion follows from the discussion preceding Proposition 1.5.2.
\qed
\enddemo

\subhead {1.6}
\endsubhead
In what follows, and unless we specify otherwise, we will assume that the field $k$ is a $p$-adic local field,
i.e. $k$ is a finite extension of $\Bbb Q_p$ for some fixed prime integer $p>0$. 

In this case, and in
the framework of the above discussion, one has the following more precise statement.

\proclaim{Proposition 1.6.1} Assume that $k$ is a $p$-adic local field.
If the section $s=s_x$ is point-theoretic, i.e. arises from a $k$-rational point $x\in X(k)$ (cf. Definition 3.1.1), then
$$\underset{i\ge 1}\to {\varprojlim}\ _n\Br (X_i[s])=0,$$ 
for every $\Sigma$-integer $n$.
In particular, if the section $s$ is point-theoretic, then we have a natural isomorphism
$$\underset{i\ge 1}\to{\varprojlim}\ \Pic (X_i[s])/n\Pic(X_i[s])\isom \underset{i\ge 1}\to{\varprojlim}\ H^2(X_i[s],M_X/nM_X).$$
\endproclaim

\demo {Proof} Follows from the proof of Proposition 1.6.3 below, and the fact that $X_i(k)\neq \varnothing$ implies that
$\period (X_i)=\index (X_i)=1$.
\qed
\enddemo

In the framework of the $p$-adic version of the Grothendieck
anabelian section conjecture (cf. 3.1) it is natural to ask the following question.

\definition{Question 1.6.2} Assume that $k$ is a $p$-adic local field. Is 
$\underset{i}\to {\varprojlim}\ _n\Br (X_i[s])=0$
for every  group-theoretic section $s:G_k\to \Pi_X$, and every $\Sigma$-integer $n$, for a given non-empty set of prime integers $\Sigma$? 
If this is the case, then 
every group-theoretic section $s:G_k\to \Pi_X$ would be a well-behaved section in the sense of Definition 1.5.1 (cf. exact sequence (1.7)). 
\enddefinition

In connection with Question 1.6.2 we have the following.

\proclaim{Proposition 1.6.3} Assume that $k$ is a $p$-adic local field. Then the following properties hold.

\noindent
{\rm (i)}\ Assume that $\Sigma\cap \{p\}=\varnothing$. Then $\underset{i\ge 1}\to {\varprojlim}\ _n\Br (X_i[s])=0$
for every $\Sigma$-integer $n$.
In particular, the section $s$ is well-behaved (hence is good by Proposition 1.5.2) in this case (compare with [Esnault-Wittenberg], Corollary 3.4).

\noindent
{\rm (ii)} \ Assume that $p\in \Sigma$. Then $\underset{i\ge 1}\to {\varprojlim}\ _{p^m}\Br (X_i[s])=0$ for every positive integer $m$,
which would imply that the section $s$ is well-behaved (cf. assertion (i)), if and only if
$$\underset{i\ge 1}\to {\varinjlim} \Pic (X_i[s])/p\Pic (X_i[s])=0,$$
where the transition homomorphisms
$$\Pic (X_i[s])/p\Pic (X_i[s])\to \Pic (X_{i+1}[s])/p\Pic (X_{i+1}[s])$$
in the inductive limit are induced by the natural pull back of line bundles homomorphisms, via the natural morphisms
$X_{i+1}[s]\to X_{i}[s]$.
\endproclaim

\demo{Proof} First, note that we have, for every $\Sigma$-integer $n$, a natural exact sequence
$$0 @>>> \underset{i\ge 1}\to{\varprojlim}\ \Pic (X_i[s])/n \Pic(X_i[s]) @>>>   \underset{i\ge 1}
\to{\varprojlim}\ H^2(X_i[s],M_X/nM_X)@>>> \underset{i\ge 1}\to{\varprojlim}\ _n
\Br (X_i[s]), \tag {$1.8$}$$
from which follows directly that the vanishing of $\underset{i\ge 1}\to{\varprojlim} _n
\Br (X_i[s])$ implies that the section $s$ is well-behaved.

Second, it follows from the Tate-Lichtenbaum duality between $\Pic(X)$ and $\Br(X)$ for a proper,
smooth, and geometrically connected curve over $k$ (cf. [Lichtenbaum]), that (for every $\Sigma$-integer $n$) the dual of 
the projective limit $\underset{i\ge 1}\to {\varprojlim}  _n \Br (X_i[s])$, where $\{X_i[s])\}_{i\ge 1}$ is a system of neighbourhoods
of the section $s$, is naturally identified with the inductive limit  
$\underset{i\ge 1}\to {\varinjlim} \Pic (X_i[s])/n\Pic (X_i[s])$.

Indeed, one can easily verify that the transition morphisms in the above projective limit correspond, by the above
duality, to the natural pull back of line bundles homomorphisms in the inductive limit. 

So in order to prove assertion 
{\rm (i)} it suffices to prove if $p\notin \Sigma$, for a $\Sigma$-integer $n$, that 
$\underset{i\ge 1}\to {\varinjlim} \Pic (X_i[s])/n\Pic (X_i[s]) =0$, for which it suffices to prove the following claim.

\proclaim{Claim 1.6.4}
Assume that $p\notin \Sigma$. Then given a line bundle $\Cal L\in  \Pic(X)$, and a positive $\Sigma$-integer $n$, there exists a positive integer $i$ such
that the pull back of $\Cal L$ to $X_i[s]$, via the natural morphism $X_i[s]\to X$, is the $n$-th power of an
element of $\Pic (X_i[s])$.
\endproclaim

\demo {Proof of Claim 1.6.4}
Let $\Cal L\in  \Pic(X)$, and $n$ a positive integer, not necessarily a $\Sigma$-integer for the moment.

Over the algebraic closure $\bar k$ of $k$, and for
every positive integer $i$,  one has natural identifications
$$H^2(X_i[s]\times _k\bar k, M_X/nM_X)\isom \Pic (X_i[s]\times _k\bar k)/n\Pic (X_i[s]\times _k\bar k)\isom \Bbb Z/n\Bbb Z,$$
via which the natural restriction morphisms
$$H^2(X_i[s]\times _k\bar k, M_X/nM_X)\to H^2(X_{i+1}[s]\times _k\bar k, M_X/nM_X)$$
are given by multiplication by the degree $\vert X_{i+1}[s]:X_i[s]\vert$ of the morphism $X_{i+1}[s]\to X_i[s]$.

Let $\overline {\Cal L}$ be the image
of $\Cal L$ in $\Pic (X\times _k\bar k)$. Then, by the above considerations, there exists a positive integer
$i$ such that the image $\overline {\Cal L}_i$
of $\overline {\Cal L}$ in $\Pic (X_i[s]\times _k\bar k)$, by pull back via the natural morphism
$X_i[s]\times _k\bar k\to X\times _ k\bar k$,
satisfies $\overline {\Cal L}_i={\overline {\Cal L'}}^{n}$
for some element $\overline {\Cal L'}\in \Pic (X_i[s]\times _k\bar k)$. 

The absolute Galois group $G_k$ of $k$ acts
on $\Pic (X_i[s]\times _k\bar k)$,
and $\overline {\Cal L}_i$ is a fixed element under this action. By observing the action of $G_k$ on
$\overline {\Cal L}_i={\overline {\Cal L'}}^{n}$, one sees that $\overline {\Cal L'}$ is fixed under the action of $G_k$,
up to multiplication by some $n$-torsion element of
$\Pic (X_i[s]\times _k\bar k)$. The $n$-torsion subgroup of $\Pic (X_i[s]\times _k\bar k)$ vanishes upon pull back to
$\Pic (X_j[s]\times _k\bar k)$ for some $j>i$. Thus, after replacing $X_i[s]\times _k\bar k$ by $X_j[s]\times _k\bar k$,
we can assume that $\overline {\Cal L'}$ is fixed under the action of $G_k$. Hence $\overline {\Cal L}_i=
{\overline {\Cal L'}}^{n}$, where
$\overline {\Cal L}_i\in H^0(G_k,\Pic (X_i[s]\times _k\bar k))$, and $\overline {\Cal L'}\in H^0(G_k,\Pic (X_i[s]
\times _k\bar k))$. 

Recall that there exists a canonical exact sequence (cf. [Lichtenbaum], 2)
$$0\to\Pic (X_i[s]) \to H^0(G_k,\Pic (X_i[s]\times _k\bar k))\to \Br (k)\to \Br (X_i[s]).$$
The image of an element of $H^0(G_k,\Pic (X_i[s]\times _k\bar k))$ in $\Br (k)$ is the obstruction
for this element to arise from an element of $\Pic (X_i[s])$. This obstruction lies in the $\Sigma'$ (torsion) primary
part of $\Br(k)$ in light of the existence of the section $s$ by a result of Stix [Stix], where $\Sigma'\defeq \Primes \setminus \Sigma$. 

Actually in [Stix], $\Sigma=\Primes$, but similar arguments as in loc. cit. 
yield a similar result. Namely, the existence of the section $s_i$ implies that the period $\period(X_i[s])$, and the index $\index(X_i[s])$, of $X_i[s]$, are divisible only by 
primes in $\{p\}\cup \{\Sigma'\}$, hence only by primes in $\Sigma'$ since we assumed $p\notin \Sigma$.

The element $\overline {\Cal L}_i$ arises from an element of $\Pic (X_i[s])$, so its Brauer obstruction vanishes.
In particular, the Brauer obstruction of $\overline {\Cal L'}$ is annihilated by $n$ in $\Br(k)$. If $n$
is a $\Sigma$-integer, this Brauer obstruction necessarily vanishes by the above mentioned result of Stix. Hence
$\Cal L_i={\Cal L'}^n$, where $\Cal L_i\in \Pic (X_i[s])$ is the pull-back of $\Cal L$, and ${\Cal L}'\in \Pic (X_i[s])$.

This finishes the proof of Claim 1.6.4. 
\qed
\enddemo

The proof of the assertion {\rm (ii)} is clear in light of the proof
of the above Claim 1.6.4. 

This finishes the proof of Proposition 1.6.3.
\qed
\enddemo

In fact one can say more in the case where $k$ is a $p$-adic local field about the group
$\underset{i\ge 1}\to {\varinjlim} \Pic (X_i[s])/p\Pic (X_i[s])$. The vanishing of this group is an obstruction for 
the section $s$ to be well-behaved in the case where $p\in \Sigma$, by Proposition 1.6.3, (ii).

\proclaim{Lemma 1.6.5} Assume that $k$ is a $p$-adic local field, and $p\in \Sigma$. 
Then there exists a natural injective homomorphism
$$\underset{i\ge 1}\to {\varinjlim} \Pic (X_i[s])/p\Pic (X_i[s])\hookrightarrow _p\Br(k).$$

In particular, exactly one of the two following cases occur.

\noindent
{\rm (i)} Either $\underset{i\ge 1}\to {\varinjlim} \Pic (X_i[s])/p\Pic (X_i[s])=0$, which is equivalent
to $\underset{i\ge 1}\to {\varprojlim}\ _p\Br (X_i[s])=0$.

\noindent
{\rm (ii)} Or $\underset{i\ge 1}\to {\varinjlim}  \Pic (X_i[s])/p\Pic (X_i[s])$ is a finite group of cardinality $p$,
which is equivalent to $\underset{i\ge 1}\to {\varprojlim}\ _p\Br (X_i[s])\isom \Bbb Z/p\Bbb Z$.
\endproclaim

\demo{Proof} The Kummer exact sequence
$$1\to \mu_p\to \Bbb G_m@>p>> \Bbb G_m\to 1$$
in \'etale topology induces, for every positive integer $i$, a natural exact sequence
$$0\to \Pic (X_i[s])/p\Pic (X_i[s])\to H^2(X_i[s],\mu_p)\to \ _p\Br (X_i[s])\to 0.$$
By passing to the inductive limit one obtains the following exact sequence
$$0\to \underset{i\ge 1}\to {\varinjlim}
\Pic (X_i[s])/p\Pic (X_i[s])\to \underset{i\ge 1}\to {\varinjlim}\ H^2(X_i[s],\mu_p)\to \underset{i\ge 1}\to {\varinjlim}
\ _p\Br (X_i[s])\to 0.$$
The assertion then follows by observing that $\underset{i}\to {\varinjlim}\ H^2(X_i[s],\mu_p)$ , which is naturally isomorphic to
$\underset{i}\to {\varinjlim}\ H^2(\Pi_X[i,s],\mu_p)$ ,
is naturally identified with
$H^2(s(G_k),\mu_p)$, which is $_p\Br(k)$, and the later is a finite group of cardinality $p$.
\qed
\enddemo

The following Proposition establishes the link between several properties of a given group-theoretic section
$s:G_k\to \Pi_X$ of $\Pi_X$, in the case where $k$ is a $p$-adic local field.

\proclaim{Proposition 1.6.6} Assume that the field $k$ is a $p$-adic local field, and $p\in \Sigma$.
Write $\Sigma'\defeq \Primes\setminus \Sigma$.
Consider the following properties.

\noindent
{\rm (i)}\ The section $s$ is a good group-theoretic section in the sense of Definition 1.4.1.

\noindent
{\rm (ii)}\ The section $s$ is a well-behaved group-theoretic section in the sense of Definition 1.5.1.

\noindent
{\rm (iii)}\ For every positive integer $i$, one has $X_i[s](k^{\tame})\neq \varnothing$, where $k^{\tame}$ is the
maximal tamely ramified extension of $k$.

\noindent
{\rm (iv)}\ For every positive integer $i$, there exists a finite extension $\ell_i$ of $k$ with
$\gcd ([\ell_i:k],p)=1$ such that  $X_i[s](\ell_i)\neq \varnothing$.

\noindent
{\rm (v)}\ For every positive integer $i$, the period $\period (X_i[s])$, and index $\index (X_i[s])$, of $X_i[s]$ are $\Sigma'$-integer, 
i.e. are divisible only by primes in $\Sigma'$.
Moreover, if $2\in \Sigma$ then $\period (X_i)=\index(X_i)$.

\noindent
{\rm (vi)}\ For every positive integer $i$, the kernel of the natural homomorphism  $\Br (k)\to \Br(X_i[s])$ is contained in the $\Sigma'$-primary (torsion) part of $\Br(k)$.

\noindent
{\rm (vii)}\ $\underset{i\ge 1}\to {\varinjlim} \Pic (X_i[s])/p\Pic (X_i[s])=0$,
where the transition homomorphisms
\newline
$\Pic (X_i[s])/p\Pic (X_i[s])\to \Pic (X_{i+1}[s])/p\Pic (X_{i+1}[s])$
in the inductive limit are induced by the natural pull back of line bundles homomorphisms via the natural morphisms
$X_{i+1}[s]\to X_{i}[s]$.

\noindent
{\rm (viii)}\ For every positive integer $i$, the elementary obstruction $\ob(X_i[s])$ vanishes (see [Wittenberg], and
[Borovoi-Colliot-Th\'el\`ene-Skorobogatov], for the definition of the elementary obstruction).

Then we have the following equivalences
$$(i) \Longleftrightarrow (ii)  \Longleftrightarrow (iii) \Longleftrightarrow (iv) \Longleftrightarrow (v) \Longleftrightarrow (vi) \Longleftrightarrow (vii).$$ 

Moreover, in the case where $\Sigma=\Primes$ we have the equivalence
$$(v) \Longleftrightarrow (viii).$$
\endproclaim 

\demo {Proof} First we prove {\rm (i)}$\Longleftrightarrow${\rm (ii)}. 

For every positive integer $i$, the natural
pairing $\Pic (X_i[s])\times \Br (X_i[s])\to \Br (k)$
is perfect (cf. [Lichtenbaum]). We have, for every positive $\Sigma$ integer $n$, an induced
pairing $\Pic (X_i[s])/n\Pic (X_i[s])\times \ _n \Br (X_i[s])\to \ _n\Br (k)$. 
Also, we have natural homomorphisms 
$s_i$ : $s_{i,n}^{\star}: H^2(X_i[s],M_X/nM_X)\to H^2(G_k,M_X/nM_X)$,
induced by pulling back cohomology classes via the sections
$s_i$, for each positive integer $i$.

The equivalence {\rm (i)}$\Longleftrightarrow${\rm (ii)} follows from the observation (cf. Lemma A.4 in the Appendix)
that the image of the Picard part of  $H^2(X_i[s],M_X/nM_X)$ in $H^2(G_k,M_X/nM_X)$ via $s_{i,n}^{\star}$
coincides with the image of
$\Pic (X_i[s])/n\Pic (X_i[s])$ in $_n\Br (k)$ via the above pairing, as a result of pairing 
the elements of $\Pic (X_i[s])/n\Pic (X_i[s])$
with the image of
$c_n(s_i)\in H^2(X_i[s],M_X/nM_X)$ in $_n \Br (X_i[s])$ via the natural homomorphism
$H^2(X_i[s],M_X/nM_X) \to _n\Br (X_i[s])$. Here $c_n(s_i)\in H^2(X_i[s],M_X/nM_X)$ denotes the image of the Chern
class $c(s_i)$ via the natural map $H^2(X_i[s],M_X)\to H^2(X_i[s],M_X/nM_X)$. 

Next, we prove {\rm (iv)} $\Longleftrightarrow$ {\rm (iii)}. 

The implication {\rm (iv)}$\Rightarrow${\rm (iii)} is clear.
We prove {\rm (iii)} $\Rightarrow$ {\rm (iv)}.

Let $\Cal X_i$ be the minimal regular model with normal crossings
of $X_i[s]$ over the ring of integers $\Cal O_k$ of $k$. The assumption {\rm (iii)} implies that there exists an
irreducible component of the special fibre of  $\Cal X_i$ whose multiplicity $e$ is prime to $p$. After eventually passing
to the finite totally ramified extension of $k$ of degree $e$ we can assume, without loss of generality,
that the special fibre of  $\Cal X_i$ has an irreducible component which is reduced. Such a component has a smooth
rational point over a finite extension of degree prime to $p$ of the residue field of $k$, as follows easily form
the Weil estimates for the number of rational points of a smooth, projective, and geometrically connected curve
over a finite field. Assertion {\rm (iv)} follows then from the fact that one can lift smooth
rational points (cf. [Grothendieck1], \'expos\'e III, Proposition 3.3). 
This shows {\rm (iii)} $\Rightarrow$ {\rm (iv)}, and hence {\rm (iv)} $\Longleftrightarrow$ {\rm (iii)}.

Next, we prove {\rm (v)}$\Longleftrightarrow${\rm (vi)}. 

This follows from the fact proved by Roquette, and
Lichtenbaum (cf. [Lichtenbaum]), that the kernel of the natural homomorphism
$\Br (k)\to \Br (X_i[s])$, which is induced by the natural morphism $X_i[s]\to \Spec k$, is a finite group whose
cardinality equals the index $\index (X_i[s])$ of $X_i[s]$. Plus the fact that the period divides the index.
The last assertion follows from the fact that $\index (X_i[s])$ divides $2\period (X_i[s])$ (cf. loc. cit).

Next, we prove {\rm (vi)}$\Rightarrow${\rm (iv)}. 

As a consequence of the Tate-Lichtenbaum duality between $\Pic(X_i[s])$
and $\Br(X_i[s])$, the natural homomorphism $\Br (X_i[s])\to \prod _{x\in X} \Br (k(x))$,
where $x$ runs over all closed points of $X_i[s]$, and $k(x)$ is the residue field at $x$, is injective
(cf. [Lichtenbaum], Proof of Theorem 5). 

Let $b\in \Br (k)$ be an element
of order $p$. Under our assumption on the kernel of the natural map $\Br(k)\to \Br (X_i[s])$ (cf. (vi)),
there exists a closed point $x\in X_i[s]$ such that the image of $b$ in $\Br (k(x))$, via the natural map $\Br(k)\to \Br(k(x)$ is non-zero.
After identifying
$\Br (k)$ and $\Br (k(x))$ with $\Bbb Q/\Bbb Z$, via the natural identification arising from local class field theory,
the natural map $\Br (k)\to \Br (k(x))$ is multiplication by the degree $\vert k(x):k \vert$ of the extension $k(x)/k$. From this follows immediately that
$\gcd (\vert k(x):k\vert,p)=1$. This shows {\rm (vi)}$\Rightarrow${\rm (iv)}.

Next we prove {\rm (iii)}$\Rightarrow${\rm (i)}.
As mentioned previously in the proof of Claim 1.6.4, similar arguments as the ones used in [Stix] 
yield the following result. The existence of the section $s_i$ implies that the period $\period(X_i[s])$, and index $\index(X_i[s])$, of $X_i[s]$, are divisible only by primes in 
$\{p\}\cup \{\Sigma'\}$. Assumption {\rm (iii)} implies that $\index (X_i[s])$ is prime to $p$, i.e. $\index (X_i[s])\in \Sigma'$. 
It then follows from the proof of Claim 1.6.4 that the section $s_i$ is good in the sense of Definition 1.4.1. 
Thus, {\rm (iii)}$\Rightarrow${\rm (i)}.

 Next, we prove {\rm (ii)}$\Rightarrow${\rm (v)}. 

The existence of the section $s_i$,
for a positive integer $i$,  implies that the period $\period (X_i[s])$, and index $\index (X_i[s])$, of $X_i[s]$, are divisible only by primes in 
$\{p\}\cup \{\Sigma'\}$ (cf. above discussion). 
If the section $s$ is well-behaved, then the Chern class
$c(s_i)$ lies in the Picard part $\Pic (X_i[s])^{\wedge,\Sigma}$ of $H^2(X_i[s],M_X)$. The image of the Chern class
$c(s_i)$ in  the group $H^2(X_i[s]\times _k \bar k,M_X)$, which is naturally isomorphic to $\hat \Bbb Z^{\Sigma}$,
via the natural restriction homomorphism $H^2(X_i[s],M_X)\to H^2(X_i[s]\times _k \bar k,M_X)$, equals 1. 

The image of the Chern class $c(s_i)$ in $H^2(X_i[s],M_X/pM_X)$ is a Picard element in the Picard part $\Pic (X_i[s])/p\Pic (X_i[s])$.
This Picard element has degree congruent to $1$ modulo $p$.
This shows that the period $\period (X_i[s])$, and index $\index (X_i[s])$, of $X_i[s]$ are prime to p. This shows {\rm (ii)}$\Rightarrow${\rm (v)}.

Thus, {\rm (i)}, {\rm (ii)}, {\rm (iii)}, {\rm (iv)}, {\rm (v)}, {\rm (vi)}, and {\rm (viii)}, are all equivalent.

The equivalence of {\rm (ii)} and {\rm (vii)} follows from Proposition 1.6.3 {\rm (ii)}, and the proof of Claim 1.6.4.

Finally, in the case where $\Sigma=\Primes$, the equivalence {\rm (vi)}$\Longleftrightarrow${\rm (viii)} is proven in
[Borovoi-Colliot-Th\'el\`ene-Skorobogatov], Theorem 2.5. Or, alternatively, the equivalence
{\rm (v)}$\Longleftrightarrow${\rm (viii)} follows from [Wittenberg], Theorem 3.2.1.
\qed
\enddemo

\definition{Definition 1.6.7 (Tame Point-Theoretic Sections over $p$-adic Local Fields)}\ \ \ 
Assume that $p\in \Sigma$.
We say that the section $s$ is tame point-theoretic if one of the equivalent
conditions  {\rm (i)}, {\rm (ii)}, {\rm (iii)}, {\rm (iv)}, {\rm (v)}, {\rm (vi)}, and {\rm (vii)}, 
in Proposition 1.6.6 is satisfied. 

We say that the section $s$ is uniformly tame point-theoretic, if for every finite extension $k'/k$
the induced section $s_{k'}:G_{k'}\to \Pi_{X_{k'}}$ is tame point-theoretic.
\enddefinition

Definition 1.6.7 is motivated by the property {\rm (iii)} in Proposition 1.6.6.

In fact a tame point-theoretic section $s$ is necessarily uniformly tame point-theoretic. More precisely, we have
the following. 

\proclaim{Proposition 1.6.8} Assume that $k$ is a $p$-adic local field, and $p\in \Sigma$. Suppose that
the section $s$ is tame point-theoretic (in the sense of Definition 1.6.7). Then $s$ is uniformly tame point-theoretic. In particular, if
the section $s$ is good then $s$ is uniformly good.
\endproclaim

\demo{Proof} This follows, in the case where $\Sigma=\Primes$,  from  Corollary 3.2.3 in [Wittenberg], together with Proposition 1.6.6,
using the property {\rm (viii)} in Proposition 1.6.6 defining a good section. 
We give another different proof with no restrictions on $\Sigma$. 

It suffices to show the following.
Given a section $s:G_{k}\to \Pi_{X}$ of $\Pi_X$, such that the image of $\Pic(X)$ in $H^2(G_k,M_X)$ via the homomorphism
$\Pic(X)\to H^2(G_k,M_X)$, induced by pulling back
Chern classes of line bundles via the section $s$, is zero, then for any finite extension $k'/k$ 
the image of $\Pic(X_{k'})$ in $H^2(G_{k'},M_X)$ via the homomorphism $\Pic(X_{k'})\to H^2(G_{k'},M_X)$,
induced by pulling back Chern classes of line bundles via
the induced section $s_{k'}:G_{k'}\to \Pi_{X_{k'}}$, is also zero. 

We have a natural commutative diagram:
$$
\CD
\Pic (X_{k'}) @>>> H^2(G_{k'},M_X)\\
@AAA           @A{\res}AA\\
\Pic (X) @>>> H^2(G_{k},M_X)
\endCD
$$
where the horizontal maps are the above ones, the left vertical map is the pull back homomorphism of line bundles, and
the right vertical map is the restriction homomorphism. 

Note that  $H^2(G_{k},M_X)$ is naturally identified with
$\hat \Bbb Z^{\Sigma}$, and via this identification the right vertical map in the above diagram is multiplication by the degree
$\vert k':k \vert$ of the extension $k'/k$. 

We can assume, without loss of generality, that $k'$ is Galois over $k$,
with Galois group $H\defeq \Gal(k'/k)$. The natural action of $H$ on $H^2(G_{k'},M_X)$ is trivial.
In particular, the above homomorphism $\Pic (X_{k'})\to H^2(G_{k'},M_X)$ factors as  $\Pic (X_{k'})\to \Pic (X_{k'})_H
\to H^2(G_{k'},M_X)$, where $\Pic (X_{k'})_H$ denotes the co-invariant group.

The image of $\Pic(X)$ in  $H^2(G_{k'},M_X)$ is trivial by assumption. Also the image of $\Pic(X)$ in $\Pic (X_{k'})^H$,
via the natural homomorphism $\Pic(X)\to \Pic (X_{k'})^H$, is a finite index subgroup 
(here $\Pic (X_{k'})^H$ denotes the invariant subgroup), and the image of
$\Pic (X_{k'})^H$ in $\Pic (X_{k'})_H$ via the natural map $\Pic (X_{k'})^H\to \Pic (X_{k'})_H$ is a finite index
subgroup.

From this we deduce that the above image of $\Pic (X_{k'})_H$ in $H^2(G_{k'},M_X)$ is torsion, hence equals $0$, since
$H^2(G_{k'},M_X)$ is torsion-free.
\qed
\enddemo

\subhead {1.7}
\endsubhead
Next, we will generalise the notion of a tame
point-theoretic section of arithmetic fundamental groups to sections of certain quotients
of the absolute Galois group of $K_X$. 

Let 
$$G\defeq G_{K_X}\defeq \Gal (K_X^{\sep}/K_X)$$
be the absolute Galois group of the function field $K_X$ of $X$, and 
$$\overline G\defeq \Gal (K_X^{\sep}/K_X.\bar k)$$
the absolute Galois group of the function field $K_X.\bar k$ of $X\times _k\bar k$. Then $G$ sits naturally
in the following exact sequence
$$1\to \overline G\to G\to G_k\to 1.\tag {$1.9$}$$

Let $\overline G \twoheadrightarrow \overline H$
be a continuous surjective homomorphism between profinite groups,
where $\overline H$ is a characteristic quotient of $\overline G$,
which inserts into the following sequence of continuous surjective homomorphisms
$$\overline G \twoheadrightarrow \overline H\twoheadrightarrow \Delta_X\defeq \pi_1(X\times _k \bar k, \bar \eta)
^{\Sigma}$$
Here $\Sigma \subseteq \Primes$ is a non-empty set of primes.

By pushing out the exact sequence (1.9) by the homomorphism $\overline G \twoheadrightarrow \overline H$ we obtain
an exact sequence of profinite groups
$$1\to \overline H \to H \to G_k \to 1,\tag{$1.10$}$$
which inserts into the following commutative diagram:

$$
\CD
1 @>>> \overline G      @>>>     G  @>>> G_k  @>>> 1\\
  @.        @VVV           @VVV              @V{\id}VV \\
1@>>> \overline H      @>>>     H  @>>> G_k    @>>> 1 \\
 @.        @VVV           @VVV              @V{\id}VV \\
1@>>> \Delta_X    @>>>   \Pi_X @>>> G_k @>>> 1
\endCD
$$
where the horizontal sequences are exact, and  the middle and left vertical homomorphisms are surjective. 

Let
$$\tilde s:G_k\to H$$
be a group-theoretic section of the natural projection $H\twoheadrightarrow G_k$. 

Assume that the closed subgroup
$\overline H$ of $H$ is topologically finitely generated. Then one can define a system of neighbourhoods
$\{H_i[\tilde s]\}_{i\ge 1}$
of the section $\tilde s:G_k\to H$ in a similar way as in 1.3.

In particular, this system of neighbourhoods corresponds to a tower of finite (possibly ramified) covers
$$...\to \Tilde X_{i+1}[\tilde s]\to \Tilde X_{i}[\tilde s]\to ...\to \Tilde X_{1}[\tilde s]\defeq X$$
defined over $k$. We will refer to the set $\{\Tilde X_i[\tilde s]\}_{i\ge 1}$
as a system of neighbourhoods of the section $\tilde s$.

\proclaim {Definition/Lemma 1.7.1} Assume that $k$ is a $p$-adic local field, and $p\in \Sigma$. We say that the 
group-theoretic section $\tilde s:G_k\to H$
is a tame point-theoretic section if the following equivalent conditions are satisfied.

\noindent
{\rm (i)}\ For every positive integer $i$, one has $\Tilde X_i[\tilde s](k^{\tame})\neq \varnothing$.

\noindent
{\rm (ii)}\ For every positive integer $i$, there exists a finite extension $\ell_i$ of $k$ with
$\gcd ([\ell_i:k],p)=1$ such that  $\Tilde X_i[\tilde s](\ell_i)\neq \varnothing$.

\noindent
{\rm iiii)}\ For every positive integer $i$, the $p$-primary part of the 
kernel of the natural homomorphism  $\Br (k)\to \Br(\Tilde X_i[\tilde s])$ is trivial.

\noindent
{\rm (iv)}\ Let $\Tilde L/K_X$ be the sub-extension of $K_X^{\sep}/K_X$ 
with Galois group $H$, and $L/K$ the sub-extension of $\Tilde L/K$
corresponding to the closed subgroup $\tilde s(G_k)$
of $H$. Then the $p$-primary part of the kernel of the natural homomorphism $\Br(k)\to \Br(L)$ is trivial.

\endproclaim

\demo{Proof} The equivalence of the above properties is proven in a similar way as in the proof of Proposition 1.6.6.
\qed
\enddemo

\definition{Remark 1.7.2} Property (iv) in Definition/Lemma 1.7.1 is independent of the choice of a system of neigbourhoods
of the section $\tilde s$. Using this property one can define the notion of  
a tame point-theoretic section  $\tilde s:G_k\to H$ without any assumption on  $\overline H$.
\enddefinition

\subhead {1.8}
\endsubhead
For the rest of this section we will assume that $k$ is a number field, i.e. $k$ is a finite extension of the
field $\Bbb Q$ of rational numbers. 

Let $v$ be a place of $k$, and denote by $k_v$ the completion of $k$ at $v$.
Let $X_v\defeq X\times _k k_v$.
Let $D_v\subseteq G_k$ be a decomposition group at $v$ ($D_v$ is only defined up to conjugation), which is naturally
isomorphic to the absolute
Galois group $G_{k_v}$ of $k_v$. By pulling back the exact sequence
$$1\to \Delta _X\to \Pi_X\to G_k\to 1,$$
by the natural homomorphism $D_v\to G_k$, we obtain the exact sequence
$$1\to \Delta _{X_v}\to \Pi_{X_v}\to G_{k_v}\to 1.$$
Note that there exists a natural isomorphism $\Delta _X\isom \Delta _{X_v}$.

In particular, the section $s:G_k\to \Pi_X$ induces naturally a section 
$$s_v: G_{k_v}\to \Pi_{X_v}$$ 
of the arithmetic pro-$\Sigma$ fundamental group $\Pi_{X_v}$, for each place $v$ of $k$.

\proclaim {Proposition 1.8.1} We use the same notations as above. Consider the following properties.

\noindent
{\rm (i)} The section $s$ is a uniformly good group-theoretic section.

\noindent
{\rm (ii)} for each place $v$ of $k$, the section $s_v$ is a good group-theoretic section.

Then we have the equivalence
$(i)\Longleftrightarrow (ii)$.

\endproclaim

\demo{Proof}
First, we prove {\rm (ii)}$\Rightarrow$ {\rm (i)}.

For every positive integer $n$, we have the following commutative diagram:
$$
\CD
\Pic (X)/n\Pic (X) @>>>    H^2(X,M_X/nM_X)     @>{s_n^{\star}}>>     H^2(G_k,M_X/nM_X)  \\
@VVV        @VVV           @VVV \\
\prod_v \Pic(X_v)/n\Pic(X_v)  @>>>   \prod _v H^2(X_v,M_X/nM_X)   @>{s_{v,n}^{\star}}>> 
\prod _v H^2(G_{k_v},M_X/nM_X)
\endCD
$$
where the product is over the set of all places of $k$, the left horizontal homomorphisms are induced
by Kummer theory, the right horizontal homomorphisms are the restrictions via the sections $s$, and $s_v$,
and the vertical homomorphisms are the natural diagonal ones. 

The implication (ii)$\Rightarrow $(i) follows from the fact that
the far right vertical diagonal homomorphism $  H^2(G_k,M_X/nM_X)\to \prod _v H^2(G_{k_v},M_X/nM_X)$ is injective,
as follows from the Brauer-Hasse-Noether principle in global class field theory.

Next, we prove {\rm (i)}$\Rightarrow${\rm (ii)}. 

The case where $v$ is a real place follows from the well-known fact that every group-theoretic section $s_v$ arises from a rational point
$x\in X(k_v)$ in this case, the so-called real section conjecture (cf. [Stix], A, for example). So we only consider the case where $v$ is a $p$-adic place,
$p>0$ being a prime integer.

First, we treat the degree $0$ line bundles.
The group $\Pic ^0(X_v)$ is an open subgroup of $\Jac_{X_v}(k_v)$, the later being a finitely generated
profinite group. Further, we have $\Jac_{X_v}(k_v)\isom \Bbb Z_p^r\times A$, 
where $A$ is a finite abelian group. We can, without loss of generality, assume that $\Pic ^0(X_v)=\Jac_{X_v}(k_v)$.

Write $k_v^h$ for the henselisation of $k$ at $v$.
Since $\Jac_{X_v}(k_v^h)$ is dense in $\Jac_{X_v}(k_v)$, the elements of 
$\Pic^0(X_v)$ can be approximated by elements in
$\Jac_X(k')$, where $k'/k$ is a finite extension. So the assertion that the images of the
Chern classes of the elements of  $\Pic ^0(X_v)$ in  $H^2(G_{k_v},M_X/nM_X)$, via the map $s_{v,n}^{\star}$, is zero follows. 

For the
non-zero degree part it suffices to consider a line bundle of a given non-zero degree. In this case one
can consider the canonical bundle $\Omega_v\in \Pic (X_v)$ of $X_v$ of degree $2g-2$, which 
arises from the canonical bundle $\Omega\in \Pic (X)$ of $X$.
\qed
\enddemo

\subhead
\S 2. Cuspidalisation of Good Sections of Arithmetic Fundamental Groups over Slim Fields
\endsubhead
In this section we will introduce, and investigate,  the problem of cuspidalisation of group-theoretic sections of arithmetic
fundamental groups. 

We follow the notations in $\S 1$. 

In particular, $X$ is a proper, smooth, geometrically connected, hyperbolic algebraic curve over the field $k$,
$\Sigma \subseteq \Primes$ is
a non-empty set of prime integers, with $\char (k)\notin \Sigma$, and we have the natural exact sequence
$$1\to \Delta_X\to \Pi_X @>{\pr_{X,\Sigma}}>> G_k\to 1,\tag {$2.1$}$$
where $\Pi_X$ is the geometrically pro-$\Sigma$ arithmetic fundamental group of $X$ (cf. exact sequence (1.2)).

\subhead {2.1}
\endsubhead
In this subsection we recall the definition of (geometrically) cuspidally central, and cuspidally abelian, 
arithmetic fundamental groups of affine hyperbolic curves,
and the definition of cupidally abelian absolute Galois groups of function fields of curves (cf. [Mochizuki], Definition 1.5).

\subhead {2.1.1}
\endsubhead
Let $U\subseteq X$ be a non-empty open subscheme of $X$. The geometric point $\eta$ of $X$ (cf. 1.1) determines
a geometric point $\eta$ of $U$, and a geometric point $\bar \eta$ of $\overline U\defeq U\times _k \bar k$.

Write 
$$\Delta_U\defeq \pi_1(\overline U,\bar \eta)^{\Sigma}$$
for the maximal pro-$\Sigma$ quotient of the fundamental group $\pi_1(\overline U,\bar \eta)$ of $\overline U$
with base point $\bar \eta$, and 
$$\Pi_U\defeq  \pi_1(U, \eta)/ \Ker  (\pi_1(\overline U,\bar \eta)\twoheadrightarrow
\pi_1(\overline U,\bar \eta)^{\Sigma})$$ 
for the quotient of  the arithmetic fundamental group
$\pi_1(U, \eta)$ by the kernel of the natural surjective homomorphism
$\pi_1(\overline U,\bar \eta)\twoheadrightarrow \pi_1(\overline U,\bar \eta)^{\Sigma}$, which is a normal subgroup
of $\pi _1(U,\eta)$. 

Thus, we have a natural exact sequence
$$1\to \Delta_U\to \Pi_U @>{\pr_{U,\Sigma}}>> G_k\to 1,$$
which inserts into the following commutative diagram:
$$
\CD
1 @>>> \Delta_U  @>>> \Pi_U  @>{\pr_{U,\Sigma}}>> G_k @>>> 1  \\
@.   @VVV     @VVV    @V{\id }VV   \\
1 @>>> \Delta_X   @>>> \Pi_X    @>{\pr_{X,\Sigma}}>> G_k  @>>> 1\\
\endCD
$$
where the left vertical homomorphisms are surjective, and are induced by the natural surjective homomorphism
$\pi_1(\overline U,\bar \eta)\twoheadrightarrow \pi_1(\overline X,\bar \eta)$.

Let
$$I_U\defeq \Ker (\Pi_U\twoheadrightarrow \Pi_X)=\Ker (\Delta_U\twoheadrightarrow \Delta_X).$$
We shall refer to $I_U$ as the cuspidal subgroup of $\Pi_U$ (cf. [Mochizuki], Definition 1.5). It is the normal
subgroup of $\Pi_U$ generated by the (pro-$\Sigma$) inertia subgroups at the geometric points of
$S\defeq X\setminus U$.
We have the following natural exact sequence
$$1\to I_U\to \Pi_U\to \Pi_X\to 1.\tag {$2.2$}$$

Let $I_U^{\ab}$ be the maximal abelian quotient of $I_U$. By pushing out the exact sequence $(2.2)$ by the natural
surjective homomorphism $I_U\twoheadrightarrow I_U^{\ab}$ we obtain a natural commutative diagram:
$$
\CD
1 @>>> I_U @>>> \Pi_U @>>> \Pi_X @>>> 1  \\
@. @VVV   @VVV   @V{\id}VV  \\
1 @>>> I_U^{\ab} @>>> \Pi_U^{\c-\ab} @>>> \Pi_X @>>> 1
\endCD
$$

We shall refer to the quotient $\Pi_U^{\c-\ab}$ of $\Pi_U$ as the maximal cuspidally abelian quotient of $\Pi_U$,
with respect to the natural homomorphism $\Pi_U\twoheadrightarrow \Pi_X$ (cf. [Mochizuki], Definition 1.5).

Similarly, we can define the maximal cuspidally abelian quotient $\Delta_U^{\c-\ab}$ of $\Delta_U$,
with respect to the natural homomorphism $\Delta_U \twoheadrightarrow \Delta_X$, which sits in a natural exact sequence
$$1 \to I_U^{\ab} \to \Delta_U^{\c-\ab} \to \Delta_X \to 1. \tag{$2.3$}$$

We have a commutative diagram:
$$
\CD
1 @>>> I_U @>>> \Delta_U @>>> \Pi_X @>>> 1  \\
@. @VVV   @VVV   @V{\id}VV  \\
1 @>>> I_U^{\ab} @>>> \Delta_U^{\c-\ab} @>>> \Pi_X @>>> 1
\endCD
$$
which is a push out diagram by the natural surjective homomorphism $I_U\twoheadrightarrow I_U^{\ab}$.

The profinite group $\Delta _X$ acts naturally by automorphisms on $I_U^{\ab}$ (cf. exact sequence (2.3)).
Write $I_U^{\cn}$ for the maximal quotient of  $I_U^{\ab}$ on which the action of $\Delta_X$ 
is trivial. By pushing out the sequence $(2.3)$
by the natural surjective homomorphism $I_U^{\ab}\twoheadrightarrow I_U^{\cn}$ we obtain a natural exact sequence
$$1\to I_U^{\cn}\to \Delta_U^{\c-\cn}\to \Delta _X\to 1.\tag {$2.4$}$$

Define 
$$\Pi_U^{\c-\cn}\defeq \Pi_U^{\c-\ab}/\Ker (I_U^{\ab}\twoheadrightarrow I_U^{\cn}),$$
which sits naturally in the following exact sequence
$$1\to I_U^{\cn}\to \Pi _U^{\c-\cn}\to \Pi  _X\to 1.\tag {$2.5$}$$

We shall refer to the quotient $\Pi_U^{\c-\cn}$ of $\Pi_U$ as the maximal (geometrically) cuspidally central
quotient of $\Pi_U$, with respect to the natural homomorphism $\Pi_U\twoheadrightarrow \Pi_X$ (cf. loc. cit.).
We have the following commutative diagram
$$
\CD
1 @>>> I_U @>>> \Pi_U @>>> \Pi_X @>>> 1\\
@. @VVV     @VVV    @V{\id}VV \\
1 @>>> I_U^{\ab} @>>> \Pi_U^{\c-\ab} @>>> \Pi_X @>>> 1\\
@. @VVV     @VVV    @V{\id}VV \\
1@>>>  I_U^{\cn}  @>>>  \Pi _U^{\c-\cn}  @>>> \Pi  _X @>>> 1\\
\endCD
$$
where the vertical homomorphisms on the left, and in the middle, are natural surjections.

\subhead {2.1.2}
\endsubhead
Similarly, we have a natural exact sequence of absolute Galois groups
$$1\to G_{\bar k.K_X}\to G_{K_X}\to G_k\to 1,$$
where $G_{\bar k.K_X}\defeq \Gal (K_X^{\sep}/\bar k.K_X)$, and $G_{K_X}\defeq \Gal (K_X^{\sep}/K_X)$.

Let 
$$\overline G_X\defeq G_{\bar k.K_X}^{\Sigma}$$ 
be the maximal pro-$\Sigma$ quotient of $G_{\bar k.K_X}$, and 
$$G_X\defeq  G_{K_X}/ \Ker  ( G_{\bar k.K_X}  \twoheadrightarrow G_{\bar k.K_X}^{\Sigma}),$$
which insert into the following commutative diagram:
$$
\CD
1 @>>>  G_{\bar k.K_X} @>>>  G_{K_X}   @>>>  G_k   @>>>  1 \\
@. @VVV   @VVV    @V{\id}VV   \\
1 @>>> \overline G_{X}  @>>>  G_X  @>{\Tilde {\pr}_{X,\Sigma}}>> G_k @>>> 1  \\
@.   @VVV     @VVV    @V{\id }VV   \\
1 @>>> \Delta_X   @>>> \Pi_X    @>{\pr_{X,\Sigma}}>> G_k   @>>>  1\\
\endCD
$$
where the left, and middle, vertical maps are the natural surjective homomorphisms. 

Let
$$\Tilde I_X\defeq \Ker (G_{X}\twoheadrightarrow \Pi_X)=\Ker (\overline G_{X}\twoheadrightarrow \Delta_X).$$
We shall refer to $\Tilde I_X$ as the cuspidal subgroup of $G_{X}$. It is the normal subgroup of
$G_{X}$  generated by the (pro-$\Sigma$) inertia subgroups at all geometric closed points of $X$.
We have the following natural exact sequence
$$1\to \Tilde I_X\to G_{X}  \to \Pi_X\to 1.$$

Let $\Tilde I_X^{\ab}$ be the maximal abelian quotient of $\Tilde I_X$. By pushing out the above sequence by the natural
surjective homomorphism $\Tilde I_X\twoheadrightarrow \Tilde I_X^{\ab}$, we obtain a natural exact sequence
$$1\to \Tilde I_X^{\ab}\to G_{X}^{\c-\ab}\to \Pi_X\to 1.\tag {$2.6$}$$

We will refer to the quotient $G_{X} ^{\c-\ab}$ as the maximal cuspidally abelian quotient of $ G_{X}$,
with respect to the natural surjective homomorphism $G_{X}\twoheadrightarrow \Pi_X$.
Note that $ G_{X} ^{\c-\ab}$ is naturally identified with the projective limit
$$\underset{U}\to {\varprojlim}\ \Pi_U^{\c-\ab},$$ 
where the limit runs over all open subschemes $U$ of $X$.

\subhead {2.2}
\endsubhead
Next, we will consider a continuous group-theoretic section $s:G_k\to \Pi_X$
of the natural projection $\pr_X\defeq \pr _{X,\Sigma}:\Pi_X\twoheadrightarrow G_k$ (cf. 1.1).

\definition{Definition 2.2.1 (Lifting of Group-Theoretic Sections)}
Let $U\subseteq X$ be a non-empty open subscheme. We say that a continuous group-theoretic
section $s_U:G_k\to \Pi_U$ of the natural projection  $\pr_U\defeq \pr_{U,\Sigma}:\Pi_U\twoheadrightarrow G_k$,
meaning that $\pr_U\circ s_U=\id _{G_k}$, is a lifting
of the section $s:G_k\to \Pi_X$, if $s_U$ fits into a commutative diagram:
$$
\CD
G_k @>{s_U}>> \Pi_U \\
@V{\id}VV   @VVV \\
G_k @>{s}>>  \Pi _X
\endCD
$$
where the right vertical homomorphism is the natural one. 

More generally, we say that a group-theoretic
section $\tilde s :G_k\to G_{X}$ of the natural projection  $\tilde {\pr_X}\defeq \tilde {\pr}_{X,\Sigma}:G_{X}\twoheadrightarrow G_k$,
meaning that $\tilde {\pr_X}\circ \tilde s=\id _{G_k}$, is a lifting of
the section $s$, if $\tilde s$ fits into a commutative diagram:
$$
\CD
G_k @>{\tilde s}>> G_{X} \\
@V{\id}VV   @VVV \\
G_k @>{s}>>  \Pi _X
\endCD
$$
where the right vertical homomorphism is the natural one.
\enddefinition

\definition {Remark 2.2.2}
One can easily verify that if the section $s$ is point-theoretic, i.e. $s=s_x$ arises from a rational point $x\in X(k)$ (cf. Defintion 3.1.1), 
then the section $s$ can be lifted to a section
$s_U:G_k\to \Pi_U$ of the natural projection $\Pi_U\twoheadrightarrow G_k$, for every open subscheme $U\subseteq X$, and can also be lifted to a section
$\tilde s : G_k\to G_{X}$ of the natural projection $G_X\twoheadrightarrow G_k$.
\enddefinition

Next, we introduce the cuspidalisation problem for sections of arithmetic fundamental groups.

\definition {The Cuspidalisation Problem for sections of Arithmetic Fundamental Groups}
Given a group-theoretic section $s:G_k\to \Pi_X$ as above, and a  non-empty open subscheme $U\subseteq X$,
is it possible to construct a lifting $s_U:G_k\to \Pi_U$ of $s$? Similarly, is it possible to construct a lifting $\tilde s : G_k\to G_{X}$ of $s$?
\enddefinition

\subhead {2.3}
\endsubhead
In this section we will investigate the problem of cuspidalisation of a group-theoretic section $s:G_k\to \Pi_X$
of the geometrically pro-$\Sigma$ arithmetic fundamental group $\Pi_X$, under the assumption
that the section $s$ is uniformly good (cf. Definition 1.4.1).  

First, we recall the definition of a slim field, and define the notion of a $\Sigma$-regular field.

\definition {Definition 2.3.1}
(i)\ We say that the field $k$ is slim, if its absolute Galois group $G_k$ is slim in the sense
of [Mochizuki], $\S0$, meaning that every open subgroup of $G_k$ is centre free. 

Examples of slim fields include number
fields, and $p$-adic local fields (cf. [Mochizuki1], Theorem 1.1.1). 

One defines in a similar way the notion of a slim profinite group $G$, meaning that
every open subgroup of $G$ is centre free.

(ii)\ We say that the field $k$ is $\Sigma$-regular, if for every prime integer $l\in \Sigma$, and every finite extension $k'/k$,
the $l$-part of the cyclotomic character $\chi _l:G_{k'}\to \Bbb Z_l^{\times}$ is not trivial; or equivalently, if for every prime integer $l\in \Sigma$
the image of the $l$-part of the cyclotomic character $\chi _l:G_k\to \Bbb Z_l^{\times}$ is infinite. 

Examples of $\Sigma$-regular fields (for every non-empty set
$\Sigma$ of prime integers) include number fields, $p$-adic local fields, and finite fields. 

The field $k$ is $\Sigma$-regular if and only if, for every finite extension $k'/k$, the $G_{k'}$-module
$M_X$ (cf. 1.2) has no non trivial fixed elements.

In the case where $\Sigma=\Primes$, and $k$ is $\Sigma$-regular, we say that $k$ is regular.

\enddefinition

Our first result concerning the cuspidalisation problem is the following.

\proclaim {Proposition 2.4 (Lifting of Uniformly Good Sections to Cuspidally Central Arithmetic Fundamental Groups
over Slim Fields)}
Assume that the field $k$ is slim, and $\Sigma$-regular (cf. Definition 2.3.1).
Let $s:G_k\to \Pi_X$ be a uniformly good section (in the sense of Definition 1.4.1) of the natural projection $\Pi_X\twoheadrightarrow G_k$.
Then the followings hold.

(i)\ 
Let $U\defeq X\setminus S$ be a non-empty open subscheme of $X$, and
$\Pi_{U}^{\c-\cn}$ the maximal (geometrically) cuspidally central quotient of $\Pi_{U}$, with respect to the natural homomorphism
$\Pi_{U}\twoheadrightarrow \Pi_X$ (cf. 2.1.1). Then there exists a section $s_U^{\c-\cn}:G_k\to \Pi_{U}^{\c-\cn}$ of the
natural projection  $\Pi_{U}^{\c-\cn}\twoheadrightarrow G_k$, which lifts the section $s$, i.e. which inserts into the following
commutative diagram:
$$
\CD
G_k @>s_U^{\c-\cn}>>  \Pi_{U}^{\c-\cn} \\
@V{\id}VV     @VVV  \\
G_k   @>{s}>> \Pi_X
\endCD
$$

Moreover, the set of all possible liftings $s_U^{\c-\cn}$ of $s$ is a torsor under the group $H^1(G_k,I_U^{\c-\cn})$.
Here the $G_k$-module structure of $I_U^{\c-\cn}$ is naturally induced by the section $s$.

(ii)\ There exists, for each non-empty open subscheme $U\defeq X\setminus S$ of $X$, a section
$s_U^{\c-\cn}:G_k\to \Pi_{U}^{\c-\cn}$ as in (i) (i.e. $s_U^{\c-\cn}$ is a lifting of $s$),
such that for every non-empty open subscheme
$V\defeq X\setminus T$ of $X$, with $U\subseteq V$, we have the following commutative diagram:
$$
\CD
G_k @>s_U^{\c-\cn}>>  \Pi_{U}^{\c-\cn} \\
@V{\id}VV     @VVV  \\
G_k   @>s_V^{\c-\cn}>> \Pi_{V}^{\c-\cn}\\
\endCD
$$
where the right vertical homomorphism is the natural one.

(iii)\ There exists a section $\tilde s^{\c-\cn}:G_k\to \underset{U}\to {\varprojlim}\ \Pi_U^{\c-\cn}$
of the natural projection $\underset{U}\to {\varprojlim}\ \Pi_U^{\c-\cn}\twoheadrightarrow G_k$, which lifts the section $s$,
i.e. which inserts in the following commutative diagram:
$$
\CD
G_k @>\tilde s^{\c-\cn}>>   \underset{U}\to {\varprojlim}\ \Pi_{U}^{\c-\cn} \\
@V{\id}VV     @VVV  \\
G_k   @>{s}>> \Pi_X
\endCD
$$
where the right vertical homomorphism is the natural one. Here the projective limit is over all open subschemes $U$ of $X$
\endproclaim

\demo{Proof}
We start by proving assertion (i).

First, we treat the case where the set $S=\{x_i\}_{i=1}^n\subseteq X(k)$ consists of finitely many $k$-rational points.
 
 For $i\in \{1,...,n\}$, write $U_{x_i}\defeq X\setminus \{x_i\}$. The maximal (geometrically) cuspidally central quotient
$\Pi_{U_{x_i}}^{\c-\cn}$ of $\Pi_{U_{x_i}}$,
with respect to the natural homomorphism
$\Pi_{U_{x_i}}\twoheadrightarrow \Pi_X$, sits naturally in the following exact sequence
$$1\to M_X\to \Pi_{U_{x_i}}^{\c-\cn} \to \Pi_X\to 1\tag {$2.7$}$$
(cf. [Mochizuki], Proposition 1.8, (iii), and [Mochizuki3], Lemma 4.2). 

By pulling back the group extension $(2.7)$ by the section $s:G_k\to \Pi_X$ we
obtain a group extension
$$1\to M_X\to s^{\star}(\Pi_{U_{x_i}}^{\c-\cn}) \to G_k\to 1,\tag {$2.8$}$$
which inserts naturally in the following commutative diagram:
$$
\CD
1 @>>> M_X  @>>>  D_{x_i}\defeq s^{\star}(\Pi_{U_{x_i}}^{\c-\cn}) @>>> G_k  @>>> 1 \\
@.    @V{\id}VV  @VVV    @V{s}VV   \\
1 @>>> M_X  @>>>  \Pi_{U_{x_i}}^{\c-\cn} @>>> \Pi_X  @>>> 1
\endCD
$$
where the right square is cartesian. 

The class in $H^2(\Pi_X,M_X)$ of the group
extension $(2.7)$ coincides, via the natural identification  $H^2(\Pi_X,M_X)\isom H^2(X,M_X)$ (cf. [Mochizuki], Proposition 1.1),
with the \'etale Chern class $c(x_i)\in H^2(\Pi_X,M_X)$ associated to the degree $1$ line bundle $\Cal O(x_i)$ (cf. [Mochizuki3], Lemma 4.2).
The class in $H^2(G_k,M_X)$ of the group
extension $(2.8)$ coincides then with the image $s^{\star} (c(x_i))$ of the Chern class $c(x_i)$ via the
(restriction) homomorphism
$s^{\star}:H^2(X,M_X)\isom H^2(\Pi_X,M_X)\to H^2(G_k,M_X)$, which is naturally induced by $s$. This image equals $0$ since the
section $s$ is assumed to be good. This follows from the very definition of goodness (cf. Definition 1.4.1).

Thus, the group extension $(2.8)$
splits. The set of all possible splittings of the extension $(2.8)$ is a torsor under the group $H^1(G_k,M_X)$, which is naturally
identified, via Kummer theory, with the $\Sigma$-adic completion ${(k^{\times})}^{\wedge,\Sigma}\defeq \underset{n\ \Sigma-\text {integer}}\to
{\varprojlim}\ k^{\times}/(k^{\times})^n$ of the multiplicative group $k^{\times}$ of $k$.

Each splitting of the exact sequence $(2.8)$ gives rise naturally to a group-theoretic section
$s_{U_{x_i}}^{\c-\cn}:G_k\to \Pi_{U_{x_i}}^{\c-\cn}$ of the natural projection $\Pi_{U_{x_i}}^{\c-\cn}\twoheadrightarrow G_k$
(cf. above commutative diagram) which necessarily lies above $s$.
Reciprocally, each section  $s_{U_{x_i}}^{\c-\cn}:G_k\to \Pi_{U_{x_i}}^{\c-\cn}$ of the
natural projection $\Pi_{U_{x_i}}^{\c-\cn}\twoheadrightarrow G_k$, which lifts the section $s$, arises from a splitting of the group extension $(2.8)$.
We have the following commutative diagram:
$$
\CD
G_k @>{s_{U_{x_i}}^{\c-\cn}}>> \Pi_{U_{x_i}}^{\c-\cn}\\
@V{\id}VV    @VVV \\
G_k @>{s}>>  \Pi_X
\endCD
$$
where the right vertical homomorphism is the natural surjective one. Moreover, all the possible sections
$s_{U_{x_i}}^{\c-\cn}$ as above, which lift the section $s$, form a torsor under the group $H^1(G_k,M_X)$. 

We have a natural identification
$\Pi_{U}^{\c-\cn}\isom \prod _{i=1}^n\Pi_{U_{x_i}}^{\c-\cn}$, where the fibre product is taken over $\Pi_X$.
More precisely, $\Pi_{U}^{\c-\cn}$ is naturally an extension of $\Pi_X$
by a product of copies of $M_X$ indexed by the points of $S$, i.e. sits in a natural exact sequence 
$$1\to \prod_{x\in S} M_X\to \Pi_{U}^{\c-\cn} \to \Pi_X\to 1$$
(cf. [Mochizuki], Proposition 1.8, (iii)). 

In particular, any collection of sections $\{s_{U_{x_i}}^{\c-\cn}:
G_k\to \Pi_{U_{x_i}}^{\c-\cn}\}_{i=1}^n$, which lift the section $s$, determine naturally a section
$$s_{U}^{\c-\cn}:G_k\to \Pi_{U}^{\c-\cn}$$
of the natural projection $\Pi_{U}^{\c-\cn}\twoheadrightarrow G_k$, which lifts the section
$s$. The set of all possible sections $s_U^{\c-\cn}$ as above, which lift the section $s$, is a torsor under
$H^1(G_k,\prod _{x\in S}M_X)$, where the product is taken over all points $x$ in $S$, by the above discussion.

Assume now that $S\defeq \{p_i\}_{i=1}^n$ consists of a finite set of closed points, which are not necessarily $k$-rational.
Let $k_S$ be the minimal Galois extension of $k$, with Galois group $\Gal (k_S/k)$, over which all points in $S$
are rational. For $i\in \{1,...,n\}$, let $S_i\defeq \{x_{i,j}\}_{j=1}^{n_i}$ be the set of points of $X_{k_S}\defeq X\times _k k_S$ above $p_i$.
Thus, $x_{i,j}\in X_{k_S}(k_S)$ is a $k_S$-rational point of  $X_{k_S}$, and $n_i\defeq \vert k(p_i):k\vert$, where $k(p_i)$ is the residue field at $p_i$. 

Let
$s_{k_S}:G_{k_S}\to \Pi_{X_{k_S}}$ be the group-theoretic section of the natural projection
$\Pi_{X_{k_S}}\twoheadrightarrow G_{k_S}$, which is naturally induced by the section $s$.

For $i\in \{1,...,n\}$, $j\in \{1,...,n_i\}$,
let $\Tilde U_{x_{i,j}}\defeq X_{k_S}\setminus \{x_{i,j}\}$, and
$\Tilde U \defeq \prod _{p_i\in S}  (\prod _{x_{i,j}\in S_i} \Tilde U_{x_{i,j}})$,
where the fibre product is taken above $X_{k_S}$.

The maximal (geometrically) cuspidally central quotient $\Pi_{\Tilde U_{x_{i,j}}}^{\c-\cn}$ of  $\Pi_{\Tilde U_{x_{i,j}}}$,
with respect to the natural surjective homomorphism $\Pi_{\Tilde U_{x_{i,j}}}\twoheadrightarrow  \Pi_{X_{k_S}}$, sits naturally
in the exact sequence:
$$1\to M_X\to \Pi_{\Tilde U_{x_{i,j}}}^{\c-\cn}\to \Pi_{X_{k_S}}\to 1.$$

Let  $\Pi_{\Tilde U}^{\c-\cn}$ be the  maximal (geometrically) cuspidally central quotient of  $\Pi_{\Tilde U}$,
with respect to the natural surjective homomorphism $\Pi_{\Tilde U}\twoheadrightarrow  \Pi_{X_{k_S}}$,
which is naturally identified with the fibre product
$\prod _{p_i\in S}(\prod _{x_{i,j}\in S_i}\Pi_{\Tilde U_{x_{i,j}}}^{\c-\cn})$ over $\Pi_{X_{k_S}}$.
Thus, $\Pi_{\Tilde U}^{\c-\cn}$ is an extension of $\Pi_{X_{k_S}}$ by a product of copies of $M_X$
which is indexed by the set of $k_S$-rational points $\{\{x_{i,j}\}_{j=1}^{n_i}\}_{i=1}^{n}$. More precisely, we have a natural exact sequence
$$1\to \prod _{p_i\in S}  (\prod _{x_{i,j}\in S_i} M_X)\to \Pi_{\Tilde U}^{\c-\cn}\to \Pi_{X_{k_S}}\to 1.$$

For $i\in \{1,...,n\}$, $j\in \{1,...,n_i\}$, 
the natural projection $\Pi_{\Tilde U_{x_{i,j}}}^{\c-\cn}\twoheadrightarrow G_{k_S}$ admits
group-theoretic sections $\tilde s_{i,j}:G_{k_S}\to \Pi_{\Tilde U_{x_{i,j}}}^{\c-\cn}$,
which lift the section  $s_{k_S}:G_{k_S}\to \Pi_{X_{k_S}}$, by the above
discussion. Here we use the fact that $s$ is uniformly good. Thus, $s_{k_S}$ is a good group-theoretic section.
The set of all such possible liftings $\tilde s_{i,j}$ of the section $s_{k_S}$ is a torsor under $H^1(G_{k_S},M_X)$. Any set of sections
$\{\{\tilde s_{i,j}:G_{k_S}\to \Pi_{\Tilde U_{x_{i,j}}}^{\c-\cn}\}_{j=1}^{n_i}\}_{i=1}^n$ as above  gives
rise to a unique section $\tilde s:G_{k_S}\to \Pi_{\Tilde U}^{\c-\cn}$, of the natural projection $\Pi_{\Tilde U}^{\c-\cn}\twoheadrightarrow G_k$,
which lifts the section $s_{k_S}$.

The Galois group $\Gal (k_S/k)$ acts naturally by outer automorphisms on both $G_{k_S}$ and
$\Pi_{\Tilde U}^{\c-\cn}$. This later action permutes the "components"
$\Pi_{\Tilde U_{x_{i,j}}}^{\c-\cn}$ of $\Pi_{\Tilde U}^{\c-\cn}$,
in a way that is compatible with the natural action of  $\Gal (k_S/k)$ on $S(k_S)$. More precisely, this action naturally
corresponds to the natural action of $\Gal (k_S/k)$ on the Chern classes $\{\{c(x_{i,j})\}_{j=1}^{n_1}\}_{i=1}^n$, where
$c(x_{i,j})\in H^2(X_{k_S},M_X)$ is the Chern class of the degree $1$ line bundle $\Cal O(x_{i,j})$.
Moreover, for $i\in \{1,...,n\}$, the outer action of $\Gal (k_S/k)$ on $\Pi_{\Tilde U}^{\c-\cn}$ stabilises 
$\prod _{x_{i,j}\in S_i}\Pi_{\Tilde U_{x_{i,j}}}^{\c-\cn}$, where the fibre product is over $X_{k_S}$.

The maximal (geometrically) cuspidally central quotient $\Pi_{U}^{\c-\cn}$ of  $\Pi_{U}$, 
with respect to the natural surjective homomorphism $\Pi_{U}\twoheadrightarrow  \Pi_{X}$, can be reconstructed
from $\Pi_{\Tilde U}^{\c-\cn}$ endowed with the outer action of  $\Gal (k_S/k)$ as follows. 

The profinite group 
$\Pi_{\Tilde U}^{\c-\cn}$ is slim (cf. [Mochizuki], Proposition 1.8, (i)). 
In loc. cit. $k$ is a $p$-adic local field, or a finite field, but the same arguments are valid if 
$G_k$ is slim, and under our assumption that $k$ is $\Sigma$-regular (cf. Definition 2.3.1). 

Thus, we have a natural exact sequence
$$1\to \Pi_{\Tilde U}^{\c-\cn} \to \Aut (\Pi_{\Tilde U}^{\c-\cn})\to \Out (\Pi_{\Tilde U}
^{\c-\cn})\to 1.$$

By pulling back this exact sequence by the natural homomorphism
$$\Gal (k_S/k)\to \Out (\Pi_{\Tilde U}^{\c-\cn}),$$ 
we obtain the exact sequence
$$1\to \Pi_{\Tilde U}^{\c-\cn} @>{\iota}>> \Pi_{U}^{\c-\cn}\to \Gal (k_S/k)\to 1.$$

Note that in order that a section $\tilde s:G_{k_S}\to \Pi_{\Tilde U}^{\c-\cn}$ as above, which lifts the section $s_{k_S}$,
descends to a section  $s_U^{\c-\cn}:G_{k}\to \Pi_{U}^{\c-\cn}$, which lifts the section $s$, it is necessary that
the section $\tilde s$ is invariant under the natural outer action of $\Gal (k_S/k)$ on $\Pi_{\Tilde U}^{\c-\cn}$.

We have the following commutative diagram of exact sequences:

$$
\CD
@.     @.   1   @.    1  \\
@.     @.   @VVV         @VVV        @.\\
1  @>>> \prod _{p_i\in S}  (\prod _{x_{i,j}\in S_i} M_X)  @>>>   \Pi_{\Tilde U}^{\c-\cn} @>>>  \Pi_{X_{k_S}} @>>>  1 \\
@.       @V{\id}VV       @V{\iota}VV     @VVV     \\
1  @>>> \prod _{p_i\in S}  (\prod _{x_{i,j}\in S_i} M_X)  @>>>   \Pi_{U}^{\c-\cn} @>>>  \Pi_{X} @>>>  1 \\
@. @. @VVV    @VVV  @. \\
@. @.  \Gal (k_S/k)  @>{\id}>> \Gal (k_S/k) @. \\
@. @. @VVV    @VVV \\
@. @. 1 @.  1\\
\endCD
$$
which by pull back via the sections $s:G_k\to \Pi_X$, and $s_{k_S}:G_{k_S}\to \Pi_{X_{k_S}}$, gives rise to the following commutative diagram of exact sequences:
$$
\CD
@.     @.   1   @.    1  \\
@.   @.         @VVV        @VVV \\
1  @>>> \prod _{p_i\in S}  (\prod _{x_{i,j}\in S_i} M_X)  @>>>   \Tilde D_{S} @>>>  G_{k_S}=s_{k_S}(G_{k_S}) @>>>  1 \\
@.       @V{\id}VV       @V{\iota}VV     @VVV     \\
1  @>>> \prod _{p_i\in S}  (\prod _{x_{i,j}\in S_i} M_X)  @>>>   D_S @>>>  G_k=s(G_k) @>>>  1 \\
@. @.   @VVV    @VVV  @. \\
@.  @. \Gal (k_S/k)  @>{\id}>> \Gal (k_S/k) @. \\
@. @. @VVV    @VVV \\
@. @. 1 @.  1\\
\endCD
$$
where $D_S$ (resp. $\Tilde D_S$) is the fibre of $s(G_k)$ (resp. fibre of $s_{k_S}(G_{k_S})$) in  
$\Pi_{U}^{\c-\cn}$ (resp. in $\Pi_{\Tilde U}^{\c-\cn}$), and the upper right square is cartesian.

In order to construct a section $s_U^{\c-\cn}:G_{k}\to \Pi_{U}^{\c-\cn}$, which lifts the section $s$,
it is equivalent to construct a section $G_k\to D_S$ of the natural projection $D_S\twoheadrightarrow G_k$.

The group extension $1\to \prod _{p_i\in S}  (\prod _{x_{i,j}\in S_i} M_X)  \to \Tilde D_{S} \to G_{k_S} \to 1$ is split by the above discussion.
We will show that the group extension $D_{S}$ splits.
To simplify notations write $M\defeq \prod _{p_i\in S}  (\prod _{x_{i,j}\in S_i} M_X)$. We will carefully distinguish
between the $G_k$, and $G_{k_S}$, module structure of $M$.

Let $[D_S]\in H^2(G_k,M)$  (resp. $[\Tilde D_S]\in H^2(G_{k_S},M)$) be the class of the group extension
$D_S$ (resp. $\Tilde D_S$).
We have a natural exact sequence arising from the Hochschild-Serre spectral sequence (associated to  the natural inclusion $G_{k_S}\subseteq G_k$)
$$0\to H^1(\Gal (k_S/k), H^1(G_{k_S}, M))\to H^2(G_k,M) \to
H^2(G_{k_S},M)^{\Gal (k_S/k)}.$$ 
Here we use the fact that $H^0(G_{k_S},M)=0$, as follows from the fact that $k$ is $\Sigma$-regular.

The image of $[D_S]$ in $H^2(G_{k_S},M)^{\Gal (k_S/k)}$, which is $[\Tilde D_S]$, is trivial by assumption.
We will show that the homomorphism  $H^2(G_k,M) \to
H^2(G_{k_S},M)^{\Gal (k_S/k)}$ is injective, or equivalently that $H^1(\Gal (k_S/k), H^1(G_{k_S}, M))=0$.

The $G_k$-module $M$ is naturally identified with $\prod _{p_i\in S}(\Ind _{G_{k(p_i)}}^{G_k}M_X)$, where $k(p_i)$ denotes the residue field at $p_i$.
The cohomology group $H^2(G_k,M)$ is thus naturally identified with $\prod _{p_i\in S} H^2(G_{k(p_i)},M_X)$, which is naturally identified
with $\prod _{p_i\in S} T_{\Sigma} \Br (k(p_i))$, where $T_{\Sigma} \Br (k(p_i))$ denotes the $\Sigma$-adic Tate module of the Brauer group
$\Br (k(p_i))$. In particular, $H^2(G_k,M)$ is torsion free. Hence the kernel of the above homomorphism $H^2(G_k,M) \to
H^2(G_{k_S},M)^{\Gal (k_S/k)}$, which is torsion by the above exact sequence, is trivial. Thus, $[D_S]=0$, and the group extension
$D_S$ splits. 

Moreover, the set of all possible splittings of the group extension $D_S$ is a torsor under 
$H^1(G_{k_S},M)^{\Gal (k_S/k)}$, as follows from the natural identification
$ H^1(G_k,M)\isom H^1(G_{k_S},M)^{\Gal (k_S/k)}$ arising from the five terms exact sequence in low degree in the Hochschild-Serre spectral sequence
(the cohomology groups $H^i(\Gal (k_S/k), H^0(G_{k_S},M_X))$ are trivial in all degrees by our assumption on $k$).

This finishes the proof of assertion (i).

Next, we prove assertion (ii). 

We will argue by induction on the cardinality of the finite set $S\defeq X\setminus U$.
The case where $S$ consists of a single closed point is treated in the above proof of assertion (i).

Suppose that $S=\{p_1,p_2,...,p_n\}$ consists of $n$ closed points (not necessarily rational over $k$), and assume that assertion (ii) holds true 
for all open subschemes $U=X\setminus S$, in the case where $S$ consists of at most $n-1$ closed points.

Let $k_S/k$ be the minimal Galois extension over which all points in $S$ are rational.
By our induction hypothesis, there exist sections $s_i^{\c-\cn}:G_k\to \Pi_{U_{p_i}}^{\c-\cn}$
of the natural projection $\Pi_{U_{p_i}}^{\c-\cn}\twoheadrightarrow G_k$,
which lift the section $s$, for $i\in \{1,...,n\}$. Here $U_{p_i}\defeq X\setminus \{p_i\}$.
Over $k_S$, the section $s_i^{\c-\cn}$ restricts to a section $\Tilde s_i^{\c-\cn}:G_{k_S}\to \Pi_{\Tilde U_{p_i}}^{\c-\cn}$
of the natural projection $\Pi_{\Tilde U_{p_i}}^{\c-\cn}\twoheadrightarrow G_{k_S}$,
where $\Tilde U_{p_i}\defeq U_{p_i}\times _k k_S$. 

We use the same notations as in the proof of assertion (i), and
we fix a base point of the torsor of splittings of the group extensions $D_S$, and $\Tilde D_S$ (cf. above commutative diagram),
which are split extensions by the proof of assertion (i).

For $i\in \{1,...,n\}$, let $S_i\defeq \{x_{i,j}\}_{i=1}^{n_i}$ be the set of points of $X_{k_S}$ above $p_i$. The section $s_i^{\c-\cn}$ is given by an element
$\sigma _i\in H^1(G_{k_S},\prod _{x_{i,j}\in S_i}  M_X)^{\Gal (k_S/k)}$, as follows from the proof of assertion (i).

The element $\{\sigma_i\}_{i=1}^n \in \prod _{p_i\in S} H^1(G_{k_S},\prod _{x_{i,j}\in S_i}  M_X)^{\Gal (k_S/k)}$ determines 
a section $\tilde s_U^{\c-\cn}:G_{k_S}\to \Pi_{U_{k_S}}^{\c-\cn}$  of the natural projection
$\Pi_{U_{k_S}}^{\c-\cn}\twoheadrightarrow G_k$, where $U_{k_S}\defeq U\times _kk_S$. This section $\tilde s_U^{\c-\cn}$
descends uniquely 
to a section $s_U^{\c-\cn}:G_k\to \Pi_U^{\c-\cn}$ of the natural projection $\Pi_U^{\c-\cn}\twoheadrightarrow G_k$ (cf. Proof of assertion (i)). 
The section  $s_U^{\c-\cn}$ constructed in this way has the property required in assertion (ii) .

This finishes the proof of assertion (ii).

Finally, assertion (iii) follows formally from assertion (ii).
\qed
\enddemo

The following description of maximal cuspidally abelian arithmetic fundamental
groups plays an important role in the proof of our main result in Theorem 2.6.

\proclaim {Proposition 2.5 (Cuspidally Abelian Arithmetic Fundamental Groups over Slim Fields)}
Assume that the field $k$ is slim, and $\Sigma$-regular (cf. Definition 2.3.1). 
Let $U\subseteq X$
be a non-empty open subscheme of $X$, and $\Pi_{U}^{\c-\ab}$ the maximal cuspidally abelian quotient of
$\Pi_{U}$, with respect to the natural surjective homomorphism $\Pi_{U}\twoheadrightarrow \Pi_X$.

For a finite \'etale Galois cover $X'\to X$, with Galois group $\Gal (X'/X)$, let $U'\defeq U\times _XX'$,
and $\Pi_{U'}^{\c-\cn}$ the maximal (geometrically) cuspidally central quotient of $\Pi_{U'}$,
with respect to the natural surjective homomorphism $\Pi_{U'}\twoheadrightarrow \Pi_{X'}$. Then
$\Pi_{U'}^{\c-\cn}$ is slim. 

The natural outer action of $Gal (X'/X)$ on $\Pi_{X'}$ extends to an outer action on  $\Pi_{U'}^{\c-\cn}$.
Denote by $\Pi_{U'}^{\c-\cn} \rtimes ^{\out}\Gal (X'/X)$ the profinite group
which is obtained by pulling back the exact sequence 
$$1\to \Pi_{U'}^{\c-\cn}\to \Aut(\Pi_{U'}^{\c-\cn})\to \Out (\Pi_{U'}^{\c-\cn})\to 1,$$ 
by the natural homomorphism
$$\Gal (X'/X)\to \Out (\Pi_{U'}^{\c-\cn}).$$ 

Thus, we have a natural exact sequence:
$$1\to \Pi_{U'}^{\c-\cn} \to \Pi_{U'}^{\c-\cn} \rtimes ^{\out}\Gal (X'/X)\to \Gal (X'/X)\to 1,$$
which inserts into the following commutative diagram of exact sequences:
$$
\CD
@.   1   @.   1 \\
@.   @VVV  @VVV \\
@.    I_{U'}^{\c-\cn}   @>{\id}>>  I_{U'}^{\c-\cn}   \\
@.  @VVV    @VVV \\
1@>>> \Pi_{U'}^{\c-\cn} @>>> \Pi_{U'}^{\c-\cn} \rtimes ^{\out}\Gal (X'/X) @>>>  \Gal (X'/X) @>>> 1\\
@.   @VVV            @VVV                @V{\id}VV  \\
1@>>> \Pi_{X'}@>>> \Pi_{X} @>>> \Gal (X'/X)@>>> 1\\
@.     @VVV  @VVV   @VVV \\
@. 1 @.   1 @. 1\\
\endCD
$$

Then we have a natural isomorphism:
$$\Pi_{U}^{\c-\ab}\isom \underset{X'\to X}\to {\varprojlim} (\Pi_{U'}^{\c-\cn} \rtimes ^{\out}\Gal (X'/X)),$$
where the projective limit is taken over all  finite \'etale Galois cover $X'\to X$.
\endproclaim

\demo{Proof} This is proven in [Mochizuki], Proposition 1.14, (i), in the case where the field $k$ is a
$p$-adic local field, or a finite field. 

The same arguments as in loc. cit. are valid in the general case where $G_k$ is assumed to be slim, and 
$k$ is $\Sigma$-regular.

More precisely, the fact that $\Pi_{U'}^{\c-\cn}$ is slim follows from
the slimness of $G_k$, and the fact that $k$ is $\Sigma$-regular (compare with [Mochizuki], 1.8. (i)). 

The rest of the assertion follows from the observation
that the kernel $I_U^{\ab}$ of the surjective homomorphism  $\Pi_{U}^{\c-\ab}\twoheadrightarrow \Pi_X$ is naturally
identified with the projective limit of the kernels of the natural surjective homomorphisms
$\Pi_{U'}^{\c-\cn} \twoheadrightarrow \Pi_{X'}$, where the projective limit is taken over all finite \'etale
Galois cover $X'\to X$, and $U'\defeq U\times _XX'$ (cf. loc. cit.).
\qed
\enddemo

The following result is our main result in this section. It shows that (uniformly) good sections of arithmetic fundamental groups
over slim, and $\Sigma$-regular, base fields behave well with respect to the cuspidalisation problem.

\proclaim {Theorem 2.6 (Lifting of Uniformly Good Sections to Cuspidally abelian Arithmetic Fundamental Groups over Slim Fields)} 
Assume that the field $k$ is slim, and $\Sigma$-regular (cf. Definition 2.3.1).
Let $s:G_k\to \Pi_X$ be a uniformly good section (in the sense of Definition 1.4.1) of the natural projection
$\Pi_X\twoheadrightarrow G_k$. Then the followings hold.

(i) Let $U\subseteq X$ be a non-empty open subscheme of $X$.
For a finite \'etale Galois cover $X'\to X$, with Galois group $\Gal (X'/X)$, let $U'\defeq U\times _XX'$,
and $\Pi_{U'}^{\c-\cn}$ the maximal (geometrically) cuspidally central quotient of $\Pi_{U'}$,
with respect to the natural surjective homomorphism $\Pi_{U'}\twoheadrightarrow \Pi_{X'}$.

Then there exists a section $s_{X'}: G_k\to \Pi_{U'}^{\c-\cn} \rtimes ^{\out}\Gal (X'/X)$ of the natural projection
$\Pi_{U'}^{\c-\cn} \rtimes ^{\out}\Gal (X'/X)\twoheadrightarrow G_k$, which lifts the section $s$, i.e. which inserts into the following commutative diagram:

$$
\CD
G_k  @>{s_{X'}}>> \Pi_{U'}^{\c-\cn} \rtimes ^{\out}\Gal (X'/X) \\
@V{\id}VV    @VVV   \\
G_k @>s>>  \Pi _X\\
\endCD
$$

(ii) With the same notations as in (i). There exists, for each finite \'etale Galois cover $X'\to X$, 
a section $s_{X'}: G_k\to \Pi_{U'}^{\c-\cn} \rtimes ^{\out}\Gal (X'/X)$ of the natural projection
$\Pi_{U'}^{\c-\cn} \rtimes ^{\out}\Gal (X'/X)\twoheadrightarrow G_k$ as in (i) (i.e. which lifts the section $s$), and such that for each factorisation 
$X'\to X''\to X$, where $X''\to X$ is Galois, we have a commutative diagram:

$$
\CD
G_k  @>{s_{X'}}>> \Pi_{U'}^{\c-\cn} \rtimes ^{\out}\Gal (X'/X) \\
@V{\id}VV    @VVV   \\
G_k  @>{s_{X''}}>> \Pi_{U''}^{\c-\cn} \rtimes ^{\out}\Gal (X''/X) \\
@V{\id}VV    @VVV   \\
G_k @>s>>  \Pi _X\\
\endCD
$$

(iii)\ Let $U\subseteq X$ be a non-empty open subscheme of $X$, and
$\Pi_{U}^{\c-\ab}$ the maximal cuspidally abelian quotient of $\Pi_{U}$, with respect to the natural surjective
homomorphism $\Pi_{U}\twoheadrightarrow \Pi_X$ (cf. 2.1.1). Then there exists a section $s_U^{\c-\ab}:G_k\to \Pi_{U}^{\c-\ab}$
of the natural projection  $\Pi_{U}^{\c-\ab}\twoheadrightarrow G_k$, which lifts the section $s$, i.e. which inserts into
the following commutative diagram:
$$
\CD
G_k @>s_U^{\c-\ab}>>  \Pi_{U}^{\c-\ab} \\
@V{\id}VV     @VVV  \\
G_k   @>{s}>> \Pi_X
\endCD
$$

Moreover, the set of all possible liftings $s_U^{\c-\ab}$ of $s$ is a torsor under the group $H^1(G_k,I_U^{\ab})$.
Here the $G_k$-module structure of $I_U^{\ab}$ is naturally induced by the section $s$.

(iv)\ There exists, for every non-empty open subscheme $U\defeq X\setminus S$ of $X$, a section
$s_U^{\c-\ab}:G_k\to \Pi_{U}^{\c-\ab}$ as in (iii) (i.e. which lifts the section $s$), such that for every non-empty open subscheme
$V\defeq X\setminus T$ of $X$, with $U\subseteq V$, we have the following commutative diagram:
$$
\CD
G_k @>s_U^{\c-\cn}>>  \Pi_{U}^{\c-\ab} \\
@V{\id}VV     @VVV  \\
G_k   @>s_V^{\c-\cn}>> \Pi_{V}^{\c-\ab}\\
\endCD
$$
where the right vertical homomorphism is the natural one.
\endproclaim

\demo{Proof} First, assertion (iii) follows formally from assertion (ii), and Proposition 2.5.

Next, we prove assertion (i).

Let $X'\to X$ be an \'etale Galois cover with Galois group $\Gal (X'/X)$, and $U'\defeq U\times _XX'$.
Let $S\defeq \{x_i\}_{i=1}^n\defeq X\setminus U$. Assume that, 
for each $i\in \{1,...,n\}$, the set of points $S_i\defeq \{x'_{i,j}\}_{j=1}^{n_i}$
in $X'$ above $x_i$ are $k$-rational.
In particular, the image of the open subgroup $\Pi_{X'}$ of $\Pi_X$ in $G_k$,
via the natural projection $\Pi_X\twoheadrightarrow G_k$, coincides with $G_k$, and 
the group-theoretic section $s:G_k\to \Pi_X$ restricts to a group-theoretic section $s':G_k\to \Pi_{X'}$ of the natural projection $\Pi_{X'}\twoheadrightarrow G_k$.

Let $\Pi_{U'}^{\c-\cn}$ be the maximal (geometrically) cuspidally central quotient of $\Pi_{U'}$, with respect to the natural
surjective homomorphism $\Pi_{U'}\twoheadrightarrow \Pi_{X'}$. The section $s'$ is uniformly good, since $s$ is uniformly good. 
This follows from the very definition of uniform goodness (cf. Definition 1.4,1). Thus, $s'$ lifts to a section
${s'}_{U'}^{\c-\cn}:G_k\to \Pi_{U'}^{\c-\cn}$ of the natural projection  $\Pi_{U'}^{\c-\cn}\twoheadrightarrow G_k$ by Proposition
2.4, (i), which is also a section $s_{X'}:G_k\to \Pi_{U'}^{\c-\cn} \rtimes ^{\out}\Gal (X'/X)$ of the natural projection
$\Pi_{U'}^{\c-\cn} \rtimes ^{\out}\Gal (X'/X)\twoheadrightarrow G_k$, and which lifts the section $s$ as required in (ii) (cf.
commutative diagram in the statement of Proposition 2.5). 

In the case where the points $S_i\defeq \{x'_{i,j}\}_{j=1}^{n_i}$
in $X'$ above the points $\{x_i\}_{i=1}^n$ in $S$ are not necessarily $k$-rational, one argues using a descent argument, similar to the one used in the proof of assertion (i)
in Proposition 2.4, to show the existence of a section $s_{X'}:G_k\to \Pi_{U'}^{\c-\cn} \rtimes ^{\out}\Gal (X'/X)$ of the natural projection
$\Pi_{U'}^{\c-\cn} \rtimes ^{\out}\Gal (X'/X)\twoheadrightarrow G_k$, which lifts the section $s$. One first constructs such a section over a field
over which all points $x'_{i,j}$ are rational, and then descend this section to $k$.

This finishes the proof of assertion (i).

Next, we prove assertion (ii). 

Let $X'\to X$ be an \'etale Galois cover, with Galois group $\Gal (X'/X)$, and $U'\defeq U\times _XX'$.
Let $S\defeq \{x_i\}_{i=1}^n\defeq X\setminus U$.

We argue by induction on the degree $\vert X':X\vert$ of the Galois cover $X'\to X$. In the case where
$\vert X':X\vert=1$, i.e. $X'=X$, we fix a section  $s_X\defeq s_U^{\c-\cn}:G_k\to \Pi_{U}^{\c-\cn}$ of the
natural projection  $\Pi_{U}^{\c-\cn}\twoheadrightarrow G_k$, which lifts the section $s$, and which exists by Proposition 2.4, (i).

Let $X'\to X$ be a finite \'etale Galois cover of degree $n>1$. We assume that assertion (ii) holds true
for every finite \'etale Galois cover $X''\to X$ of degree less or equal to $n-1$.

First, we treat the case where the Galois group $\Gal (X'/X)$ is a simple group, i.e. there are no intermediate Galois covers $X'\to X''\to X$,
with $1< \vert X'':X\vert  < \vert X':X\vert$. We have to show in this case that the section $s_X\defeq s_U^{\c-\cn}:G_k\to \Pi_{U}^{\c-\cn}$
can be lifted to a section $s_{X'}:G_k\to \Pi_{U'}^{\c-\cn} \rtimes ^{\out}\Gal (X'/X)$ of the natural projection
$\Pi_{U'}^{\c-\cn} \rtimes ^{\out}\Gal (X'/X)\twoheadrightarrow G_k$, i.e. we have to construct a section $s_{X'}$ which fits into the following
commutative diagram

$$
\CD
G_k  @>{s_{X'}}>>  \Pi_{U'}^{\c-\cn} \rtimes ^{\out}\Gal (X'/X) \\
@V{\id}VV    @VVV  \\
G_k @>{s_X}>> \Pi_{U}^{\c-\cn} \\
\endCD
$$

Assume that, for each $i\in \{1,...,n\}$, the set of points $S_i\defeq \{x'_{i,j}\}_{j=1}^{n_i}$
in $X'$ above the point $x_i\in S$ are $k$-rational.
In particular, the image of the open subgroup $\Pi_{X'}$ of $\Pi_X$ in $G_k$,
via the natural projection $\Pi_X\twoheadrightarrow G_k$, coincides with $G_k$, and 
the group-theoretic section $s:G_k\to \Pi_X$ restricts to a section $s':G_k\to \Pi_{X'}$.

We have a natural commutative diagram:
$$
\CD
1  @>>> \prod _{x_i\in S}  (\prod _{x'_{i,j}\in S_i} M_X)  @>>>   \Pi_{U'}^{\c-\cn} \rtimes ^{\out}\Gal (X'/X) @>>>  \Pi_{X} @>>>  1 \\
@.       @VVV       @VVV     @V{\id}VV     \\
1  @>>> \prod _{x_i\in S}  M_X  @>>>   \Pi_{U}^{\c-\cn} @>>>  \Pi_{X} @>>>  1 \\
@.    @VVV       @VVV \\
@. 1 @.  1 \\
\endCD
$$
where the horizontal lines are exact, and the two left vertical maps are surjective. 

The above diagram induces by pull back via the sections $s:G_k\to \Pi_X$, 
and $s':G_{k}\to \Pi_{X'}$, the following commutative diagram of exact sequences:
$$
\CD
@.   1   @. 1  \\
@.  @VVV    @VVV  \\
@.   M @>{\id}>>   M     \\
@. @VVV     @VVV   \\
1  @>>> \prod _{x_i\in S}  (\prod _{x'_{i,j}\in S_i} M_X)  @>>>   D_{S'} @>>>  G_k=s'(G_k) @>>>  1 \\
@.       @VVV       @VVV     @V{\id}VV     \\
1  @>>> \prod _{x_i\in S}  M_X  @>>>   D_S @>>>  G_k=s(G_k) @>>>  1 \\
@. @VVV    @VVV \\
@. 1 @.  1\\
\endCD
$$
where $D_S$ (resp. $D_{S'}$) is the fibre of $s(G_k)$ in  
$\Pi_{U}^{\c-\cn}$ (resp. fibre of $s'(G_k)$ in  $\Pi_{U'}^{\c-\cn}$), and $M$ is the kernel of the natural homomorphism $D_S\to D_{S'}$.

The above commutative diagram is a push out diagram via the natural morphism (of $G_k$-modules)
$\prod _{x_i\in S}  (\prod _{x'_{i,j}\in S_i} M_X)\to \prod _{x_i\in S}  M_X$, on the left of the above diagram, which maps the copy of $M_X$
indexed by the point $x'_{i,j}$ identically to the copy of $M_X$ which is indexed by $x_i$.

The group extensions $D_S$, and $D_{S'}$, are split, since the sections $s$ and $s'$ are (uniformly) good sections  
(cf. proof of Proposition 2.4, (i)).
We fix base points of the torsors of splittings of the group extensions $D_S$, and $D_{S'}$, which are compatible with
the morphism $D_{S'}\to D_{S}$.

The section $s_X\defeq s_U^{\c-\cn}:G_k\to \Pi_{U}^{\c-\cn}$ of the
natural projection  $\Pi_{U}^{\c-\cn}\twoheadrightarrow G_k$, which lifts the section $s$ (cf. above induction hypothesis),
is uniquely determined by a section $\tilde s:G_k\to D_S$
of the natural projection $D_S\twoheadrightarrow G_k$. Such a section $\tilde s$ is uniquely determined by an element
of $H^1(G_k,\prod _{x_i\in S}  M_X)$. Similarly, a section of the natural projection $D_{S'}\to G_k$ is uniquely determined by 
an element of $H^1(G_k,\prod _{x_i\in S}  (\prod _{x'_{i,j}\in S_i} M_X))$.

Moreover, in order to construct a section $s_{X'}\defeq s_{U'}^{\c-\cn}:G_{k}\to \Pi_{U'}^{\c-\cn}\rtimes ^{\out}\Gal (X'/X)$, which lifts the section $s_X\defeq s_U^{\c-\cn}:G_k
\to \Pi_U^{\c-\cn}$,
it is equivalent to construct a section $\tilde s':G_k\to D_{S'}$ of the natural projection $D_{S'}\twoheadrightarrow G_k$, which lifts the section $\tilde s$, i.e.
such that we have a commutative diagram:
$$
\CD
G_k @>{\tilde s'}>>  D_{S'} \\
@V{\id}VV     @VVV  \\
G_k   @>{\tilde s}>> D_S\\
\endCD
$$

Such a section $\tilde s'$ exists, since the natural map $H^1(G_k,\prod _{x_i\in S}  (\prod _{x'_{i,j}\in S_i} M_X))\to H^1(G_k,\prod _{x_i\in S}  M_X)$
is surjective. The set of all such possible liftings $\tilde s'$ of the section $\tilde s$ is a torsor under the group $H^1(G_k,M)$.

In the general case where, 
for $i\in \{1,...,n\}$, the points 
in $X'$ above the point $x_i\in S$ are not necessarily $k$-rational, one argues in a similar way as in the proof of Proposition 2.4, (i), i.e.
by passing to a finite extension where all these points are rational, and then descend to $k$. 

More precisely, let $k_S/k$ be the minimal Galois extension, with Galois group $\Gal (k_S/k)$, over which all points in $X'$ above the points in $S$
are rational. Let $\Tilde S\defeq \{\tilde x_{i'}\}_{i'=1}^{n'}\subset X\times _kk_S$ be the set of points above the points in $S$.

For $i'\in \{1,...,n'\}$, let $\Tilde S_{i'}\defeq \{\tilde x'_{i',j}\}_{j=1}^{n_{i'}}$ be the set of points 
of $X'_{k_S}\defeq X'\times _kk_S$ above $\tilde x_{i'}$,
which are $k_S$-rational points.

With the same notations as above, write $M_1\defeq \prod _{\tilde x_{i'}\in \Tilde S}  (\prod _{\tilde x'_{i',j}\in \Tilde S_{i'}} M_X))$, $M_2\defeq \prod _{\tilde x_{i'}\in \Tilde S}  M_X$,
and $M\defeq \Ker(M_1\to M_2)$.

In order to perform descent, one only needs the fact that the natural map 
$H^1(G_{k_S},M_1)^{\Gal (k_S/k)}\to H^1(G_{k_S},M_2)^{\Gal (k_S/k)}$
is surjective, as follows from the above proof in the case where $k_S=k$ (compare with the descent argument in the proof of Proposition 2.4, (i)). 

This later map is surjective. Indeed, we have a long exact cohomology sequence
$$0\to H^1(G_{k_S},M)^{\Gal (k_S/k)}\to 
H^1(G_{k_S},M_1)^{\Gal (k_S/k)}\to H^1(G_{k_S},M_2)^{\Gal (k_S/k)}$$
$$\to H^1(\Gal (k_S/k), H^1(G_{k_S},M))\to ...,$$
which arises from the exact sequence 
$$0\to H^1(G_{k_S},M) \to H^1(G_{k_S},M_1)\to H^1(G_{k_S},M_2)\to 0$$
of $\Gal (k_S/k)$-modules. 
Moreover, $H^1(\Gal (k_S/k), H^1(G_{k_S},M))$ vanishes as follows from the following 
commutative diagram of exact sequences:

$$
\CD
1   @>>>  H^1(\Gal (k_S/k),H^1(G_{k_S},M))   @>>>   H^2(G_k,M)  @>>>   H^2(G_{k_S},M)^{\Gal (k_S/k)}\\
@.         @VVV       @VVV    @VVV  \\
1   @>>>  H^1(\Gal (k_S/k),H^1(G_{k_S}, M_1))   @>>>   H^2(G_k, M_1)  
@>>>   H^2(G_{k_S}, M_1)^{\Gal (k_S/k)}\\
\endCD
$$
and the facts that the map $H^2(G_k,M)\to H^2(G_k,M_1)$ is injective, and 
\newline
$H^1(\Gal (k_S/k),H^1(G_{k_S}, M_1))$
vanishes (cf. proof of Proposition 2.4, (i)).

Or alternatively, we have natural identifications $H^1(G_{k_S}, M_1))^{\Gal (k_S/k)}\isom H^1(G_{k}, M_1))$, 
$H^1(G_{k_S}, M_2))^{\Gal (k_S/k)}\isom H^1(G_{k}, M_2))$ (cf. proof of Proposition 2.4, (ii)), and the natural
map $H^1(G_k,M_1)\to H^1(G_k,M_2)$ is surjective.

In the general case, where $\Gal (X'/X)$ is not simple, let $X'\to X_i\to X$, $i\in \{1,...,n\}$, be the intermediate Galois covers of $X'\to X$, with 
$1<\vert X_i:X\vert <\vert X':X\vert$.
We can assume, without loss of generality, that $X'=\prod _{i=1}^nX_i$ where the fibre product is over $X$.
We can also assume, without loss of generality, that $U\defeq X\setminus \{x\}$ is the complement of a single point $x\in X$. 

For $i\in \{1,...,n\}$, assume that the set of points $\{x_{j}\}_{j\in J_i}$ of $X_i$ above $x$ are $k$-rational. The 
general case is treated using a descent argument similar to the one used above.

For $i\in \{1,...,n\}$, let $U_i\defeq U\times _XX_i$. 
By induction hypothesis, there exist group-theoretic sections $s_{X_i}\defeq s_{U_i}^{\c-\cn}:G_{k}\to \Pi_{U_i}^{\c-\cn}\rtimes ^{\out}\Gal (X_i/X)$ of the natural projection
$\Pi_{U_i}^{\c-\cn}\rtimes ^{\out}\Gal (X_i/X)\twoheadrightarrow G_k$, which lift the section $s_X\defeq s_U^{\c-\cn}:G_k\to \Pi_U^{\c-\cn}$.
The section $s_X$ is uniquely determined by an element $\alpha\in H^1(G_k,M_X)$ (cf. proof of assertion (i) in Proposition 2.4).
The section $s_{X_i}$ is uniquely determined by an element $\alpha_i\in \prod _{j\in J_i} H^1(G_k,M_X)$, which lifts the element $\alpha$ via the natural homomorphism
$\prod _{j\in J_i} H^1(G_k,M_X) \to H^1(G_k,M_X)$ (cf. proof of assertion (ii)). 

One has to construct a section $s_{X'}\defeq s_{U'}^{\c-\cn}:G_{k}\to 
\Pi_{U'}^{\c-\cn}\rtimes ^{\out}\Gal (X'/X)$ of the natural projection $\Pi_{U'}^{\c-\cn}\rtimes ^{\out}\Gal (X'/X)\twoheadrightarrow G_k$, which lifts the section $s_{X_i}$, for
$i\in \{1,...,n\}$. It follows from the proof of assertion (ii), that such a section is uniquely determined by an element
$\alpha'\in \prod _{i=1}^{n} (\prod _{j\in J_i} H^1(G_k,M_X))$, which lifts the element $\alpha_i$ via the natural homomorphism
$\prod _{i=1}^{n} (\prod _{j\in J_i} H^1(G_k,M_X))\to \prod _{j\in J_i} H^1(G_k,M_X)$, for $i\in \{1,...,n\}$. Such an element $\alpha'$ exists,
since the product $\prod _{i=1}^{n} (\prod _{j\in J_i} H^1(G_k,M_X))$ is a fibre product of the various 
$\prod _{j\in J_i} H^1(G_k,M_X)$ above $H^1(G_k,M_X)$, via the various
homomorphisms 
$\prod _{j\in J_i} H^1(G_k,M_X)\to H^1(G_k,M_X)$, by our assumption that $X'$ is the fibre product over $X$ of the various $X_i$.

This finishes the proof of assertion (ii).

Finally, the proof of assertion (iv) is similar to the proof of assertion (ii) in Proposition 2.4, and is done by induction on the cardinality
of the set of points $S\defeq X\setminus U$. The case where $S$ consists of a single point follows from assertion (ii). 
\qed
\enddemo

As a consequence of Theorem 2.6 one obtains the following.

\proclaim {Theorem 2.7 (Lifting of Uniformly Good Sections to Cuspidally abelian Galois Groups over Slim Fields)} Assume that
the field $k$ is slim, and $\Sigma$-regular. Let $s:G_k\to \Pi_X$ be a group-theoretic section of the natural projection
$\Pi_X\twoheadrightarrow G_k$. Then the section $s:G_k\to \Pi_X$ is uniformly good (in the sense of Definition 1.4.1),
if and only if there exists a section $s^{\c-\ab}:G_k\to G_{X}^{\c-\ab}$ of the
natural projection  $G_{X}^{\c-\ab}  \twoheadrightarrow G_k$, which lifts the section $s$,
i.e. which inserts into the following commutative diagram:
$$
\CD
G_k @>s^{\c-\ab}>>   G_{X}^{\c-\ab}\\
@V{\id}VV     @VVV  \\
G_k   @>{s}>> \Pi_X
\endCD
$$

Moreover, the set of all possible liftings $s^{\c-\ab}$ of $s$ is a torsor under the group $H^1(G_k,\Tilde I_X^{\ab})$.
Here the $G_k$-module structure of $\Tilde I_X^{\ab}$ is naturally induced by the section $s$.
\endproclaim

\demo{Proof} The if part follows from the various definitions, and from the geometric interpretation of the various
$\Pi_{U'}^{\c-cn}$ involved in the proof of Theorem 2.6 in terms of Chern classes of line bundles.
The only if part follows formally from Theorem 2.6, (iv), using the natural identification
$G_{X} ^{\c-\ab}\isom \underset{U}\to {\varprojlim}\ \Pi_U^{\c-ab}$,
where the projective limit is over all non-empty open subschemes of $X$.
\qed
\enddemo

\subhead
\S 3. Applications to the Grothendieck Anabelian Section Conjecture
\endsubhead
In this section we apply the theory of cuspidalisation of (uniformly) good sections of arithmetic fundamental
groups, which was investigated in $\S2$,
to the Grothendieck anabelian section conjecture.

We follow the notations in $\S1$. 

In particular, $X$ is a proper, smooth, geometrically connected, and hyperbolic algebraic curve over the field $k$,
$\Sigma \subseteq \Primes$ is
a non-empty set of prime integers, with $\char (k)\notin \Sigma$, and we have the natural exact sequence
$$1\to \Delta_X\to \Pi_X @>{\pr_{X,\Sigma}}>> G_k\to 1,\tag {$3.1$}$$
where $\Pi_X$ is the geometrically pro-$\Sigma$ arithmetic fundamental group of $X$ (cf. exact sequence (1.2)).

\subhead
{3.1}
\endsubhead
Sections of the natural projection $\pr _{X,\Sigma}:\Pi_X\twoheadrightarrow G_k$ arise naturally from $k$-rational points of $X$.

More precisely, let $x\in X(k)$ be a rational point. Then $x$ determines a decomposition subgroup
$D_x\subset \Pi_X$, which is only defined up to conjugation by the elements of $\Delta _X$, and which maps isomorphically
to $G_k$ via the projection $\Pi_X\twoheadrightarrow G_k$. Hence, the decomposition group $D_x$ determines a group-theoretic
section 
$$s_x:G_k\to \Pi_X$$ 
of the natural projection $\Pi_X \twoheadrightarrow G_k$, which is only defined up to conjugation by the elements of $\Delta _X$.

Let $\overline {\Sec} _{\Pi_X}$ be the set of conjugacy classes of all continuous group-theoretic sections $G_k\to \Pi_X$
of the natural projection $\Pi_X \twoheadrightarrow G_k$, modulo inner conjugation by the elements of $\Delta _X$.
We have a natural set-theoretic map
$$\varphi_X\defeq \varphi_{X,\Sigma}: X(k)\to  \overline {\Sec} _{\Pi_X},$$
$$x\mapsto \varphi_X (x)\defeq [s_x],$$
where $[s_x]$ denotes the image, i.e. conjugacy class, of the section $s_x$ in $\overline {\Sec} _{\Pi_X}$.

\definition {Definition 3.1.1}\ Let $s:G_k\to \Pi_X$ be a group-theoretic section of the natural projection
$\Pi_X \twoheadrightarrow G_k$. We say that the section $s$ is point-theoretic, if the conjugacy class $[s]$ of $s$ in
$\overline {\Sec} _{\Pi_X}$ belongs to the image of the map $\varphi_X$.
\enddefinition

The following conjecture, which was formulated by Grothendieck, 
is the main anabelian conjecture concerning sections of arithmetic fundamental groups (cf. [Grothendieck]).
 
\definition {Grothendieck's Anabelian Section Conjecture (GASC)}  Assume that $k$ is finitely generated over the prime field
$\Bbb Q$, and that $\Sigma =\Primes$. Then the map $\varphi_X: X(k)\to  \overline {\Sec} _{\Pi_X}$ is bijective.
\enddefinition

The injectivity of the map $\varphi_X$ under the assumptions in the GASC is well-known (cf. for example [Mochizuki2], Theorem 19.1). 
So the statement of this conjecture is equivalent to the
surjectivity of $\varphi_X$, i.e. that every group-theoretic section of the natural projection
$\Pi_X\twoheadrightarrow G_k$ is point-theoretic, under the above assumptions. 

Note that this conjecture can be formulated over any field,
but one can not expect it  to be true in general. For example, the analog of this conjecture doesn't hold over
finite fields. 

Indeed, over a finite field the natural
projection $\Pi_X \twoheadrightarrow G_k$ admits group-theoretic sections, since the profinite group
$G_k$ is free in this case. On the other hand there are proper, geometrically connected, hyperbolic, and
smooth curves over finite fields with no rational points. 

One can formulate an analog of the GASC over $p$-adic
local fields.

\definition {A $p$-adic Version of Grothendieck's Anabelian Section Conjecture ($p$-adic GASC)} Assume that $k$ is a $p$-adic local field,
i.e. a finite extension of the field $\Bbb Q_p$ for some prime integer $p$, and 
that $\Sigma=\Primes$. Then the map
$\varphi_X: X(k)\to  \overline {\Sec} _{\Pi_X}$ is bijective.
\enddefinition

The map $\varphi_X$ is also known to be injective in the case where $k$ is a $p$-adic local field, and $p\in \Sigma$ (cf. [Mochizuki2], Theorem 19.1).
Thus, the statement of the $p$-adic GASC is equivalent to the surjectivity of the map $\varphi_X$ in this case.

The following observation, du to Tamagawa, is crucial in investigating the Grothendieck anabelian section conjecture.

\proclaim{Lemma 3.1.2}\ Assume that $k$ is finitely generated over the prime field $\Bbb Q$, or that $k$ is a
$p$-adic local field. Let $s:G_k\to \Pi_X$ be a group-theoretic section of the natural projection
$\Pi_X \twoheadrightarrow G_k$, and $\{X_i[s]\}_{i\ge 1}$ a system of neighbourhoods of the section $s$ (cf. 1.3). Then
$s$ is point-theoretic if and only if $X_i[s](k)\neq \varnothing$, for every $i\ge 1$.
\endproclaim

\demo {Proof}\ See [Tamagawa], Proposition 2.8, (iv).
\qed
\enddemo

Lemma 3.1.2 reduces the proof of the Grothendieck anabelian section conjecture, in the case where
$k$ is finitely generated over the prime field $\Bbb Q$, or that $k$ is a
$p$-adic local field, and $\Sigma =\Primes$, to proving the following implication
$$\{\overline {\Sec} _{\Pi_X}\neq \emptyset\} \Longrightarrow \{X(k)\neq \emptyset\}.$$

\subhead {3.2}
\endsubhead
One can also formulate a birational version of the Grothendieck anabelian section conjecture as follows. See also [Pop].

Let 
$$\overline G_X\defeq \Gal (K_X^{\sep}/K_X.\bar k)^{\Sigma}$$ 
be the maximal pro-$\Sigma$ quotient of the absolute
Galois group $\Gal (K_X^{\sep}/K_X.\bar k)$, and
$$G_X\defeq \Gal (K_X^{\sep}/K_X)/\Ker (\Gal (K_X^{\sep}/K_X.\bar k)\twoheadrightarrow \Gal (K_X^{\sep}/K_X.\bar k)^{\Sigma})$$
the maximal geometrically pro-$\Sigma$ Galois group of the function field $K_X$. Thus, $G_X$ sits naturally
in the following exact sequence
$$1\to \overline G_X\to G_X\to G_k\to 1.$$

Let $x\in X(k)$ be a rational point. Then $x$ determines a decomposition subgroup
$D_x\subset G_X$, which is only defined up to conjugation by the elements of $\overline G_X$, and which maps surjectively
to $G_k$ via the natural projection $G_X\twoheadrightarrow G_k$. More precisely, $D_x$ sits naturally in the following exact sequence
$$1\to M_X\to D_x\to G_k\to 1.$$

The above group extension is known to be split. 
Indeed, the field extension of the completion ${\hat K}_{X,x}$ of $K_X$ at $x$ obtained by extracting $n$-th roots of a uniformising parameter
at $x$, for all $\Sigma$-integers $n$, determines a splitting of the above sequence.

The set of all possible splittings of this extension, i.e. sections
$G_k\to D_x$ of the natural projection $D_x\twoheadrightarrow G_k$, is a torsor under $H^1(G_k,M_X)$. Each section
$G_k\to D_x$ of the natural projection $D_x\twoheadrightarrow G_k$
determines naturally a section $G_k\to G_X$ of the natural projection $G_X\twoheadrightarrow G_k$, whose image is contained in $D_x$.

\definition {The Birational Grothendieck Anabelian Section Conjecture (BGASC)} 
\newline
Assume that $k$ is finitely generated over the
prime field $\Bbb Q$, and $\Sigma =\Primes$. Let $s:G_k\to G_X$ be a group-theoretic section of the natural projection
$G_X \twoheadrightarrow G_k$. Then the image $s(G_k)$ is contained in a decomposition subgroup $D_x$ associated to a unique
rational point $x\in X(k)$. In particular, the existence of the section $s$ implies that $X(k)\neq \varnothing$.
\enddefinition

One can also formulate an analog of the BGASC over $p$-adic local fields.

\definition {A $p$-adic Version of the Birational Grothendieck Anabelian Section Conjecture ($p$-adic BGASC)} 
Assume that $k$ is a $p$-adic local field, i.e. $k$ is a finite extension of $\Bbb Q_p$, for some prime integer $p$, and $\Sigma=\Primes$. 
Let $s:G_k\to G_X$ be a group-theoretic section of the natural projection
$G_X \twoheadrightarrow G_k$. Then the image $s(G_k)$ is contained in a decomposition subgroup $D_x$ associated to a unique
rational point $x\in X(k)$. In particular, the existence of the section $s$ implies that $X(k)\neq \varnothing$.
\enddefinition

\definition {Remark 3.2.1}\ An affirmative answer to the cuspidalisation problem for sections of 
arithmetic fundamental groups of hyperbolic curves (cf. 2.2), over a field which is finitely generated over the prime field $\Bbb Q$, 
plus the validity of the birational Grothendieck anabelian section conjecture (BGASC), implies (using lemma 3.1.2) 
the validity of the Grothendieck anabelian section conjecture (GASC) for $\pi_1$, in the case where $\Sigma=\Primes$.
\enddefinition

A similar Remark holds for the $p$-adic GASC.

\subhead {3.3}
\endsubhead 
The following result of Koenigsmann concerning the $p$-adic BGASC is fundamental (see [Koenigsmann]).

\proclaim {Theorem 3.3.1\ (Koenigsmann)} The $p$-adic BGASC holds true. More precisely,
assume that $k$ is a $p$-adic local field, and $\Sigma=\Primes$.
Let $s:G_k\to G_X$ be a group-theoretic section of the natural projection
$G_X \twoheadrightarrow G_k$. Then the image $s(G_k)$ is contained in a decomposition subgroup $D_x$ associated to a
unique rational point $x\in X(k)$. In particular, $X(k)\neq \varnothing$.
\endproclaim

This result has been strengthened by Pop who proved the following (see [Pop]). For a profinite group
$H$, and a prime integer $p$,  we denote by $H''$ the maximal $\Bbb Z/p\Bbb Z$-metabelian quotient of $H$.
Thus, $H''$ is the second quotient of the $\Bbb Z/p\Bbb Z$-derived series of $H$.

\proclaim {Theorem 3.3.2 \ (Pop)} Assume that $k$ is a $p$-adic local field, which contains a primitive $p$-th root of $1$, and assume $p\in \Sigma$.
Let $s:G_k''\to G_{X}''$ be a group-theoretic section of the natural projection
$G_{X}'' \twoheadrightarrow G_k''$. Then the image $s(G_k'')$ is contained in a decomposition subgroup $D_x\subset G_{X}''$
associated to a unique rational point $x\in X(k)$. In particular, $X(k)\neq \varnothing$. Here the $(\ \ )''$ of the various profinite groups
are with respect to the prime $p$, i.e. the second quotients of the $\Bbb Z/p\Bbb Z$-derived series.
\endproclaim

The above Theorem of Pop can be viewed as a very ``minimalistic'' version of the birational Grothendieck anabelian section
conjecture over $p$-adic local fields. Note that the quotient $G_k''$ of $G_k$ is finite, in the case where $k$ is a $p$-adic local field. 

Our main result concerning the Grothendieck anabelian section conjecture over $p$-adic local fields is the following.

\proclaim {Theorem 3.3.3} Assume that $k$ is a $p$-adic local field, which contains a primitive $p$-th root of unity,
and $p\in \Sigma$. Let $s:G_k\to \Pi_X$ be a
group-theoretic section of the natural projection  $\Pi_X \twoheadrightarrow G_k$. Then $s$ is a good section
(in the sense of Definition 1.4.1) if and only if there exists a section $s^{\c-\ab}:G_k\to G_{X}^{\c-\ab}$ of the
natural projection $ G_{X}^{\c-\ab}  \twoheadrightarrow G_k$ which lifts the section $s$.
Furthermore, if the section  $s^{\c-\ab}:G_k\to G_{X}^{\c-\ab}$ is tame point-theoretic, in the sense of Definition
1.7.1, then  $X(k)\neq \varnothing$.
\endproclaim

\demo{Proof}
The first assertion follows from Proposition 1.6.8, and  Theorem 2.7. 

For a profinite group $H$, denote by $H'$
the maximal quotient of $H$ which is abelian and annihilated by $p$. Thus, $H'$ is the first quotient of the
$\Bbb Z/p\Bbb Z$-derived series of $H$. 

The existence of the section  $s^{\c-\ab}:G_k\to G_{X}^{\c-\ab}$
implies the existence of a section $s':G_k'\to  G_X'$ of the natural projection $G_X'\twoheadrightarrow G_k'$. Let $\Tilde L/K_X$
be the sub-extension of $K_X^{\sep}/K_X$ with Galois group $G_X'$, and $L/K_X$
the sub-extension of $\Tilde L/K_X$ which corresponds to the closed subgroup  $s'(G_k')$ of $G_X'$. 

The assumption
that the section  $s^{\c-\ab}:G_k\to G_{X}^{\c-\ab}$ is tame point-theoretic, in the sense of Definition
1.7.1, implies that the natural homomorphism $\Br k\to \Br L$ is injective. Under this assumption (which is implied by the lifting property
of the section  $s':G_k'\to  G_X'$ to a section $s'':G_k''\to  G_X''$ which is imposed in [Pop]) Pop proves that the image
$s'(G_k')$ is contained in a decomposition subgroup $D_x\subset G_{K_X}'$ associated to a unique rational point $x\in X(k)$ (cf. loc. cit.).
In particular,  $X(k)\neq \varnothing$.
\qed
\enddemo

In light of Theorem 3.3.3, and in connection with the $p$-adic GASC, we conjecture the following.

\definition {Conjectures 3.3.4} Let $k$ be a $p$-adic local field, i.e. $k$ is a finite extension of $\Bbb Q_p$,
and $X$ a proper, smooth, geometrically connected, and hyperbolic curve over $k$.
\enddefinition
\definition {Conjecture A} Assume that $\Sigma=\Primes$. Let $s:G_k\to \Pi_X$ be a group-theoretic section of
the natural projection $\Pi_X\twoheadrightarrow G_k$. Then $s$ is good in the sense of Definition 1.4.1.
\enddefinition
\definition {Conjecture B} Assume that $\Sigma=\Primes$. 
Let $\tilde s:G_k\to G_X^{\c-\ab}$ be a group-theoretic section
of the natural projection $G_X^{\c-\ab}\twoheadrightarrow G_k$. Then $\tilde s$
is tame-point theoretic in the sense of Definition 1.7.1.
\enddefinition

As an immediate consequence of Theorem 3.3.3, and Lemma 3.1.2, one deduces the following, which links Conjecture A, and Conjecture B,
to the $p$-adic GASC.

\proclaim {Corollary 3.3.5} Assume that $k$ is a $p$-adic local field, and $\Sigma=\Primes$. 
Assume that Conjecture A, and Conjecture B, in 3.3.4, hold true for every proper, smooth, geometrically connected, and 
hyperbolic curve over $k$.
Then the $p$-adic GASC holds true.

In other words, we have the following implication
$$\{\Conjecture\ A+\Conjecture\ B\}\Longrightarrow \{p-\adic \ \GASC\}.$$
\endproclaim

In Theorem 3.3.3 we assumed that $k$ contains a primitive $p$-th root of unity. This condition is omitted in Corollary 3.3.5.
Indeed, one can easily verify that a section $s:G_k\to \Pi_X$ of the natural projection $\Pi_X\twoheadrightarrow G_k$ is point-theoretic
if and only if its restriction to an open subgroup of $G_k$ is point-theoretic, in the case where $\Sigma=\Primes$, by using a descent argument which resorts to 
the injectivity of the map $\varphi_X: X(k)\to  \overline {\Sec} _{\Pi_X}$.

\subhead {3.4}
\endsubhead
In [Esnault-Wittenberg1] sections of geometrically abelian absolute Galois groups 
of function fields of curves over number fields were investigated.

It is shown in loc. cit. that the existence of such a section implies (in fact is equivalent to) the existence of degree $1$ divisors on the curve, 
under a finiteness condition of the Tate-Shafarevich group of the jacobian of the curve.

As an application of our results on the cuspidalisation of sections of arithmetic fundamental groups, we can prove an analogous result
for good sections of arithmetic fundamental groups.

Our main result concerning the Grothendieck anabelian section conjecture over number fields is the following.

\proclaim {Theorem 3.4.1} Assume that $k$ is a number field, i.e. a finite extension of the prime field $\Bbb Q$,
and $\Sigma=\Primes$. 
Let $s:G_k\to \Pi_X$ be a group theoretic section of the natural projection  $\Pi_X \twoheadrightarrow G_k$. 
Assume that $s$ is a uniformly good group-theoretic section (in the sense of definition 1.4.1), and that the jacobian variety
of $X$ has a finite Tate-Shafarevich group. Then there exists a divisor of degree $1$ on $X$.
\endproclaim

\demo{Proof} Follows formally from Theorem 2.7, and Theorem 2.1 in [Esnault-Wittenberg1].
\qed
\enddemo

\subhead
\S 4. Sections of Geometrically pro-$\Sigma$ Arithmetic Fundamental Groups of Curves over
$p$-adic Local Field: $p\notin \Sigma$
\endsubhead
In $\S4$ we investigate sections of geometrically pro-$\Sigma$ arithmetic fundamental groups of hyperbolic algebraic
curves over $p$-adic local field, in the case where $p\notin \Sigma$. 

We give
examples of such sections which are not point-theoretic, and 
give another direct proof of the fact that such sections are good sections (in the sense of Definition 1.4.1), 
which doesn't resort to Proposition 1.5.2, and Proposition 1.6.3 (i).
We will use the notations in $\S1$.

\subhead {4.1}
\endsubhead
One of the difficulties in investigating the Grothendieck anabelian section conjecture is that, for the time being, one doesn't 
know how to construct sections of arithmetic fundamental groups, and hence test the 
validity of the conjecture on concrete examples.
One way to construct such sections is as follows. 

Let $X$ be a proper,
smooth, geometrically connected, and hyperbolic algebraic curve over a field $k$, 
and $\Sigma \subseteq \Primes$ a non-empty set of prime integers, with $\char (k)\notin \Sigma$.
Consider the exact sequence
$$1\to \Delta_X\to \Pi_X @>{\pr}>> G_k\to 1,\tag {$4.1$}$$
where $\Pi_X$ is the geometrically pro-$\Sigma$ arithmetic fundamental group of $X$. 

Note that the exact sequence $(4.1)$ induces a natural homomorphism
$$\rho _{X,\Sigma}: G_k\to \Out (\Delta _X),$$
where 
$$\Out (\Delta _X)\defeq \Aut (\Delta_X)/\Inn (\Delta_X)$$ 
is the group of outer automorphisms of $\Delta_X$.
For $g\in G_k$, its image $\rho _{X,\Sigma}(g)$ is the class of the automorphism of $\Delta_X$ obtained by
lifting $g$ to an element $\tilde g\in \Pi_X$, and letting $\tilde g$ act on $\Delta _X$ by inner conjugation.

If $G$ is a slim profinite group, then we have a natural exact sequence
$$1\to G\to \Aut G \to \Out G \to 1,$$
where the homomorphism $G\to \Aut G$ sends an element $g\in G$ to the corresponding inner
automorphism $h\mapsto ghg^{-1}$. 

Moreover, if the profinite group $G$ is finitely generated then the groups 
$\Aut(G)$, and $\Out(G)$, are naturally endowed with a profinite topology, and the above sequence is an exact 
sequence of profinite groups.

\proclaim {Lemma 4.1.1} The profinite group $\Delta_ X$ is slim. In particular, the exact sequence
$(4.1)$ is obtained form the following exact sequence
$$1\to \Delta _X\to \Aut (\Delta _X) \to \Out (\Delta _X)\to 1, \tag  {$4.2$}$$
by pull back via the natural continuous homomorphism $\rho _{X,\Sigma}:G_k\to \Out (\Delta _X)$. 

More precisely, we have a commutative diagram:
$$
\CD
1 @>>>  \Delta _X     @>>> \Aut (\Delta _X)     @>>> \Out (\Delta _X)  @>>> 1\\
  @.        @A{\id}AA           @AAA              @A{\rho _{X,\Sigma}}AA \\
1 @>>>  \Delta _X       @>>> \Pi_X       @>>> G_k   @>>> 1
\endCD
\tag {$4.3$}
$$
where the horizontal arrows are exact, and the right square is cartesian.
\endproclaim

\demo {Proof} Well known (cf. for example [Tamagawa], Proposition 1.11).
\qed
\enddemo

To construct a continuous group-theoretic section $s:G_k\to \Pi_X$ of the natural projection
$\Pi_X\twoheadrightarrow G_k$, it is equivalent to construct a continuous homomorphism 
$${\tilde \rho_{X,\Sigma}}: G_k\to \Aut (\Delta _X),$$ 
which lifts the continuous homomorphism  ${\rho_{X,\Sigma}}:G_k\to \Out (\Delta _X)$ above, 
i.e. such that the following diagram commutes
$$
\CD
G_k@  >{\tilde \rho_{X,\Sigma}}>> \Aut (\Delta _X) \\
@V{\id}VV         @VVV \\
G_k  @>{\rho_{X,\Sigma}}>>  \Out (\Delta _X) \\
\endCD
$$
as follows directly from the fact that the right square in the above diagram (4.3) is cartesian.

\subhead {4.2}
\endsubhead
Let $p>0$ be a prime integer. For the rest of this section we will assume that
the field $k$ is a $p$-adic local field, i.e. $k$ is a finite extension of $\Bbb Q_p$. 
Write $O_k$ for the ring of integers of $k$, and
$F$ for the residue field of $k$, which is a finite field.

One can not expect the $p$-adic version of the Grothendieck anabelian section conjecture (cf. 3.1, $p$-adic GASC) to hold,
if in the statement of the conjecture one considers a set of prime integers $\Sigma$ not containing $p$.

Indeed, first of all if $p\notin \Sigma$ the map
$\varphi_X\defeq \varphi_{X,\Sigma}: X(k)\to  \overline {\Sec} _{\Pi_X}$ is not injective in general, 
and one can prove that it may not be surjective. More precisely, we have the following.

\proclaim{Proposition 4.2.1}  Let $k$ be a $p$-adic local field, and $\Sigma$ a non-empty set of prime integers, with $p\notin
\Sigma$. Then there exists a smooth, proper, and geometrically connected  hyperbolic curve $X$ over $k$, which has good reduction
over $\Cal O_k$, and such that the natural map
$\varphi_X\defeq \varphi_{X,\Sigma}: X(k)\to  \overline {\Sec} _{\Pi_X}$ (cf. 3.1) is not surjective.
\endproclaim

Before proving Proposition 4.2.1, let's assume that the hyperbolic $k$-curve $X$ has good reduction over $O_k$, i.e. $X$ 
extends to a smooth, and proper,
relative curve $\Cal X$ over $O_k$. Let $\Cal X_s\defeq \Cal X\times _{O_k}F$ be the special fiber of $\Cal X$.

Let $\xi$ be a geometric point of $\Cal X_s$ above the generic point of $\Cal X_s$.
Then $\xi$ determines naturally an algebraic closure $\overline F$
of $F$, and a geometric point $\bar {\xi}$ of $\overline {\Cal X_s} \defeq \Cal X_s\times _F \overline F$.

There exists a canonical exact sequence of profinite groups
$$1\to \pi_1(\overline {\Cal X_s},\bar \xi)\to \pi_1(\Cal X_s, \xi) @>>> G_F\to 1.$$

Here, $\pi_1(\Cal X_s, \xi)$ denotes the arithmetic \'etale fundamental group of $\Cal X_s$ with base
point $\xi$, $\pi_1(\overline {\Cal X_s},\bar \xi)$ the \'etale fundamental group of $\overline {\Cal X_s}\defeq 
\Cal X_s\times _F \overline F$ with base
point $\bar \xi$, and $G_F\defeq \Gal (\overline F/F)$ the absolute Galois group of $F$.

Write 
$$\Delta_{\Cal X_s}\defeq \pi_1(\overline {\Cal X_s},\bar \xi)^{\Sigma}$$
for the maximal pro-$\Sigma$ quotient of $\pi_1(\overline X,\bar \xi)$,
and
$$\Pi_{\Cal X_s}\defeq  \pi_1(\Cal X_s, \xi)/ \Ker  (\pi_1(\overline {\Cal X_s},\bar \xi)\twoheadrightarrow
\pi_1(\overline {\Cal X_s},\bar \xi)^{\Sigma})$$
for the quotient of  $\pi_1(\Cal X_s, \xi)$ by the kernel of the natural surjective homomorphism
$\pi_1(\overline {\Cal X_s},\bar \xi)\twoheadrightarrow \pi_1(\overline {\Cal X},\bar \xi)^{\Sigma}$, 
which is a normal subgroup
of $\pi _1(\Cal X,\xi)$. Thus, we have an exact sequence of profinite groups
$$1\to \Delta_{\Cal X_s}\to \Pi_{\Cal X_s} @>>> G_F\to 1.\tag {$4.4$}$$

Moreover, after a suitable choice of the base points $\xi$ and $\eta$, there exists a natural commutative specialisation
diagram:

$$
\CD
1 @>>>  \Delta _X     @>>> \Pi_X     @>>> G_k @>>> 1\\
  @.        @VVV           @V{\Sp}VV              @VVV \\
1 @>>>    \Delta _ {\Cal X_s}     @>>> \Pi_{\Cal X_s}      @>>> G_F   @>>> 1
\endCD
\tag {$4.5$}
$$
where the left vertical homomorphism $\Sp:\Delta _X \to \Delta _ {\Cal X_s}$ is an 
isomorphism, since we assumed $p\notin \Sigma$,
and the right vertical
homomorphisms are surjective, as follows easily from the specialisation theory for fundamental groups
of Grothendieck (cf. [Grothendieck1]). 

In fact one has the following more precise statement.

\proclaim {Lemma 4.2.2} In the above commutative diagram (4.5) the right square is cartesian.
\endproclaim

\demo {Proof}\ The slimness of $\Delta _X$ implies that we have the following commutative diagram:
$$
\CD
1 @>>>  \Delta _X     @>>> \Aut (\Delta _X)     @>>> \Out (\Delta _X)  @>>> 1\\
  @.        @A{\id}AA           @AAA              @A{\rho _{X,\Sigma}}AA \\
1 @>>>  \Delta _X       @>>> \Pi_X       @>>> G_k   @>>> 1
\endCD
$$
where the horizontal arrows are exact, and the right square is cartesian. 

The right vertical homomorphism
${\rho _{X,\Sigma}}:G_k\to \Out (\Delta _X)$ factors as $G_k\twoheadrightarrow G_F \to \Out (\Delta _X)$, where
$G_k\twoheadrightarrow G_F$ is the natural projection, and the homomorphism $G_F \to \Out (\Delta _X)$
is naturally deduced from the exact sequence (4.4) (i.e. equals $\rho _{\Cal X_s,\Sigma}$ followed by the natural identification
$\Out (\Delta _X)\isom \Out (\Delta _{\Cal X_k})$ induced by the specialisation isomorphism
$\Sp:\Delta _X \to \Delta _ {\Cal X_s}$, since $X$ has good
reduction over $O_k$, and $p\notin \Sigma$, as is well-known). 

Thus, we have the following commutative diagram:
$$
\CD
1 @>>>  \Delta _X     @>>> \Aut (\Delta _X)     @>>> \Out (\Delta _X)  @>>> 1\\
@.        @A{\id}AA           @AAA              @AAA \\ 
1@>>>   \Delta _X   @>>>      \Tilde {\Pi}_X   @>>>          G_F @>>> 1 \\
@.        @A{\id}AA           @AAA              @AAA \\
1 @>>>  \Delta _X       @>>> \Pi_X       @>>> G_k   @>>> 1
\endCD
$$
where the horizontal arrows are exact, and the right squares are cartesian. 

Here, $\Tilde {\Pi}_X$ is by definition
the pull back of the upper exact sequence by the above natural homomorphism $G_F\to \Out (\Delta _X)$. Note that the specialisation
isomorphism $\Sp:\Delta _X \to \Delta _ {\Cal X_s}$ induces a natural isomorphism
$\Tilde {\Pi}_X\isom \Pi_{\Cal X_s}$. From this follows our claim, since the lower right square in the above commutative
diagram is cartesian.
\qed
\enddemo

\demo{Proof of Proposition 4.2.1}
First, there exists a proper, smooth, and geometrically connected hyperbolic curve $X$ over $k$, which has good reduction over $O_k$,
and such that $\Cal X_s (F)=\varnothing$, where $\Cal X_s$ is the special fibre of a proper, and smooth, model $\Cal X$ of $X$.

Indeed, there exists a proper, smooth, hyperbolic, and geometrically connected curve  $\Cal X_s$ over $F$, with
$\Cal X_s (F)=\varnothing$, and $\Cal X_s$ can be lifted to a proper, smooth, and
geometrically connected curve over $\Cal O_k$, with generic fiber $X$ over $k$.  

Recall the commutative diagram:
$$
\CD
1 @>>>  \Delta _X     @>>> \Pi_X     @>>> G_k @>>> 1\\
  @.        @V{\Sp}VV           @VVV              @VVV \\
1 @>>>    \Delta _ {\Cal X_s}     @>>> \Pi_{\Cal X_s}      @>>> G_F   @>>> 1
\endCD
\tag {$4.5$}
$$

The natural projection $\Pi_{\Cal X_s}\twoheadrightarrow G_F$ admits sections, since the Galois group $G_F$ is pro-free.
Let $\bar s:G_F\to \Pi_{\Cal X_s}$ be such a section. Note that $\bar s$ can not be point-theoretic, since $\Cal X_s (F)
=\varnothing$. 

The section $\bar s$ can be lifted to a section $s:G_k\to \Pi_X$ of the natural projection 
$\Pi_X\twoheadrightarrow G_k$, since the right square in the above diagram is cartesian (cf. Lemma 4.2.2). 
Thus, we have a commutative diagram:

$$
\CD
G_k @>s>>   \Pi_X \\
@VVV          @V{\Sp}VV \\
G_F   @>{\bar s}>>  \Pi_{\Cal X_s}\\
\endCD
$$

Finally, the section $s$ can not be point-theoretic, for otherwise the section $\bar s$ would be point-theoretic, as is easily verified.
\qed
\enddemo

\subhead {4.3}
\endsubhead
Although group-theoretic sections of geometrically pro-$\Sigma$ arithmetic fundamental groups of curves over p-adic local fields
may not be point-theoretic if $p\notin \Sigma$, they are uniformly good in the sense of definition 1.4.1, as 
we already proved in Proposition 1.5.2, and Proposition 1.6.3, (i).  Next, we give another direct proof of this fact.

\proclaim {Proposition 4.3.1 (Uniform Goodness of Sections of Geometrically pro-$\Sigma$ Arithmetic Fundamental Groups of Curves over
$p$-adic Local Field: $p\notin \Sigma$)} Let $k$ be a $p$-adic local field, and $\Sigma$ a non-empty set of prime integers
with $p\notin \Sigma$. Let $X$ be a smooth, proper, geometrically connected, and hyperbolic curve $X$ over $k$. Let $\Pi_X$ be the
geometrically pro-$\Sigma$ arithmetic fundamental group of $X$, and $s:G_k\to \Pi_X$ a group-theoretic section of the natural
projection $\Pi_X\twoheadrightarrow G_k$. Then $s$ is uniformly good in the sense of Definition 1.4.1.
\endproclaim

\demo {Proof} First, uniform goodness is equivalent to goodness (cf. Proposition 1.6.8). So it is enough to prove that
$s$ is a good section. 

Let $s^{\star}:H^2(X,M_X)\isom H^2(\Pi_X,M_X)\to H^2(G_k,M_X)$ be the natural restriction homomorphism, which is induced by the section $s$.
Let $\{X_i[s]\}_{i\ge 1}$ be a system of neighbourhoods of the section $s$, and $\Pi_X[i,s]$ the geometrically pro-$\Sigma$ arithmetic
fundamental group of $X_i[s]$. 

For every positive integer $i$, let $s_i:G_k\to \Pi_X[i,s]$ be the induced group-theoretic section,
and (bearing in mind the natural identifications  $H^2(\Pi_X[i,s],M_X)\isom H^2(X_i[s],M_X)$
(cf. [Mochizuki], Proposition 1.1)) $s_i^{\star}: H^2(X_i[s],M_X)\to H^2(G_k,M_X)$ the natural (restriction) homomorphism.
We will show that $s_i^{\star}$ annihilates the Picard part $\Pic(X_i)^{\wedge, \Sigma}$ of  $H^2(X_i[s],M_X)$.

We have a natural isomorphism $H^2(G_k,M_X)\isom \hat \Bbb Z^{\Sigma}$. In particular, in order to show that
$s_i^{\star}$ annihilates $\Pic(X_i)^{\wedge, \Sigma}$, it suffices to show that the image
$s_i^{\star} (\Pic(X_i)^{\wedge, \Sigma})$
of $\Pic(X_i)^{\wedge, \Sigma}$ in $H^2(G_k,M_X)$ is torsion. For this it suffices to show that the image of the degree zero
part of $\Pic(X_i)^{\wedge, \Sigma}$ in $H^2(G_k,M_X)$ is torsion. 

Indeed, for each positive integer $i$, we have a
natural exact sequence:
$$0 @>>> \Pic^0(X_i)^{\wedge,\Sigma}  @>>> \Pic(X_i)^{\wedge,\Sigma} @>{\deg}>> \hat \Bbb Z^{\Sigma}$$
where $\Pic^0(X_i)^{\wedge,\Sigma}$ is the $\Sigma$-adic completion of the degree $0$ part $\Pic^0(X_i)$ of $\Pic(X_i)$,
and the right map is induced by the degree homomorphism $\deg:\Pic(X_i)\to \Bbb Z$. Passing to the direct limit we obtain a
natural exact sequence
$$0 @>>> \underset{i\ge 1}\to{\varinjlim}\Pic^0(X_i)^{\wedge,\Sigma}  @>>>  \underset{i\ge 1}\to{\varinjlim}\Pic(X_i)^{\wedge,\Sigma}
@>{\deg}>>  \underset{i\ge 1}\to{\varinjlim} \ \hat \Bbb Z^{\Sigma}.$$
By passing to the direct limit, we also obtain a natural restriction homomorphism
$$\tilde s\defeq \underset{i\ge 1}\to{\varinjlim}s_i^{\star}:                 
\underset{i\ge 1}\to{\varinjlim} H^2(X_i[s],M_X)\to H^2(G_k,M_X).$$

If the homomorphism $\tilde s$ annihilates $\underset{i\ge 1}\to{\varinjlim}\Pic^0(X_i)^{\wedge,\Sigma}$, then it will induce
a natural homomorphism 
$\tilde s:\Im (\deg)\to  H^2(G_k,M_X)$, 
where $\Im (\deg)\subseteq \underset{i\ge 1}\to{\varinjlim} \ \hat \Bbb Z^{\Sigma}$ is the image of the above degree map.
The homomorphism $\Tilde s$ is then necessarily $0$,
since the various maps $\hat \Bbb Z^{\Sigma}\to \hat \Bbb Z^{\Sigma}$ in the inductive limit
$\underset{i\ge 1}\to{\varinjlim} \ \hat \Bbb Z^{\Sigma}$
are multiplication by the various degrees $\vert X_{i+1}:X_i\vert$
of the finite \'etale covers $X_{i+1}\to X_i$, and $\tilde s (\Im (\deg))$ 
would be divisible.

Thus, it suffices to show that $\tilde s(\underset{i\ge 1}
\to {\varinjlim}\Pic^0 (X_i[s])^{\wedge,\Sigma})=0$. 

But $\tilde s(\underset{i\ge 1}
\to {\varinjlim}\Pic^0 (X_i[s])^{\wedge,\Sigma})$ is torsion, since $\Pic^0 (X_i[s])^{\wedge,\Sigma}$ is torsion
for every $i$, as follows from the well-known structure of $J_{X_i[s]}(k)$, where $J_{X_i[s]}$ is the jacobian of $X_i[s]$, and the fact that
$p\notin \Sigma$.

This finishes the proof of Proposition 4.3.1.
\qed
\enddemo

\definition{Remark 4.3.2} Assume that $k$ is a $p$-adic local field. The condition $p\in \Sigma$ does not guarantee the point-theorecity of a group-theoretic section
$s:G_k\to \Pi_X$ of the natural projection $\Pi_X\twoheadrightarrow G_k$.
Indeed, Hoshi has recently constructed non geometric sections $s:G_k\to \Pi_X$ in the case where $\Sigma =\{p\}$ (cf. [Hoshi]).
However, it is still not known weather every group-theoretic section is point-theoretic, in the case where $\Sigma =\Primes$.
\enddefinition

\subhead
\S 5. Appendix
\endsubhead
In this appendix we give an alternative definition of the natural pairing between the Brauer-Grothendieck
group and the Picard group of a smooth, proper, and geometrically connected algebraic curve over a field.

Let $l$ be a prime integer.
Let $k$ be a field of characteristic $l\ge 0$, and $X$ a smooth, proper, and geometrically connected curve over $k$.

Let $\Br(X)\defeq H^2_{\et}(X,\Bbb G_m)$ be the Brauer-Grothendieck group of $X$, and $\Pic(X)\defeq
H^1_{\et}(X,\Bbb G_m)$ the Picard group of $X$. 
Let $\Sigma \subset \Primes$ a non-empty subset of the set $\Primes$ of all prime integers,
with $l\notin \Sigma$.

The evaluation map of the Brauer classes in  $\Br(X)$ at the closed points of
$X$, followed by the corestriction map in cohomology, induces a natural pairing (cf. [Lichtenbaum], 3)
$$<,>:\Pic(X)\times \Br(X)\to \Br(k).$$

More precisely, for a closed point $x\in X$, and a Brauer class $\alpha\in \Br (X)$, the value $<\Cal O(x),\alpha>$,
where $\Cal O(x)$ is the line bundle associated to $x$, is the image
of $\alpha$ via the natural map $\Br (X)\to \Br(k(x))@>\cor>> \Br (k)$.

This pairing induces for every positive $\Sigma$-integer $n$, meaning that $n$ is an integer which is only divisible by
primes in $\Sigma$, a natural pairing
$$<,>_n:\Pic(X)/n\Pic(X)\times _n \Br(X)\to _n \Br(k),$$
where $_n \Br(\ )$ denotes the part of $\Br (\ )$ which is annihilated by $n$. Moreover, for positive $\Sigma$-integers
$n$ and $m$, with $n$ divides $m$, the above pairing inserts into a commutative diagram:

$$
\CD
 \Pic(X)/m\Pic(X)\times _m \Br(X)   @> {<,>_m}>>   _m \Br(k)\\
                  @VVV     @VVV  \\
\Pic(X)/n\Pic(X)\times _n \Br(X)   @>{<,>_n}>>   _n \Br(k)\\
\endCD
$$
where the vertical maps are the natural homomorphisms. 

In particular, by passing to the
projective limit we obtain a natural pairing
$$<,>: \underset{n\ \Sigma-\text {integer}}\to{\varprojlim}  \Pic (X)/n \Pic (X)\times
\underset{n\ \Sigma-\text {integer}}\to{\varprojlim} _n\Br (X)\to
\underset{n\ \Sigma-\text {integer}}\to{\varprojlim} _n\Br (k).$$

We will denote by $\Pic(X)^{\wedge, \Sigma}\defeq \underset{n\ \Sigma-\text {integer}}\to{\varprojlim} 
\Pic (X)/n \Pic (X)$ the $\Sigma$-adic completion of the Picard group $\Pic(X)$, and
$T _{\Sigma}\Br (X)$ (resp. $T _{\Sigma}\Br (k)=H^2(G_k,M_X)$)
the $\Sigma$-Tate module of the Brauer group $\Br (X)$ (resp. the $\Sigma$-Tate module of $\Br (k)$).

Next, we will give an alternative definition of the above pairing 
$$<,>_n:\Pic(X)/n\Pic(X)\times _n \Br(X)\to _n \Br(k).$$

For a positive $\Sigma$-integer $n$, the Kummer exact sequence in \'etale topology
$$1\to \mu_n\to \Bbb G_m @>n>> \Bbb G_m\to 1$$
induces naturally an exact sequence of abelian groups
$$0 \to \Pic (X)/n \Pic (X) \to  H^2(X,\mu_{n}) \to _{n} \Br (X)
\to 0.$$
We shall refer to the subgroup $\Pic (X)/n \Pic (X)$ of $H^2(X,\mu_{n})$ as the Picard part of
the \'etale cohomology group $H^2(X,\mu_{n})$.

We will first define a natural pairing
$$(.)_n:\Pic(X)/n\Pic(X)\times H^2(X,\mu_n)\to _n \Br(k),$$
$$([\Cal L],\alpha)\mapsto ([\Cal L].\alpha)_n,$$
which will induce a natural pairing $\Pic(X)/n\Pic(X)\times _n \Br(X)\to _n \Br(k)$.

Let $\alpha \in  H^2(X,\mu_n)$, and $[\Cal L]\in  \Pic(X)/n\Pic (X)$ a Picard element.
We will define the element $([\Cal L].\alpha)_n \in  H^2(G_k,\mu_n)$. Clearly it suffices to define
the element $([\Cal L].\alpha)_n$ in the special case where a representative $\Cal L\in \Pic(X)$ of
$[\Cal L]$ is of the form
$\Cal L\defeq \Cal O(x)$, i.e.  is the line bundle associated to a closed point $x\in X$, since the line bundles of the form $\Cal O(x)$ 
 generate $\Pic(X)$.

Let $1\to \mu_n \to E \to \Pi_X\to 1$ be a group extension whose class in
$H^2(\Pi_X,\mu_n)$ coincides, via the natural identification $H^2(\Pi_X,\mu_n) \isom H^2(X,\mu_n)$, with $\alpha$.
Let $k(x)$ be the residue field of $X$ at $x$, and $D_x\subset \Pi_X$ a decomposition subgroup of $\Pi_X$
which is associated to $x$ ($D_x$ is only defined up to conjugation). Thus, $D_x$ maps isomorphically to the open
subgroup $G_{k(x)}$
of $G_k$ via the natural projection $\Pi_X \twoheadrightarrow G_k$. 

By pulling back the group extension $E$ by the
natural inclusion $D_x\hookrightarrow \Pi_X$, and using the natural identification $D_x\isom G_{k(x)}$ arising from the natural projection $\Pi_X\twoheadrightarrow G_k$,
we obtain a group extension $1\to \mu_n \to \Tilde E_x \to G_{k(x)}\to 1$
whose class in $H^2(G_{k(x)},\mu _n)$ maps, via the corestriction homomorphism
$\cor:H^2(G_{k(x)},\mu _n)\to H^2(G_k,\mu _n)$, to an extension class in $H^2(G_k,\mu _n)$ which is the desired
class $([\Cal L].\alpha)_n$.

\proclaim {Lemma A.1} The class $([\Cal L].\alpha)_n\in  H^2(G_k,\mu _n)$ defined above is well defined.
Hence we obtain a natural pairing
$$(.)_n:\Pic(X)/n\Pic(X)\times H^2(X,\mu_n)\to _n \Br(k),$$
$$([\Cal L],\alpha)\mapsto ([\Cal L].\alpha)_n.$$
\endproclaim

\demo{Proof}
Indeed, One easily verifies that this definition doesn't depend on the choice of the group extension $E$ whose class equals
$\alpha$, as well as the choice of the representative $\Cal L$ modulo $n \Pic(X)$ of the class $[\Cal L]$.
\qed
\enddemo

\proclaim {Lemma A.2} The above pairing
$$(.)_n:\Pic(X)/n\Pic(X)\times H^2(X,\mu_n)\to _n \Br(k),$$
induces a natural pairing
$$<,>_n\Pic(X)/n\Pic(X)\times _n \Br(X)\to _n \Br(k).$$
\endproclaim

\demo{Proof}
We will show that the restriction 
$$(.)_n:\Pic(X)/n\Pic(X)\times \Pic(X)/n\Pic(X)\to _n \Br(k)$$ 
of the above pairing
$(.)_n:\Pic(X)/n\Pic(X)\times H^2(X,\mu_n)\to _n \Br(k)$ to the subgroup $\Pic(X)/n\Pic(X)\subseteq H^2(X,\mu_n)$ is trivial.

Suppose that $\Cal L=\Cal O(x)$, and the extension class $\alpha=[\Cal L']$ is a Picard element. Then as above we can assume, without loss of
generality, that $\Cal L'=\Cal O(x')$ is the line bundle associated to a closed point $x'$. Also we can assume, without loss of generality,
that both $x$, and $x'$, are $k$-rational points. We will show that $([\Cal L],[\Cal L'])_n=0$.

Let $U_{x'}\defeq X\setminus \{x'\}$, and $\Pi_{U_{x'}}^{\c-\cn}$
the maximal (geometrically) cuspidally
central quotient of $\Pi_{U_{x'}}$, with respect to the natural projection  $\Pi_{U_{x'}}\twoheadrightarrow \Pi_X$
(cf. 2.1).
Then we have a natural exact sequence $1\to M_X\to \Pi_{U_{x'}}^{\c-\cn} \to \Pi_X\to 1$,  whose extension class in
 $H^2(\Pi_X,M_X)$ coincides, via the natural identification $H^2(\Pi_X,M_X) \isom H^2(X,M_X)$, with the Chern class of
the line bundle $\Cal L'$ (cf. [Mochizuki3], Lemma 4.2). 

The pull back $1\to M_X\to E_x \to G_k\to 1$
of the group extension  $\Pi_{U_{x'}}^{\c-\cn}$
via the natural inclusion $D_x\hookrightarrow \Pi_X$ is a group extension
whose class in $H^2(G_k,M_X)$ is trivial. Indeed, any decomposition group $D_x \subset  \Pi_{U_{x'}}^{\c-\cn}$
which is associated to $x$ defines a splitting of the group extension $E_x$. Furthermore,
the push forward $1\to M_X/nM_X\to (E_x)_n \to G_k\to 1$
of the group extension $E_x$ via the natural projection $M_X \twoheadrightarrow M_X/nM_X$ gives rise to an extension class
in $H^2(G_k,\mu_n)$, which is trivial, and which coincides by definition with $([\Cal L].\alpha)_n$.

From this we deduce that $([\Cal L].\alpha)_n=0$.

Thus, the natural pairing 
$$(.)_n:\Pic(X)/n\Pic(X)\times H^2(X,\mu_n)\to _n \Br(k)$$
induces a natural pairing
$$<,>_n:\Pic(X)/n\Pic(X)\times _n \Br(X)\to _n \Br(k),$$
$$([\Cal L],\beta)\mapsto <[\Cal L],\beta>\defeq ([\Cal L].\alpha),$$
where $\alpha\in  H^2(X,\mu_n)$ is an element in the pre-image of $\beta\in _{n}\Br (X)$.
\qed
\enddemo

\proclaim {Lemma A.3} Let $m$ and $n$ be $\Sigma$ integers, with $n$ divides $m$. Then we have a commutative diagram:
$$
\CD
 \Pic(X)/m\Pic(X)\times _m \Br(X)   @>{<,>_m}>>   _m \Br(k)\\
                  @VVV     @VVV  \\
\Pic(X)/n\Pic(X)\times _n \Br(X)   @>{<,>_n}>>   _n \Br(k)\\
\endCD
$$
where the horizontal maps are the above defined pairings, and the vertical maps are the natural homomorphisms.
In particular, by passing to the projective limit we obtain a natural pairing
$$<,>:\Pic(X)^{\wedge, \Sigma}\times  T_{\Sigma} \Br (X)\to T _{\Sigma}\Br (k).$$
\endproclaim

\demo{Proof}
Clear.
\qed
\enddemo

Let $\Pi_X$ be the geometrically pro-$\Sigma$ arithmetic fundamental group of $X$, and 
$s:G_k\to \Pi_X$ a continuous group-theoretic section of the natural projection $\Pi_X\twoheadrightarrow G_k$.

By pulling back cohomology classes via the section $s$, and bearing in mind the natural identifications 
$H^2(\Pi_X,M_X)\isom H^2(X,M_X)$ (cf. [Mochizuki], Proposition 1.1), we obtain a natural restriction
homomorphism
$$s^{\star}: H^2(X,M_X)\to H^2(G_k,M_X).$$

Recall that we have a natural exact sequence
$$0\to\Pic (X)^{\wedge,\Sigma}\to H^2(X,M_X)\to T_{\Sigma}\Br(X)\to 0,$$
where the subgroup $\Pic (X)^{\wedge,\Sigma}$ is the Picard part of $H^2(X,M_X)$ (cf. exact sequence (1.6)).

\proclaim {Lemma A.4} The above natural pairing $\Pic(X)^{\wedge, \Sigma}\times  T_{\Sigma} \Br (X)
\to T_{\Sigma} \Br (k)$ induces a natural homomorphism
$$\delta:\Pic(X)^{\wedge, \Sigma} \to T_{\Sigma} \Br (k),$$
$$L\mapsto \delta(L) \defeq <L,[c(s)]>,$$
where $<L,[c(s)]>$ is the result of pairing $L\in \Pic(X)^{\wedge, \Sigma}$ with the image $[c(s)] \in  T_{\Sigma}\Br(X)$
of the Chern class $c(s)\in H^2(X,M_X)$ which is associated to the section $s$ (cf. Definition 1.3.2),
via the natural homomorphism $H^2(X,M_X)\to T_{\Sigma}\Br(X)$.

Then the two homomorphisms $\delta$, and $s^{\star}$, are identical on  $\Pic(X)^{\wedge, \Sigma}$.
\endproclaim

\demo{Proof} In a similar way as above, we can define for every $\Sigma$-integer $n$ natural homomorphisms
$s^{\star}_n: H^2(X,\mu_n)\to H^2(G_k,\mu_n)$, which are the restriction maps in cohomology induced by the section
$s:G_k\to \Pi_X$, and $\delta_n:\Pic(X)/n\Pic(X) \to _n\Br (k)$ defined by $\delta_n([\Cal L])=<[\Cal L],c(s)_n>$,
where $c(s)_n$ is the image of the Chern class $c(s)$ in  $H^2(X,\mu_n)$, and it suffices to show that
$s^{\star}_n$ and $\delta_n$ are identical on $\Pic(X)/n\Pic (X)$. 

For this it suffices, without loss of generality, to consider the case where a representative $\Cal L$ of $[\Cal L]$ is of the form
$\Cal L=\Cal O(x)$, where $x\in X(k)$, i.e. is the line bundle associated to a closed $k$-rational point $x\in X(k)$. 

Let
$$1\to \mu_n \to \Cal D \to \Pi_{X\times X}\to 1$$
be a group extension whose class in  $H^2(X\times X,\mu_n)$ coincides with the Chern class $c(\eta)_n$
of the diagonal embedding $\iota:X\to X\times X$. 

Then it follows from the various definitions that
$\delta_n([\Cal L])\defeq <[\Cal L],[c(s)]>_n=s^{\star}_n([\Cal L])$, and coincide with the extension class
of the group extension $1\to \mu_n\to \Cal D_{s,x}\to G_k\to 1$ which is obtained by pulling back
the above group extension $\Cal D$ by the homomorphism $(s,i_x):G_k\times_{G_k} D_s\to  \Pi_{X\times X}$,
where $i_x:D_x \hookrightarrow \Pi_X$ is the natural inclusion.
\qed
\enddemo

$$\text{References.}$$

\noindent
[Borovoi-Colliot-Th\'el\`ene-Skorobogatov] Borovoi, M., Colliot-Th\'el\`enene, J.-L., Skorobogatov, A. N.,
The elementary obstruction and homogeneous spaces,  Duke Math. J.  141  (2008),  no. 2, 321--364.

\noindent
[Esnault-Wittenberg] Esnault, H., Wittenberg, O., Remarks on cycle classes of sections of the fundamental
group, Mosc. Math. J. 9 (2009), no. 3, 451-467.

\noindent
[Esnault-Wittenberg1] Esnault, H., Wittenberg, O., On abelian birational sections, Journal of the American Mathematical society, 
Volume 23, Number 3, July 2010, Pages 713-724.

\noindent
[Grothendieck] Grothendieck, A., Brief an G. Faltings, (German), with an
english translation on pp. 285-293,
London Math. Soc. Lecture Note Ser., 242, Geometric Galois actions, 1,
49-58, Cambridge Univ. Press,
Cambridge, 1997.

\noindent
[Grothendieck1] Grothendieck, A., Rev\^etements \'etales et groupe fondamental, Lecture 
Notes in Math. 224, Springer, Heidelberg, 1971.

\noindent
[Harari-Stix] Harari, D., Stix, J., Descent obstructions and fundamental exact sequence, arxiv 1005.1302.

\noindent
[Harari-Szamuely] Harari, D., Szamuley. T., Galois sections for abelianized fundamental groups (With an appendix by E. V. Flynn), 
Math. Ann. 344 (2009), no. 4, 779--800.

\noindent
[Hoshi] Hoshi, Y., Existence of nongeometric pro-$p$ Galois sections of hyperbolic curves, to appear in Publications of RIMS.

\noindent
[Koenigsmann] Koenigsmann, J., On the section conjecture in anabelian geometry,  J. Reine Angew. Math.  588 
(2005), 221--235.

\noindent
[Lichtenbaum] Lichtenbaum, S., Duality theorems for curves over $p$-adic fields,  Invent. Math.  7  (1969), 120--136.

\noindent
[Mochizuki] Mochizuki, S., Absolute anabelian cuspidalizations of proper hyperbolic curves,  J. Math. Kyoto
Univ.  47  (2007),  no. 3, 451--539.

\noindent
[Mochizuki1] Mochizuki, S., Topics surrounding the anabelian geometry of hyperbolic curves, 
Galois groups and fundamental groups,  119--165, Math. Sci. Res. Inst. Publ., 41, Cambridge Univ. Press,
Cambridge, 2003.

\noindent
[Mochizuki2] Mochizuki, S., The local pro-$p$ anabelian geometry of curves, Invent. Math.  138  (1999), 
no. 2, 319--423.

\noindent
[Mochizuki3] Mochizuki, S., Galois sections in absolute anabelian geometry,  Nagoya Math. J. 179 (2005), 17-45.

\noindent
[Pop] Pop, F., On the birational $p$-adic section Conjecture, Compos. Math. 146 (2010), no. 3, 621-637.

\noindent
[Stix] Stix, J., On the period-index problem in light of the section conjecture, Amer. J. Math. 132 (2010), no. 1, 157-180.

\noindent
[Tamagawa] Tamagawa, A., The Grothendieck conjecture for affine curves,  Compositio Math.  109  (1997),  no. 2, 135--194.

\noindent
[Wittenberg] Wittenberg, O., On Albanese torsors and the elementary obstruction,  Math. Ann.  340  (2008), 
no. 4, 805--838.
\bigskip

\noindent
Mohamed Sa\"\i di

\noindent
College of Engineering, Mathematics, and Physical Sciences

\noindent
University of Exeter

\noindent
Harrison Building

\noindent
North Park Road

\noindent
EXETER EX4 4QF 

\noindent
United Kingdom

\noindent
M.Saidi\@exeter.ac.uk

\end
\enddocument